\documentclass[final,3p]{elsarticle}
\usepackage{graphicx}
\usepackage{amssymb}
\usepackage{tikz}
\usetikzlibrary{positioning}
\usepackage[ruled,vlined]{algorithm2e}
\usepackage{amsthm}
\usepackage{amsfonts}
\usepackage{amsmath}
\usepackage{mathrsfs}
\usepackage{BOONDOX-cal}
\usepackage{xfrac}
\usepackage{relsize}
\usepackage{float}
\usepackage{fixltx2e}
\usepackage[justification=centering]{caption}
\usepackage{hyphenat}
\usepackage{epsfig}
\usepackage{epstopdf}
\usepackage{lineno}
\usepackage{subcaption}
\usepackage{url}
\modulolinenumbers[5]
\usepackage{changes}
\usepackage{soul}
\usepackage{xcolor}
\usepackage[]{algorithm2e}
\usepackage{hyperref}
\hypersetup{
    colorlinks=true,
    linkcolor=blue,
    filecolor=magenta,      
    urlcolor=cyan,
}

\DeclareMathOperator*{\argmin}{arg\,min}

\journal{Journal}

\begin{document}

\begin{frontmatter}


\author[rvt]{Stefanos Nikolopoulos\corref{cor1}}
\ead{stefnikolopoulos@mail.ntua.gr}

\author[rvt]{Ioannis Kalogeris}
\ead{ikalog@central.ntua.gr}

\author[rvt]{Vissarion Papadopoulos}
\ead{vpapado@central.ntua.gr}
\ead[url]{mgroup.ntua.gr}

\cortext[cor1]{corresponding author}
\address[rvt]{MGroup, Engineering Simulations Lab, Institute of Structural Analysis and Antiseismic research, National Technical University of Athens, Zografou Campus, 9 Iroon Polytechniou str, 15780 Zografou}

\title{Non-intrusive Surrogate Modeling for Parametrized Time-dependent PDEs using Convolutional Autoencoders}

\begin{abstract}
This work presents a non-intrusive surrogate modeling scheme based on Machine Learning technology for predictive modeling of complex systems, described by parametrized time-dependent PDEs. For this type of problems, typical finite element solution approaches involve the spatiotemporal discretization of the PDE and the solution of the corresponding linear system of equations at each time step. Instead, the proposed method utilizes a Convolutional Autoencoder in conjunction with a feed forward Neural Network to establish a low-cost and accurate mapping from the problem's parametric space to its solution space. The aim is to evaluate directly the entire time history solution matrix through these interpolation schemes. For this purpose, time history response data are collected by solving the high-fidelity model via FEM for a reduced set of parameter values. Then, by applying the Convolutional Autoencoder to this data set, a low-dimensional representation of the high-dimensional solution matrices is provided by the encoder, while the reconstruction map is obtained by the decoder. Using the latent representation given by the encoder, a feed forward Neural Network is efficiently trained to map points from the problem's parametric space to the compressed version of the respective solution matrices. This way, the encoded time-history response of the system at new parameter values is given by the Neural Network, while the entire high-dimensional system's response is delivered by the decoder. This approach effectively bypasses the need to serially formulate and solve the governing equations of the system at each time increment, thus resulting in a significant computational cost reduction and rendering the method ideal for problems requiring repeated model evaluations or 'real-time' computations. The elaborated methodology is demonstrated on the stochastic analysis of time-dependent PDEs solved with the Monte Carlo method, however, it can be straightforwardly applied to other similar-type problems, such as sensitivity analysis, design optimization, etc. 
\end{abstract}


\end{frontmatter}


\section{Introduction} \label{section1}

Recent advances in the field of computational mechanics has allowed researchers to develop high-fidelity models of complex physical systems that emulate their behavior. With this approach, the response of a system under investigation can be efficiently predicted via computer simulations in lieu of costly and time-consuming experiments. However, certain applications of practical interest such optimization, uncertainty quantification and parameter identification require a large number of model runs. For detailed complex models described by time-dependent PDEs, the computational cost for a single run may range from a few seconds to several hours, hence, the cost of performing these types of analyses becomes unduly expensive. Computational handling of such problems necessitates the development of highly efficient and accurate solution techniques. In this direction, surrogate modeling techniques have emerged over the past years as an effective approach for reducing the computational burden associated with predictive modeling of complex large-scale problems. Surrogate models, also referred to as metamodels, are approximations of the original model that are cheap to compute and can mimic the system's behavior with a controlled loss of accuracy. These models are typically constructed by using some assumptions about the functional shape of the model based on information about the model's response in the form of data, and for this reason they are also known as data-driven models. 

Reduced basis (RB) methods belong to this family of metamodeling techniques and are widely applied as surrogates for parametrized large scale systems \cite{LUCIA2004,Ezvan2016,JENSEN2016,MIGNOLET2013}. The idea behind RB methods is to find a suitable low-dimensional subspace of the system's high-dimensional solution space and project the governing equations onto this reduced space, where they can be solved more efficiently. The most popular linear reduced basis technique is Proper Orthogonal Decomposition (POD) \cite{RATHINAM2003,AMSALLEM2008,FARHAT2015, SENGUPTA20042693}, also known as Karhunen-Lo\'{e}ve expansion or Principal Component Analysis in certain contexts. POD is typically applied to a collection of solution vectors (snapshots) and identifies an appropriate basis for a lower dimensional subspace. The main advantage of POD stems from its ability to optimally truncate the basis such that it represents only the most energetic modes contained in the snapshots. Other linear basis construction methods include proper generalized decomposition \cite{Almeida2013, CHINESTA2011}, balanced truncation \cite{Moore1981, Safonov1989} and rational interpolation \cite{Baur2011}.

While linear RB methods have been demonstrated to work optimally on linear problems, this is not the case for general nonlinear problems with non-affine dependence on parameters \cite{Nguyen2008}. This is because in such cases the system configuration needs to be updated at each nonlinear iteration or at any new parameter value and this process can only be performed on the full-order model. Therefore, every time the system changes, the reduced system of equations needs to be re-derived using Galerkin projections, which translate to multiple inner product evaluations. However, these evaluations are costly and they significantly diminish the computational gains of linear RB methods. To address nonlinear problems with non-affine parameter dependence, several RB schemes based on the empirical interpolation method \cite{Chaturantabut2010, NEGRI2015} or subspace-angle interpolation \cite{AMSALLEM2008, AMSALLEM2012} have been proposed, but these are also intrusive in nature and their generalization to other nonlinear problems is not straightforward.

Recently, the combination of RB techniques with data-driven machine learning models has resulted in non-intrusive approaches for the solution of large-scale complex systems. The advantage of these methods is that they do not need to access and modify the governing equations of the original high-fidelity model. For instance, in \cite{HESTHAVEN2018, Park2013} it has been proposed to combine POD and Neural Networks (NNs) producing a hybrid POD-NN approach, where the NN was trained to produce the low-dimensional projection coefficients of the RB model. In this frame, the use of different interpolation schemes instead of NNs, such as Gaussian Process Regression (GPR) \cite{GUO2018} and radial basis functions \cite{XIAO2017, DEBOER2007784} were also shown to be very efficient for interpolating over the POD coefficients. Despite the fact that these methods are highly efficient, their main drawback is that for general nonlinear problems, they often require a higher number of model evaluations than intrusive methods to construct a reliable surrogate in the first place.

Motivated by the inability of linear reduction methods such as POD to capture complex response surfaces, nonlinear manifold learning methods (e.g. Kernel PCA \cite{ZHOU2020106358}, Hessian eigenmaps \cite{YE2015197}, Laplacian eigenmaps \cite{BELKIN2003}, local tangent space alignment \cite{ZHANG2004}, the diffusion maps algorithm \cite{COIFMAN20065}) gained more attention over the past few years. The main assumption in manifold learning is that the data points, which correspond to system solutions in this setting, lie on a low-dimensional manifold embedded in an ambient higher-dimensional Euclidean space. The goal is to identify the manifold's intrinsic dimensionality, that is, the parameters that describe it, and thus obtain low-dimensional representations of the data set. This approach can remedy the problems associated with the curse of dimensionality when dealing with high-dimensional data sets and, consequently, enable the development of efficient interpolation schemes. For instance, in \cite{Lataniotis2018}, the kernel PCA algorithm was employed for the purposes of dimensionality reduction and in conjunction with Kriging and polynomial chaos expansion surrogates, a cost-efficient metamodel was constructed. Similarly, in \cite{Kalogeris2020} the diffusion maps algorithm has been investigated as an alternative to POD.

Despite the effectiveness of above-mentioned algorithms in providing low-dimensional representations for high-dimensional data sets, their main drawback stems from the fact that they do not provide an analytic relation for decoding the compressed data back to their high-dimensional representations in the original space. This problem is known in the literature as the pre-image problem and several elaborate interpolation schemes have been employed to address it, such as the geometric harmonics \cite{COIFMAN2006b} and Laplacian pyramids \cite{Burt1983}. However, a more versatile solution to this problem can be provided by the Autoencoders \cite{LOPEZPINAYA2020193}. An Autoencoder (AE) is a specific type of an unsupervised Neural Network that learns how to efficiently compress and encode data and then learns how to reconstruct (decode) them, that is, to map them from their encoded representation to a representation as close to the original input as possible. The encoder and decoder parts of an Autoencoder are trained jointly, yet can be used separately. An extension of ordinary Autoencoders are the so called Convolutional Autoencoders (CAEs), which have been developed primarily for spatial field data compression but have proven particularly useful in several applications dealing with high-dimensional data sets. Some of these applications pertain to the fields of computer vision \cite{KOLOTOUROS2019}, pattern recognition \cite{OYEDOTUN2016} and time series data prediction \cite{ZHAO2011}. Similarly to ordinary AEs, CAEs also consist of an encoder and a decoder part but they are constructed using different types of layers, called convolutional and deconvolutional layers \cite{Guo2016}. 

In this work, a non-intrusive surrogate modeling strategy is proposed for the solution of problems described by parametrized time-dependent PDEs. This scheme relies on the powerful dimensionality reduction properties of CAEs, which are exploited as a means of encoding and decoding the high-dimensional solution data sets. Furthermore, feed forward NNs (FFNNs) are used to establish a mapping between the problem's parametric space to its encoded solution space. With this approach, the encoded time-history response of the system at a new parameter value is given by the FFNN, while its representation in the original high-dimensional space is obtained by the decoder. A similar approach can be found in \cite{XU2020113379}, where the authors suggest the use of 3 levels of NNs, namely a CAE, a temporal CAE \cite{oord2016wavenet} and a FFNN to perform parameter and future state prediction. On the other hand, the surrogate scheme proposed herein requires only 2 levels of NNs, a FFNN and a CAE, rendering it very easy to implement. Furthermore, in terms of performance, our investigation indicated that the optimal CAE's architecture is based on 1-D convolutional filters scanning the input data only in the time axis along with 1-D average pooling layers. This way, the proposed approach has reduced offline and online computational requirements, while at the same time achieves very accurate results. Therefore, it is capable of providing remarkably fast and accurate evaluations of the complete system's response, effectively bypassing the need to serially formulate and solve the governing equations of the system at each time increment, as is typically required by FEM methods. The elaborated methodology is applied to the stochastic analysis of time-dependent PDEs, parametrized by the system's random variables and solved in the frame of the Monte Carlo method.


The paper is organised as follows. In section 2 the general framework of Autoencoders is presented with emphasis on CAEs. In section 3,  the proposed hybrid CAE-FFNN surrogate modeling scheme for parametrized time-dependent PDEs is introduced. In section 4, numerical examples for testing the method are provided and Section 5 concludes with a brief summary and possible extensions. Codes and datasets accompanying this manuscript are available at \url{https://github.com/visten92/CAE-FFNN}.

\section{General framework of Autoencoders} \label{section2}
\subsection{Autoencoder} \label{section2.1}
The Autoencoder (AE) concept was introduced in \cite{RUMELHART1986} and is regarded as a Neural Network that learns from an unlabeled data set in an unsupervised manner. The aim of an Autoencoder is to learn a reduced representation for a set of data, referred to as encoding, and then to learn how to reconstruct the original input from the encoded input with the minimum possible error. The latter part of the AE is called decoder. 

In particular, let $\boldsymbol{X}$ be a subset of $\mathbb{R}^d$ with $\boldsymbol{x}\in \boldsymbol{X}$ denoting an element of the set. Then, the AE's encoder and decoder are defined as transition maps $\phi$, $\psi$ such that:

\begin{align}
    &\phi : \boldsymbol{X} \subseteq \mathbb{R}^d \rightarrow \boldsymbol{H}\subseteq \mathbb{R}^l \\
    &\psi : \boldsymbol{H} \subseteq \mathbb{R}^l \rightarrow \boldsymbol{X}\subseteq \mathbb{R}^d \\
    &\phi ,\psi = \argmin_{\phi,\psi} \| \boldsymbol{X}-(\psi \circ \phi)\boldsymbol{X} \|^2
\end{align}

\noindent with the dimension $l$ typically being much smaller than $d$.

Now, let us consider the simplest case, where the encoder has only one hidden layer. It takes an input $\boldsymbol{x}\in\mathbb{R}^{d}$ and sends it to $\boldsymbol{h}=\phi(\boldsymbol{x})\in\mathbb{R}^l$, which in this case can also be written as

\begin{equation} \label{encoderEq}
\boldsymbol{h}=\sigma(\boldsymbol{W}\boldsymbol{x} + \boldsymbol{b}) 
\end{equation}

\noindent with $\sigma$ being an activation function (eg. $tanh, \ ReLU$, etc), $\boldsymbol{W}$ a weight matrix and $\boldsymbol{b}$ a bias vector. The image $\boldsymbol{h}$ of $\boldsymbol{x}$ is the latent or encoded representation of $\boldsymbol{x}$ and $\boldsymbol{H}$ is the latent or feature space.

\begin{figure}[H]
    \centering
    \includegraphics[width=0.55\textwidth]{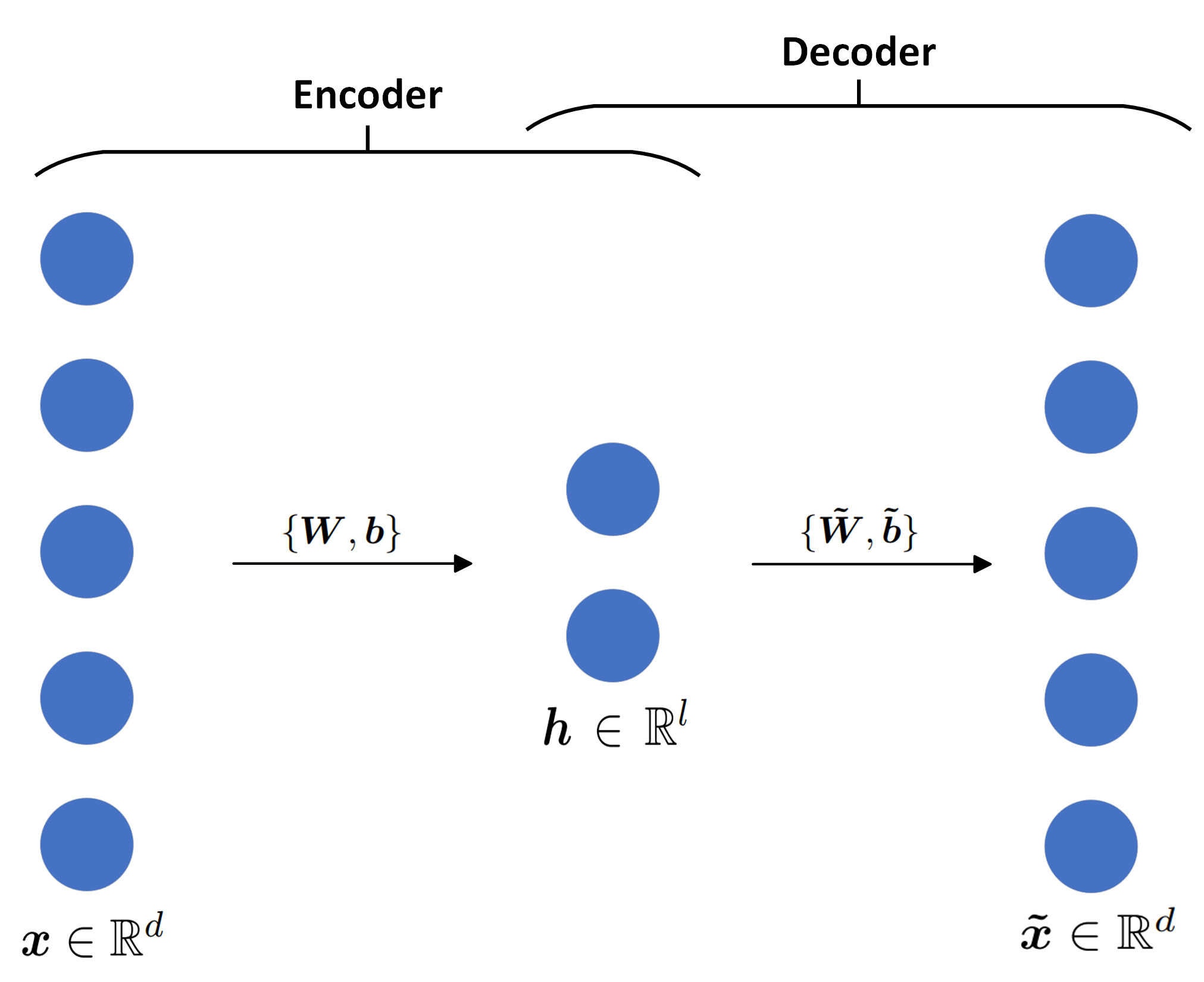}
    \caption{Schematic representation of a basic Autoencoder}
    \label{fig:AEconcept}
\end{figure}

The decoder's task is to establish the inverse mapping $\psi$ that will reconstruct the input $\boldsymbol{x}$, given its latent representation $\boldsymbol{h}$. Again, considering a one-hidden layer, the reconstructed point $\tilde{\boldsymbol{x}}=\psi(\boldsymbol{h})$ is given by

\begin{equation} \label{decoderEq}
\tilde{\boldsymbol{x}}=\tilde{\sigma}(\tilde{\boldsymbol{W}}\boldsymbol{h} + \tilde{\boldsymbol{b}}) 
\end{equation}

\noindent where $\tilde{\sigma}$, $\tilde{\boldsymbol{W}}$ and $\tilde{\boldsymbol{b}}$ may be unrelated to those of encoder. Also, the network's architecture selected for the encoder can be different than the decoder's and the number of hidden layers can be greater than one, leading to the so-called deep Autoencoders. The general concept and architecture of an Autoencoder is schematically presented in figure \ref{fig:AEconcept}. 

AEs are trained by a back propagation algorithm \cite{buscema}, which is the most commonly used algorithm for the training of NNs. Back propagation computes the gradient of the loss function with respect to network's weights very efficiently with the aid of automatic differentiation (AD) \cite{Baydin}. AD involves a set of techniques developed to numerically evaluate the gradient of a function specified by a computer program. It exploits the fact that every operation performed by the program, no matter how complicated, executes a sequence of elementary arithmetic operations (addition, subtraction, multiplication, division, etc.) and elementary functions ($exp$, $log$, $sin$, $cos$, etc.). By applying the chain rule to these operations, derivatives of arbitrary order can be computed to working precision. Thus, gradient based optimization methods such as stochastic gradient descent, adaptive moment, etc. can be applied for training multilayer NNs by updating weights such as to minimize loss. 

In the context of AEs, the loss function becomes the reconstruction error between the input points $ \boldsymbol{x}_i $ and their respective output $\tilde{\boldsymbol{x}}_i $. It is usually expressed as the mean-square error:

\begin{equation}\label{eq2.1}
\mathcal{L} = \frac{1}{N}\sum_{i=1}^{N}||\boldsymbol{x}_{i} - \tilde{\boldsymbol{x}}_{i}||_{2}^{2}
\end{equation}

\noindent with $\| \cdot \|_2$ denoting the $L^2$-norm and $N$ being the number of points in the training data set. It should be explicitly mentioned that even though the minimization of the reconstruction error implies that the encoder and decoder are trained jointly, however, they can be used separately. 

\subsection{Convolutional Autoencoder}
\label{section2.2}
 
 Despite their powerful dimensionality reduction properties, AEs face significant challenges when dealing with very high-dimensional inputs, due to the fact that the number of trainable parameters increases drastically with an increase in the input's dimensionality. In addition, AEs are not capable of capturing the spatial features of the input (e.g. when dealing with images) nor the sequential information in the input (e.g. when dealing with sequence data). 
 
 To remedy these issues, a new type of Autoencoders has emerged, that of Convolutional Autoencoders (CAEs) \cite{MASCI2011}. Similarly to AEs, CAEs also consist of an encoder and a decoder that are trained to minimize the loss function of eq. \eqref{eq2.1}, but they are built from different layer types. Specifically, in CAEs the encoder part is built using a combination of convolutional layers, fully connected layers, pooling layers and normalization layers, while the decoder is built from deconvolutional layers and unpooling layers along with fully connected and normalization layers.  Intuitively, CAEs can be viewed as extensions of ordinary AEs in the same way that Convolutional Neural Networks \cite{KRIZHEVSKY2011} are extensions of fully-connected Neural Networks. These concepts are illustrated in the following sections.

\subsubsection{Convolutional and Deconvolutional Layers}
\begin{figure}[H]
    \centering
    \includegraphics[width=0.6\textwidth]{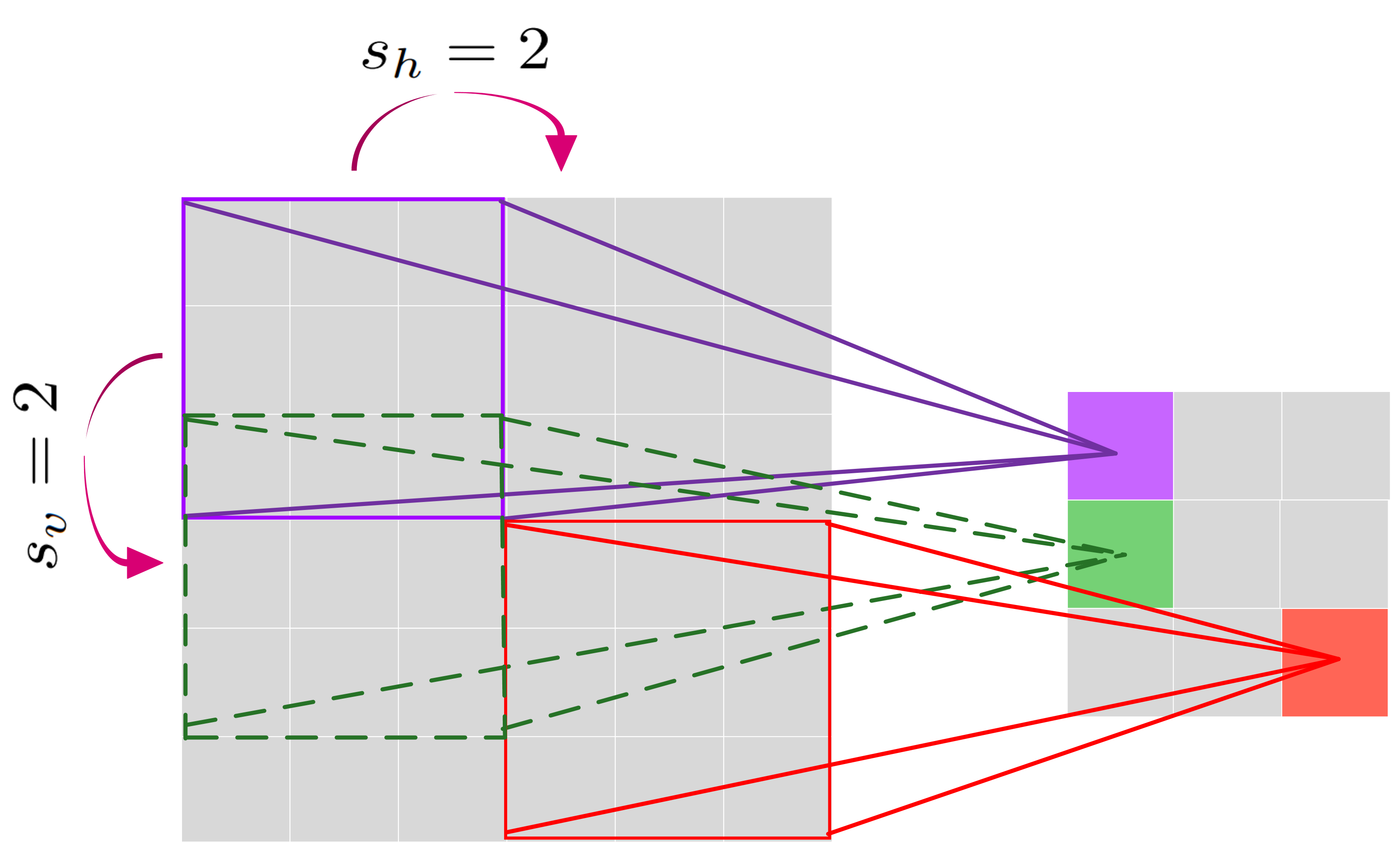}
    \caption{Schematic representation of a 2-D Convolutional filter with strides $s_{h} = 2$ and $s_{v} = 2$.}
    \label{fig:convLayer}
\end{figure}

Convolutional layers take as input a $n$-D array $\boldsymbol{M}$ and apply a filter $\mathcal{F}$ (a.k.a. kernel) of specified size to the elements of $\boldsymbol{M}$ in a moving window fashion. This process is schematically depicted in figure \ref{fig:convLayer}. Essentially, the objective of the convolution operation is to extract the most important features from the input and use them to encode it. To better clarify this process, let us consider a $2-D$ array $\boldsymbol{M}=[m_{ij}]$ and its encoded version $\boldsymbol{M}^{enc}=[\mu_{ij}]$, called feature map, which is obtained after applying a filter $\boldsymbol{W}=[w_{ij}]$ of size $f_h \times f_w $, moving with vertical stride $s_v$ and horizontal stride $s_h$. The element $\mu_{ij}$ of $\boldsymbol{M}^{enc}$ is given by the equation:

\begin{equation}\label{eq2.2}
\mu_{ij}=\sum_{u=1}^{f_h}\sum_{v=1}^{f_w}m_{i'j'}\cdot w_{uv}+b \ \ \text{with} \begin{cases} i'=i\times s_v+u \\ j'=j \times s_h+v
\end{cases}
\end{equation} 

\noindent where $b$ is the bias term and $w_{uv}$ is the element of the filter $\boldsymbol{W}$ that gives the connection weight between elements of $\boldsymbol{M}^{enc}$ and the elements of $\boldsymbol{M}$ within the respective window.

This layer architecture is significantly more economical than that of a fully connected layer since the parameters involved are only the $f_h \times f_w $ elements of the filter $w_{ij}$ and the bias term $b$. The filter parameters do not require to be manually defined, instead the Convolutional layer will automatically learn the most appropriate filter for the task.
Also, a Convolutional layer can have multiple filters, in which case it outputs one feature map $\boldsymbol{M}^{enc}_k$ per each filter $k$. This enables it to detect multiple features anywhere in its inputs. Additionally, several convolutional layers can be stacked in order to build deep architectures which allow the network to concentrate on small low-level features in this first layer and progressively assemble them into larger higher-level features in the subsequent layers. In this more general case, the element $\mu_{ijk}$ at the $q$-th Convolutional layer, corresponding to row $i$, column $j$ of the $k$ feature map $\boldsymbol{M}^{enc}_k$ is obtained as:

\begin{equation}\label{eq2.2generic}
\mu_{ijk}=\sum_{u=1}^{f_h}\sum_{v=1}^{f_w}\sum_{k'=1}^{f_{n'}}m_{i'j'k'}\cdot w_{uvk'k}+b_k \ \ \text{with} \begin{cases} i'=i\times s_v+u \\ j'=j \times s_h+v
\end{cases}
\end{equation} 

\noindent where now $f_{n'}$ is the number of feature maps in the previous layer (layer $q-1$), $m_{i'j'k'}$ the value located in row $i'$, column $j'$ of the $q-1$ layer's feature map $k'$ and $b_k$ is the bias term for the $k$-th feature map (in layer $q$). Also,
$w_{uvk'k}$ is the connection weight between the values in feature map $k$ of layer $q$ and its input located at row $u$, column $v$ at the window of the $k'$ feature map. To simplify the notation, the application of several convolutional layers, with multiple filters each, to an array $\boldsymbol{M}$ will be expressed as

\begin{equation}\label{eq2.3}
\boldsymbol{M}^{enc}=ConvNN(\boldsymbol{M})
\end{equation} 

\noindent with $ConvNN(\cdot)$ denoting the mapping from the initial input space to its encoded representation. 

\begin{figure}[H]
    \centering
    \includegraphics[width=0.6\textwidth]{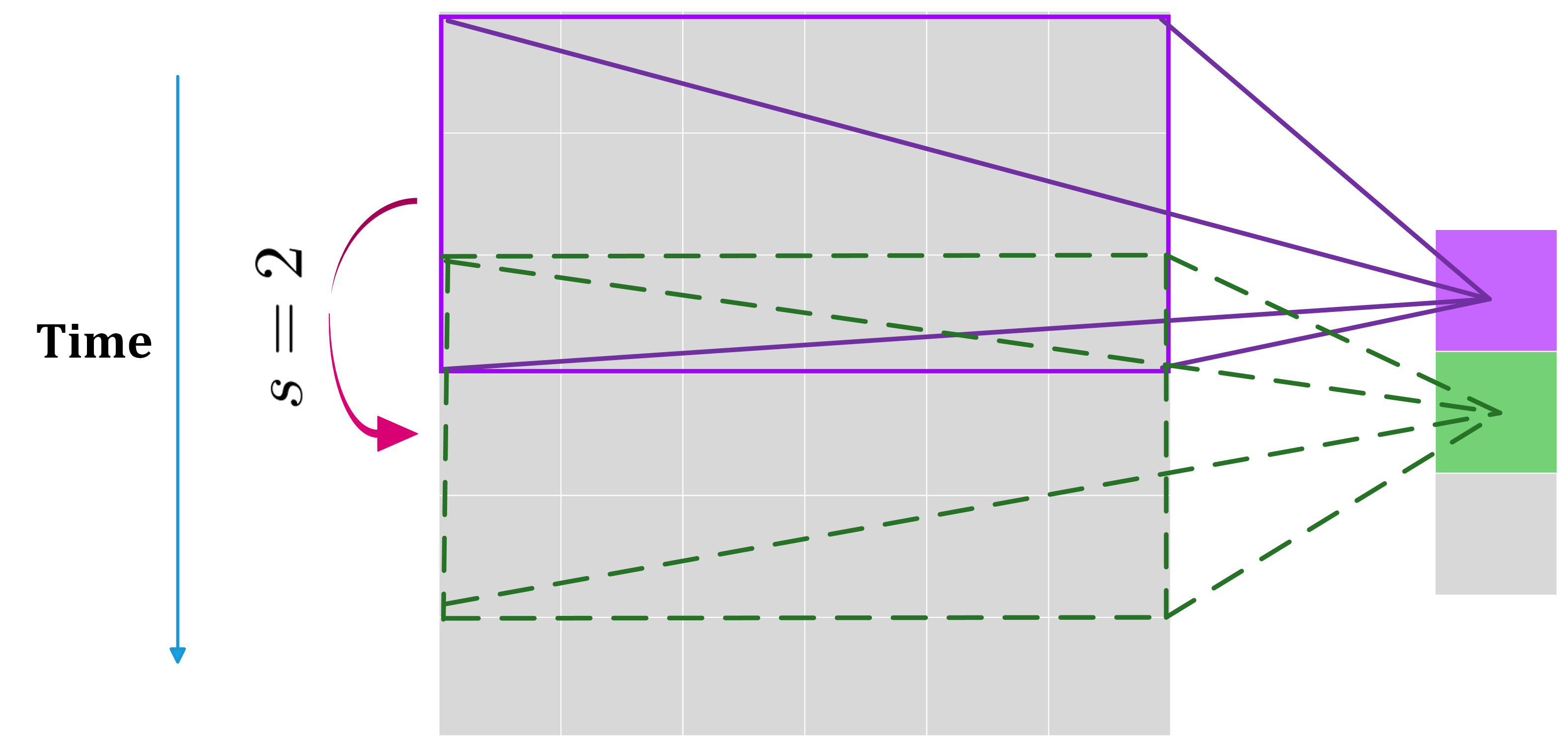}
    \caption{Schematic representation of a 1-D Convolutional filter with stride $s = 2$.}
    \label{fig:1DconvLayer}
\end{figure}

Depending on the application, the convolutional filters can either be one, two or three dimensional with the difference between them being the way they slide across the data. In this work, the focus is on processing time series data, therefore 1-D convolutional filters, such as the one depicted in figure \ref{fig:1DconvLayer}, were used to scan the data only in the time axis.  

On the other hand, a Deconvolutional layer performs the reverse operation of convolution, called deconvolution, and it is used to construct decoding layers. Their function is to multiply each input value by a filter elementwise. For instance, a 2D $f_h \times f_w$ deconvolution filter maps an $1\times 1$ spatial region of the input to an $f_h \times f_w$ region of the output. Thus, the filters learned in the deconvolutional layers create a base used for the reconstruction of the inputs' shape, taking into consideration the required shape of the output. As before, a Deconvolutional layer can have multiple filters, while several deconvolutional layers can be stacked for building deep architectures for CAEs \cite{HASAN2016, Guo2016}. The decoding procedure can be represented as:

\begin{equation}\label{eq2.6}
\boldsymbol{M}=DeconvNN(\boldsymbol{M}^{enc})
\end{equation} 

\begin{figure}[H]
    \centering
    \includegraphics[width=1\textwidth]{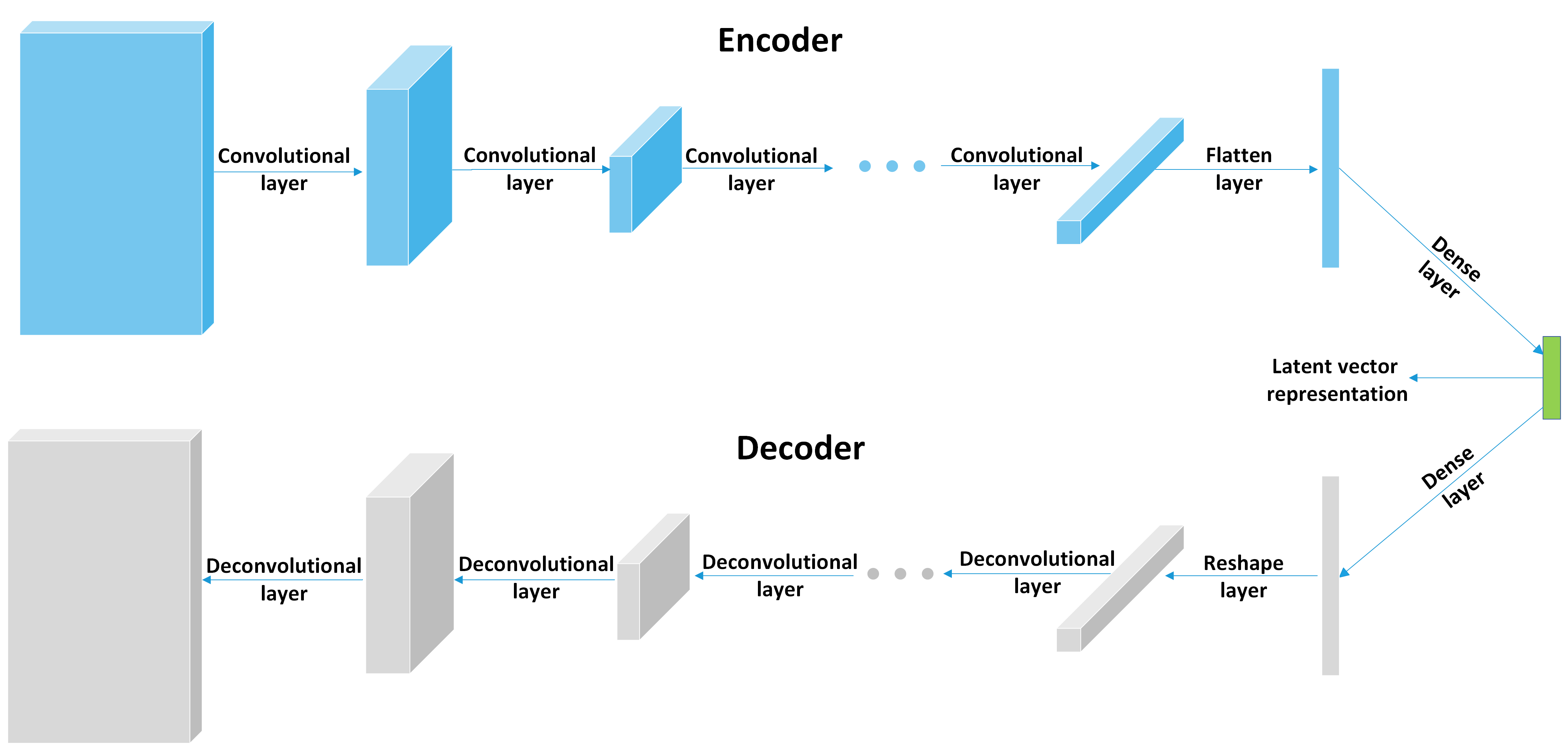}
    \caption{Schematic representation of a deep Convolutional Autoencoder.}
    \label{fig:CAE_Ex}
\end{figure}

Based on the above, the CAE's architecture consists of Convolutional, Deconvolutional and dense layers and is typically used for dimensionality reduction and reconstruction purposes. In practice, the CAE's encoder uses a number of convolutional layers to compress the input and once the desirable level of reduction has been achieved, the encoded matrix is flattened into a vector. Then, a dense layer is employed to map this vector to its latent representation. In the reverse direction, the decoder starts by taking the latent representation and transforming it into a vector through a denser layer. Subsequently, the input reconstruction is achieved by the deconvolutional layers. In accordance to eq. \eqref{eq2.1} the loss for CAEs becomes:

\begin{equation}\label{eq2.1CAE}
\mathcal{L} = \frac{1}{N}\sum_{i=1}^{N}||\boldsymbol{M}_{i} - \tilde{\boldsymbol{M}}_{i}||_{2}^{2}
\end{equation}

\noindent where $\boldsymbol{M}_{i}$ denotes the input arrays used for training and $\tilde{\boldsymbol{M}}_i=DeconvNN(ConvNN(\boldsymbol{M}_i))$ the corresponding CAE's output. In figure \ref{fig:CAE_Ex}, a schematic representation of a deep CAE is presented.  

\subsubsection{Pooling and Unpooling Layers}
Aside of Convolutional, Deconvolutional and dense layers, two other important layer types often employed in CAEs are those of pooling and unpooling. Pooling layers are quite similar to convolutional layers in the sense that they downsample the input in order to decrease its size, however, they do not involve any trainable parameters. Their goal is to reduce the computational load, the memory usage, and the number of parameters. The latter is particularly useful since it also limits the risk of overfitting. Each neuron in a pooling layer is linked to a limited number of neurons in the previous layer, located within a small window. The window's size and stride are user defined. 

Common types of pooling layers include the max pooling layer and the average pooling layer. The first outputs the maximum value from the portion of the input covered by the filter and all other inputs are neglected. Accordingly, average pooling layers return the average from the portion of the input. Aside from its dimensionality reduction properties, the pooling operation can be useful for extracting dominant features of the input such as translational, rotational and scale invariance. Nevertheless, caution should be exercised regarding the usage of pooling layers because the corresponding accuracy loss might outweigh the benefits they provide.

\begin{figure}[H]
    \centering
    \includegraphics[width=0.80\textwidth]{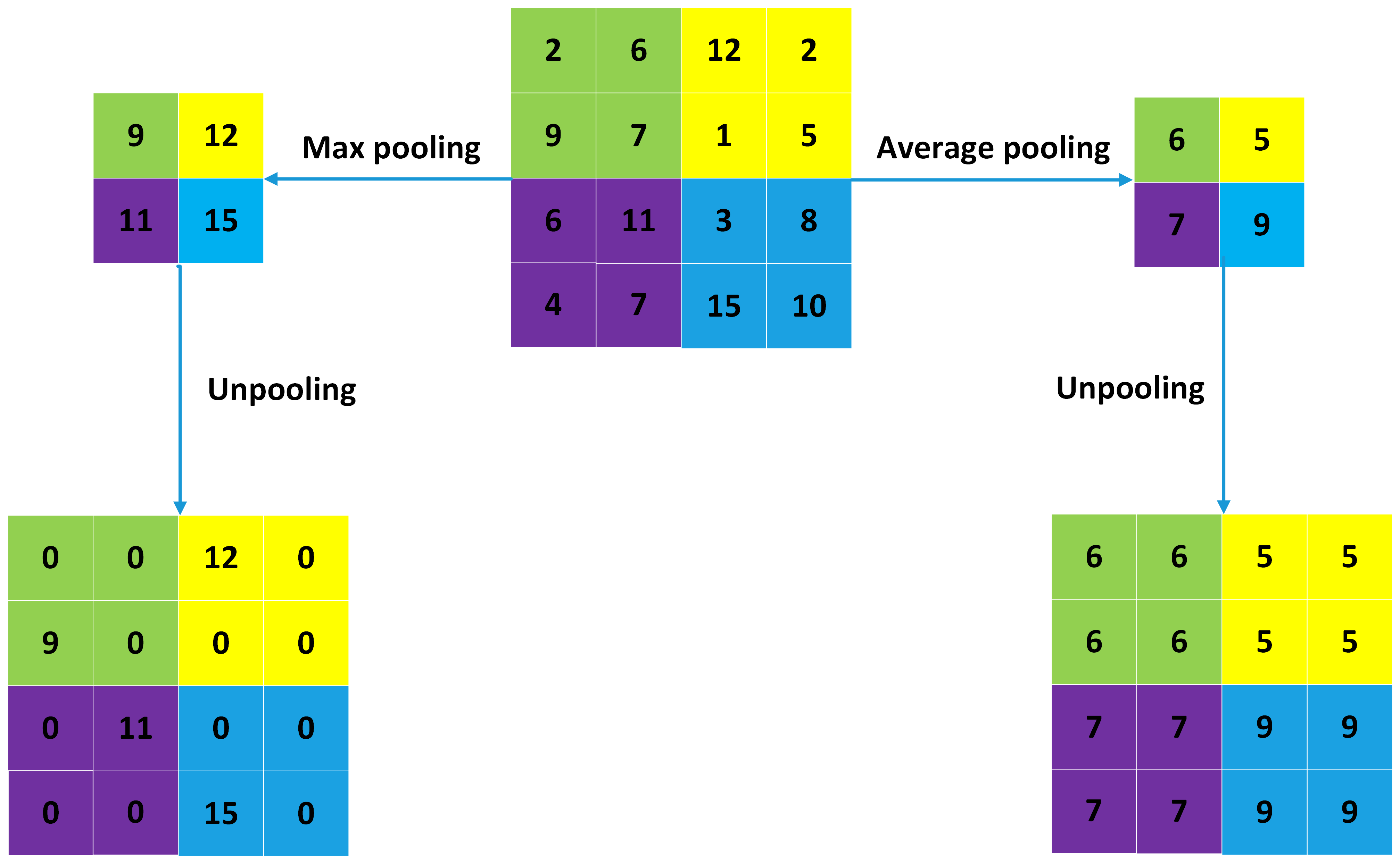}
    \caption{Examples of pooling and unpooling.}
    \label{fig:pooling}
\end{figure}

On the other hand, unpooling layers perform the reverse operation of pooling and their aim is to reconstruct the original size of each rectangular patch. During the max pooling operation, a matrix is created which records the location of the maximum values selected during pooling. This matrix is then employed in the unpooling operation in order to place each value back to its original pooled location, while setting all other values to zero. In the case of average unpooling, it assigns the same mean value to all elements of the output window. A schematic representation of max pooling, average pooling and unpooling is given in figure \ref{fig:pooling}.

\section{Surrogate modeling of parametrized time - dependent PDEs using Convolutional Autoencoders} \label{section3}
Consider the modeling of a parametrized physical system governed by partial differential equations:

\begin{align} \label{eq3.1}
\frac{\partial u\left(\boldsymbol{x},t;\boldsymbol{\theta}\right)}{\partial t} + \mathcal{N}\left( u\left(\boldsymbol{x},t;\boldsymbol{\theta}\right) \right)&=f(\boldsymbol{x},t;\boldsymbol{\theta}), \ \ \  \boldsymbol{x} \in \Omega, t\in[0,T], \boldsymbol{\theta}\in\Theta \nonumber \\
\mathcal{B}(u(\boldsymbol{x},t;\boldsymbol{\theta})&=b(\boldsymbol{x},t;\boldsymbol{\theta}), \ \ \  \boldsymbol{x} \in \partial\Omega, t\in[0,T], \boldsymbol{\theta}\in\Theta \\
\mathcal{C}(u(\boldsymbol{x},0;\boldsymbol{\theta})&=c(\boldsymbol{x};\boldsymbol{\theta}), \ \ \ \  \boldsymbol{x} \in \partial\Omega, \boldsymbol{\theta}\in\Theta \nonumber
\end{align}

\noindent where $u\left(\boldsymbol{x},t;\boldsymbol{\theta}\right)$ is the field of interest, $\mathcal{N}$ is a general differential operator that involves spatial derivatives, and $f(\boldsymbol{x},t;\boldsymbol{\theta})$ is a source field. Furthermore, $\mathcal{B}$ is the operator for the boundary conditions defined on the boundary $\partial \Omega$ of the domain $\Omega$, $\mathcal{C}$ is the operator for the initial conditions at $t=0$ and $\boldsymbol{\theta}\in \Theta$ is a vector of uncertain parameters that include randomness in the system parameters, loading or boundary conditions.

The discrete solution to the above set of equations for a given parameter value $\boldsymbol{\theta}$ can be obtained through the semidiscrete Galerkin method. Specifically, the spatial part of the solution is obtained through the finite element method on a discrete space $\mathcal{V}_h$ spanned by basis functions $\varphi_i(\boldsymbol{x})$, $i=1,2,...,d$, with $d$ being the number of degrees of freedom. To take into account the time-dependence, temporal derivatives are approximated by finite differences. Thus, the FE solutions $\mathbf{u}_h$ are expressed as:

\begin{equation} \label{eq3.2}
\mathbf{u}_h(\boldsymbol{x};\boldsymbol{\theta},t)=\sum_{i=1}^d \left(\boldsymbol{u}_h(\boldsymbol{\theta},t)\right)_i\varphi_i(\boldsymbol{x})
\end{equation}

\noindent with $\boldsymbol{u}_h(\boldsymbol{\theta},t) \in \mathbb{R}^d$ being the expansion coefficient vector at time $t$ for a given parameter value $\boldsymbol{\theta}$. Then, the complete time-history response of the system is given by $\boldsymbol{U}_h (\boldsymbol{\theta})=\lbrace \boldsymbol{u}_h(\boldsymbol{\theta},t_1), \cdots, \boldsymbol{u}_h(\boldsymbol{\theta},t_{N_t}) \rbrace$ $\in \mathbb{R}^{d\times N_t}$, where $N_t$ is the number of time increments in the temporal discretization.

To quantify the probabilistic characteristics of the solution in eq. \eqref{eq3.1} the most versatile approach is the brute force Monte Carlo (MC) simulation. In this setting, a large number, $N_{MC}$, of  parameter realizations $\lbrace\boldsymbol{\theta}_j\rbrace_{j=1}^{N_{MC}}$ is generated according to their joint probability distribution and the corresponding PDEs are solved with the FEM procedure in order to obtain an accurate estimate of the system's stochastic behaviour. Namely, for a set of parameter values, the PDEs are discretized as described above and the corresponding linear system of equations is solved at each time step either directly or iteratively. Then, the system responses are statistically processed to extract the probabilistic characteristics of the response. Evidently, MC analysis of this type is associated with increased computational requirements, espesially when handling large - scale problems where the cost of each model run may range from several minutes to several hours. 

To alleviate this computational burden, a surrogate modeling approach is proposed herein based on the powerful dimensionality reduction capabilities of CAEs. To this purpose, the PDEs are solved with the classic FEM procedure for a small, yet sufficient number, $N$, of parameter values in order to obtain a data set of time history matrices $\lbrace\boldsymbol{U}_i\rbrace_{i=1}^{N}$. The CAE (encoder and  decoder) is trained over this data set minimizing the reconstruction mean square error. The encoded representation of each time history solution matrix $\boldsymbol{U}_i$ is a low dimensional vector $\boldsymbol{z}_i \in \mathbb{R}^l$ ($l<<d$), which allows a feed forward Neural Network (FFNN) to be trained accurately and efficiently in order to construct a mapping between the PDE's parametric space and the encoded solution space. After the training phase is completed the proposed surrogate scheme works as follows. For a new input parameter vector, the encoded vector representation of the time history solution matrix is calculated by the FFNN and, subsequently, the entire time history matrix is delivered by the CAE's decoder. This way a large number of MC simulations can be performed afterwards at a minimum computational cost.  The implementation steps of the proposed approach can be divided into two phases, namely the offline and the online phase, and these are the following:

\

\textbf{\textit{Offline phase}}
\begin{enumerate}
     \item[\indent {}] \textit{Step 1}: Generate $N$ vectors of parameter values $\boldsymbol{\theta}_{i} \in\mathbb{R}^{n}$ with $i=1,2,....N$  according to their probability distribution and solve the corresponding time-dependent PDEs with the FEM procedure. Collect the solutions in a three-dimensional array $N$ x $d$ x $N_{t}$ , where $d$ is the number of degrees of freedom and $N_t$ the number of time increments.
    \item[\indent {}] \textit{Step 2}: Train a CAE over the $N$ time history solutions matrices $\boldsymbol{U}_{i} \in\mathbb{R}^{d \times N_{t}}$, collected in step 1, to obtain the encoded low dimensional vector representations $\boldsymbol{z}_{i}\in\mathbb{R}^{l}$ of these matrices along with the reconstruction map.
    \item[\indent {}] \textit{Step 3}: Train a FFNN to establish a mapping from the parametric space $\boldsymbol{\theta}_{i}$ to the low dimensional encoded space $\boldsymbol{z}_{i}$.
\end{enumerate}

\begin{figure}[H]
    \centering
    \includegraphics[width=0.75\textwidth]{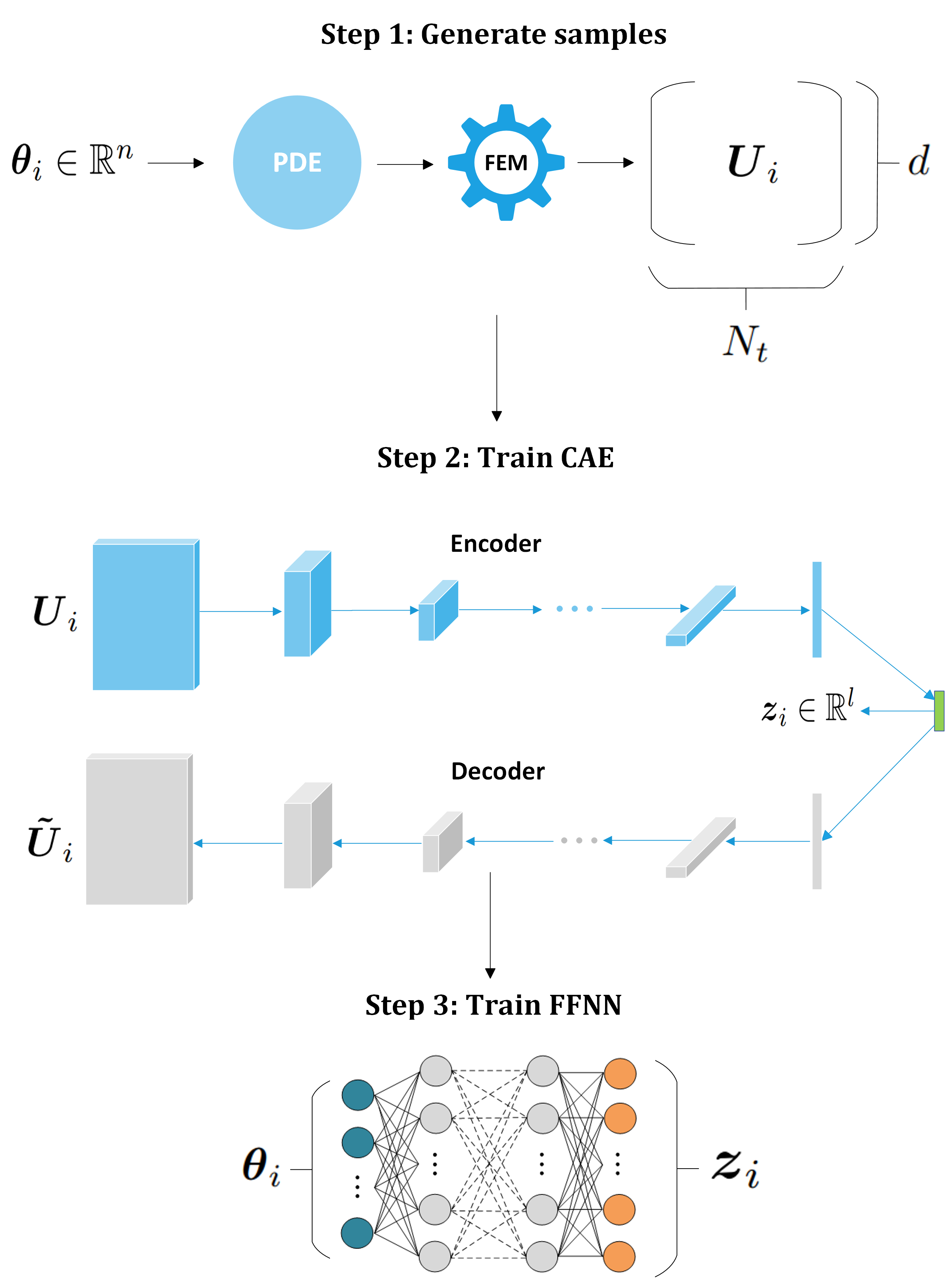}
    \caption{Offline phase of the proposed surrogate modeling method}
    \label{fig:Offline}
\end{figure}

Steps 1-3 of the offline phase are illustrated in fig. \ref{fig:Offline}.
 
\newpage

\textbf{\textit{Online phase}}
\begin{enumerate}
    \item[\indent {}] \textit{Step 1}: For $N_{MC}$ new realizations of parameter vectors $\boldsymbol{\theta}_{j}$ with $j=1,2,....N_{MC}$, generated from the same joint probability distribution, use the trained FFNN to obtain the encoded vector representations of the solution matrices $\boldsymbol{z}_{j}$.
    \item[\indent {}] \textit{Step 2}: The CAE's decoder is used to produce the solution matrices $\boldsymbol{U}_{j}$ based on their encoded representations $\boldsymbol{z}_{j}$ in the previous step. 
\end{enumerate}

\begin{figure}[H]
    \centering
    \includegraphics[width=0.65\textwidth]{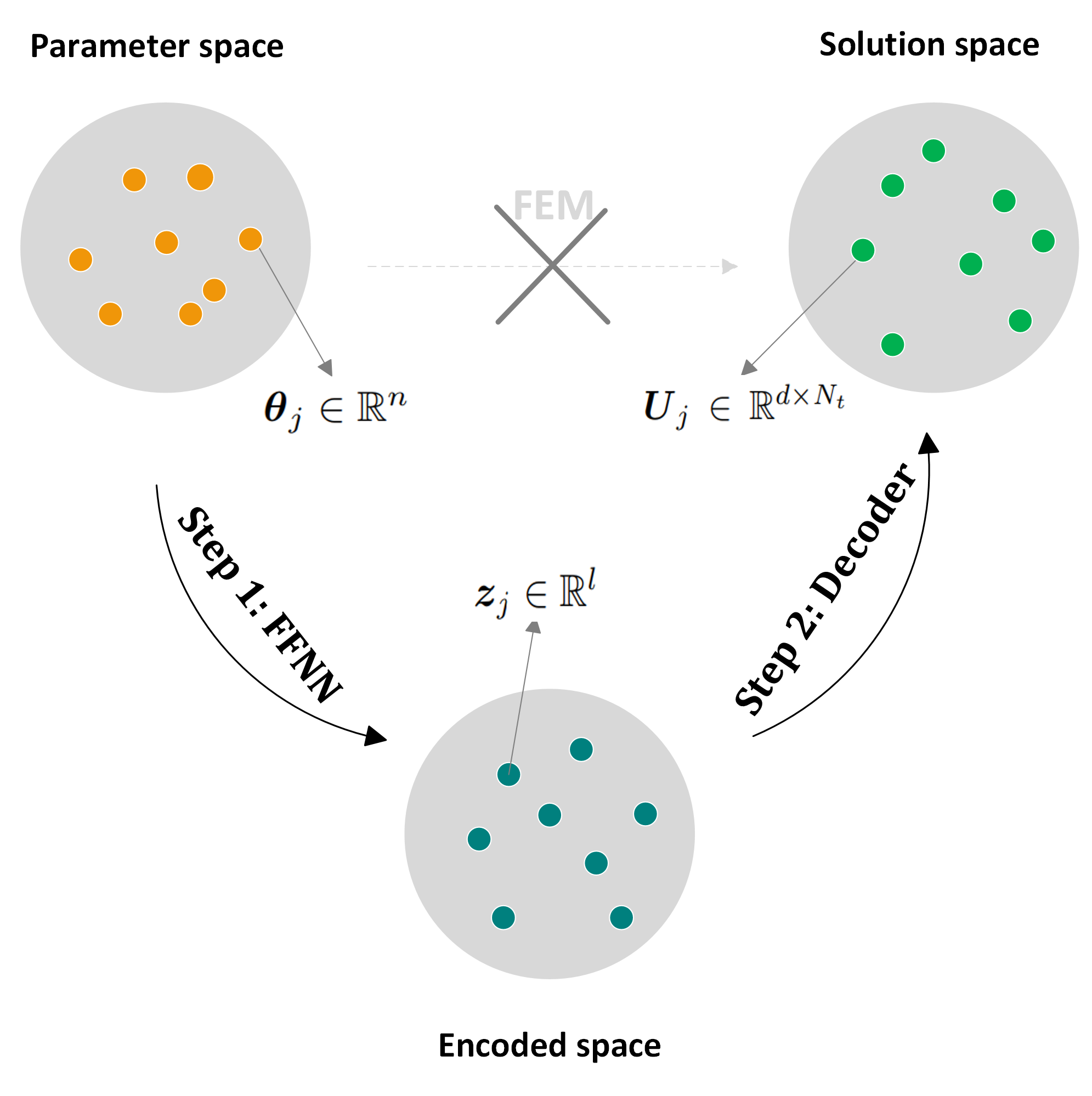}
    \caption{Online phase of the proposed surrogate modeling method}
    \label{fig:Online}
\end{figure}

Steps 1 and 2 of the online phase are schematically represented in fig. \ref{fig:Online}.

\section{Numerical Tests}

We first implement the proposed methodology in the academic case of a 1-D non-linear Burgers' equation, in order to illustrate the applicability. The efficiency and accuracy of the method are assesed subsequently on a more complex structural problem governed by the equations of 2D linear elasticity. All analyses were performed on the open-source computational mechanics software platform MSolve \cite{MSolve}.

\subsection{Burgers' Equation}

Burgers' equation is occurring in many fields of engineering and applied mechanics, such as fluid mechanics \cite{BAKER2000230} and non-linear acoustics \cite{LOMBARD2013406}. It is a convection-diffusion PDE of the following form:

\begin{equation} \label{Burgers}
\frac{\partial \boldsymbol{u}}{\partial t} + \nabla \boldsymbol{u} \cdot \boldsymbol{u} =\nu  \nabla^{2} \boldsymbol{u}
\end{equation}

\noindent where $\boldsymbol{u} \equiv \boldsymbol{u}(\boldsymbol{x},t) $ is the velocity field of the fluid, $\nabla \boldsymbol{u}$ is its gradient and $\nu$ is the fluid's viscosity. For simplicity we choose to demonstrate the 1-D version of equation \eqref{Burgers} which may be written as:

\begin{equation} \label{Burgers1D}
\frac{\partial u}{\partial t} +u \frac{\partial u}{\partial x}=\nu  \frac{\partial^2 u}{\partial x^2}
\end{equation} 

The initial conditions were taken as $u(x, 0) = -sin(\pi x)$ with $x \in [-1, 1]$ and the boundary conditions $u(\pm 1, t) = 0$ with $t \in [0,  5]$. It is well know that, as $\nu \to 0$ the solution exhibits steep gradients as time evolves, while as $\nu \to 1$ it becomes smoother. In this model, $\nu$ is considered a random variable following the uniform distribution between $[0, 1]$ to include all possible trends of the solution. In order to obtain exact solutions of eq. \eqref{Burgers1D}, a finite difference scheme is employed in both time and space domains, using a time step of $\Delta t = 0.0505 \ sec$ and spatial discretization $\Delta x = 0.0101 \ m$, leading to 100 and 200 time and spatial points, respectively. 

As explained in the previous section, the first step to apply the proposed surrogate modeling scheme is the generation of a sufficient number of training samples. To this purpose, Burgers' equation is solved numerically for $N = 100$ values of $\nu$ within the range $[0, 1]$. Subsequently, these solution snapshots are stored in a 3-D matrix $\mathbf{S} = [\boldsymbol{U}_{1}, \boldsymbol{U}_{2},...,\boldsymbol{U}_{N}]\in\mathbb{R}^{100 \times 200 \times 100}$, where $\boldsymbol{U}_{i} \in\mathbb{R}^{200 \times 100}$  is the velocity matrix of the $i$-th solution of equation \eqref{Burgers1D}. Then, a CAE is trained over this data set for 2000 epochs with learning rate equal to 1e-4 and a batch size of 16. An adaptive moment optimizer (Adam) \cite{Kingma} is utilized for the loss minimization, with the loss function being the mean square error of eq. \eqref{eq2.1CAE}.

\begin{figure}[H]
    \centering
    \includegraphics[width=1.00\textwidth]{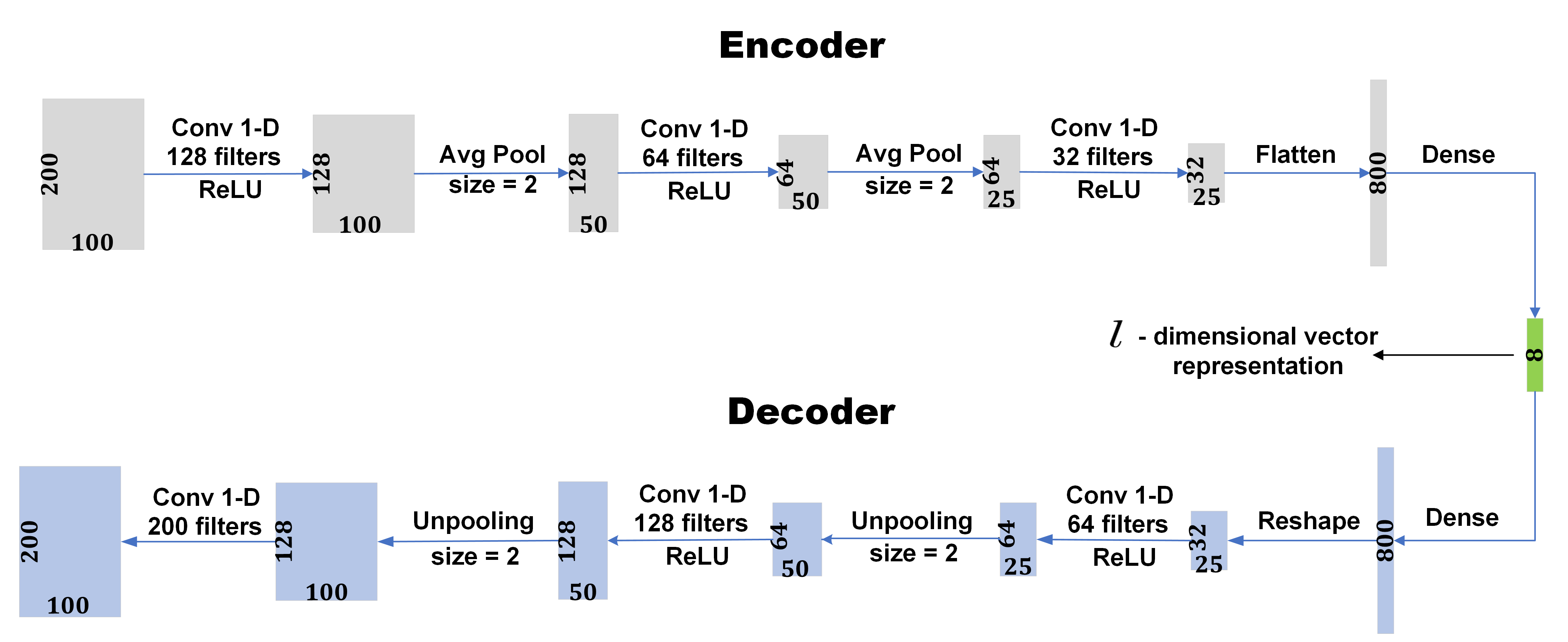}
    \caption{CAE architecture for the solution of Burgers' equation}
    \label{fig:CAE_burgers}
\end{figure}

 After CAE's training phase, an encoded data matrix $\mathbf{S}^{e} = [\boldsymbol{z}_{1}, \boldsymbol{z}_{2},...,\boldsymbol{z}_{N}]$ is obtained, where each column $\boldsymbol{z}_{i}$ is the 8 $\times$ 1 latent vector representation of the solution matrix $\boldsymbol{U}_{i}$. The selected CAE's architecture is presented in figure \ref{fig:CAE_burgers}. The final step of the training procedure is the training of the FFNN in order to establish the mapping from the problem's parameters $\nu_{i}$ to the encoded vector representations $\boldsymbol{z}_{i}$. As shown in table \ref{table:1}, the network's architecture consists of 4 hidden layers with 32 nodes per layer. The ReLU activation function \cite{Nwankpa} is being used in each node, while the Adam optimizer is again utilized to minimize the mean square error loss function. The FFNN was trained for 30000 epochs with learning rate 1e-4. 
 
 \begin{table}[H]
 \centering
\begin{tabular}{ |c|c|c|c| } 
\hline
\textbf{Layer} & \textbf{Nodes} & \textbf{Activation}  \\
\hline
Input & 1 & - \\ 
\hline
Hidden 1 & 32 & ReLU  \\ 
\hline
Hidden 2 & 32 & ReLU  \\ 
\hline
Hidden 3 & 32 & ReLU  \\
\hline
Hidden 4 & 32 & ReLU  \\
\hline
Output & 8 & -  \\ 
\hline
\end{tabular}
\caption{FFNN architecture for the solution of Burgers' equation}
\label{table:1}
\end{table}

\begin{figure}[H]
\centering
\begin{subfigure}[b]{0.48\textwidth}
         \centering
         \includegraphics[width=\textwidth, height = 0.3\textheight]{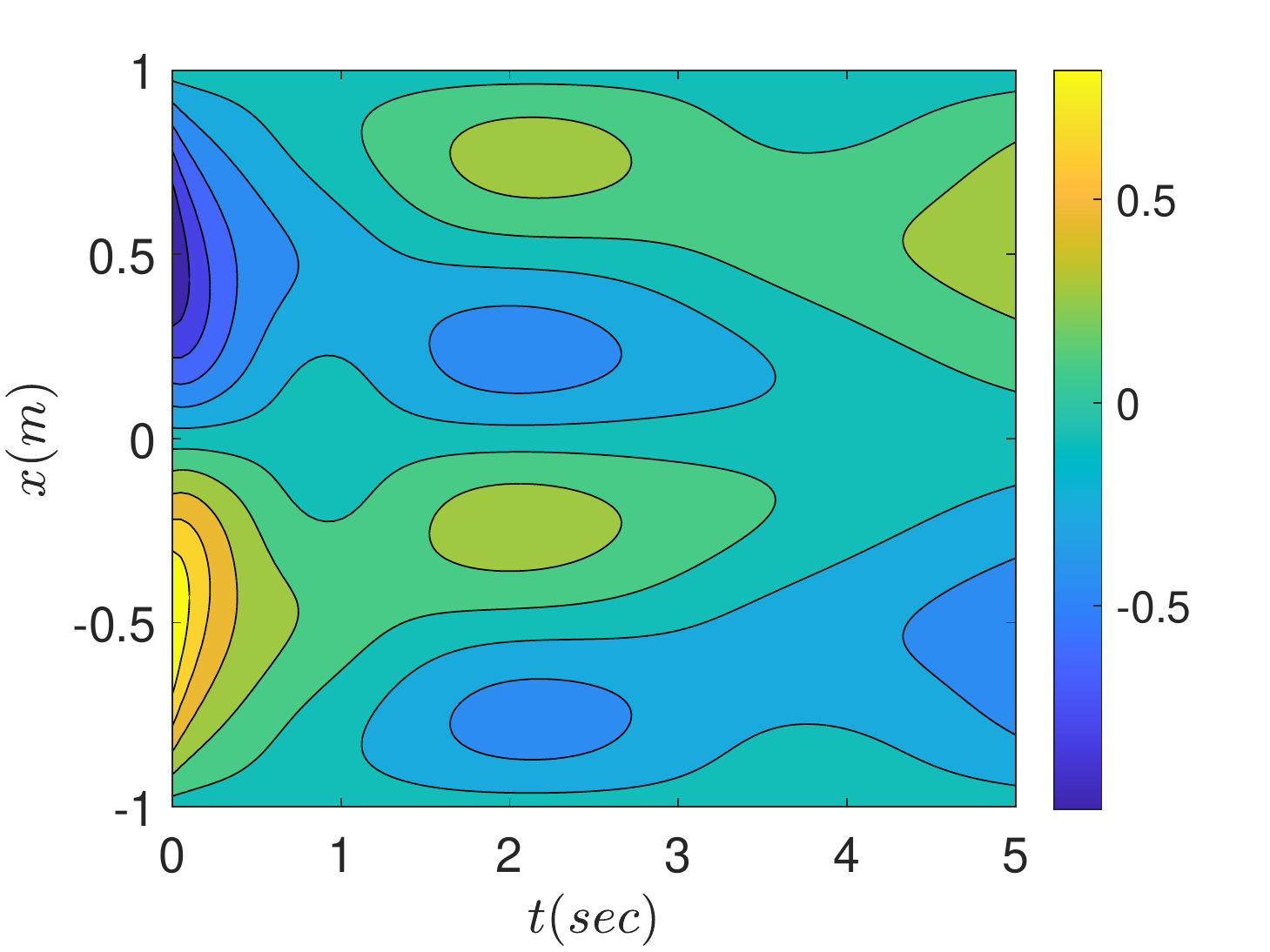}
         \caption{exact model}
     \end{subfigure}
     \hfill
     \begin{subfigure}[b]{0.48\textwidth}
         \centering
         \includegraphics[width=\textwidth, height = 0.3\textheight]{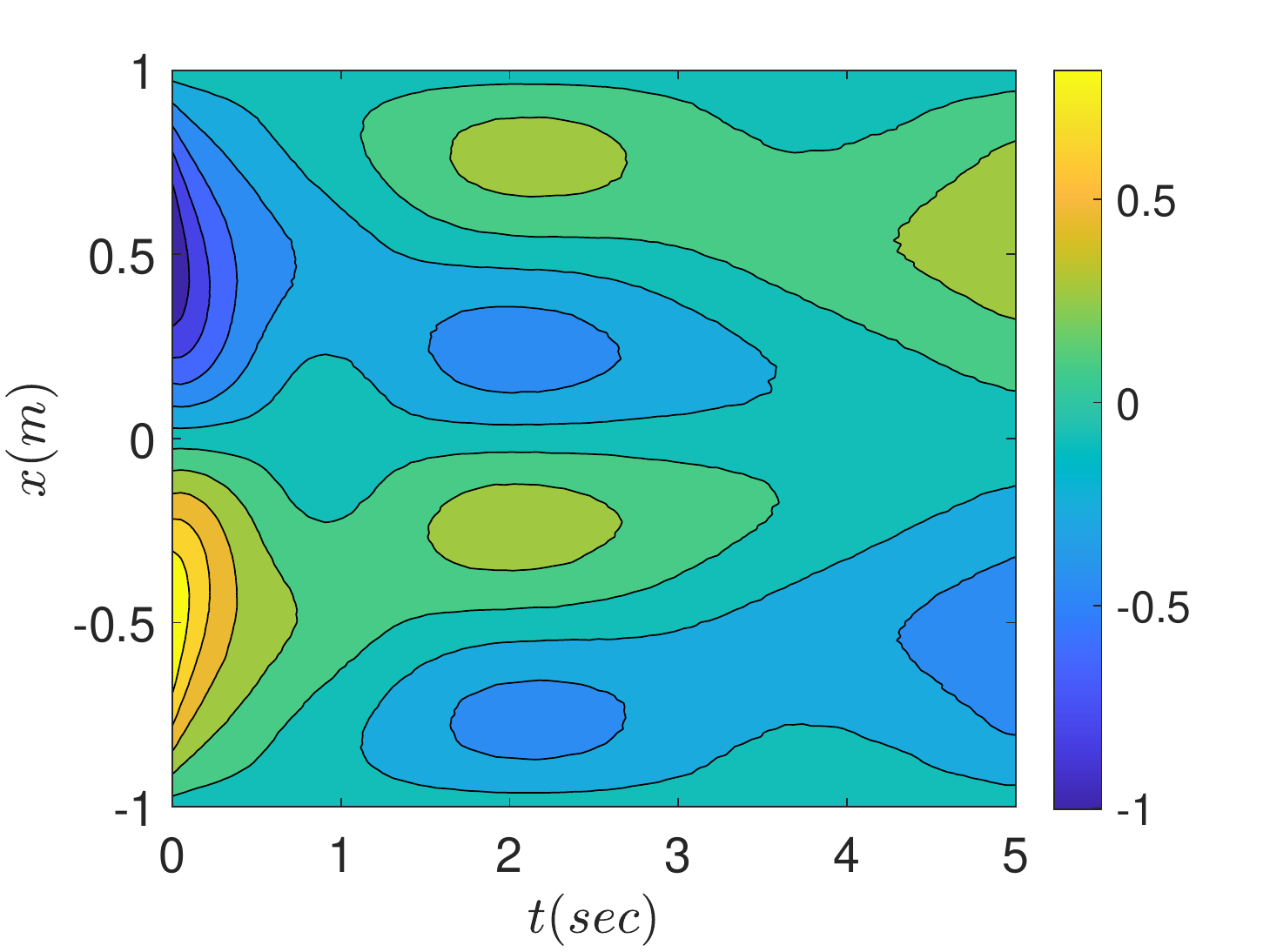}
         \caption{surrogate model}
     \end{subfigure}  
    \caption{Solution profile $u(x,t)$ for $\nu = 0.2$ predicted by (a) the exact model and (b) the surrogate model}
    \label{fig:solution02}
\end{figure}

\begin{figure}[H]
    \centering
    \includegraphics[width=0.82\textwidth]{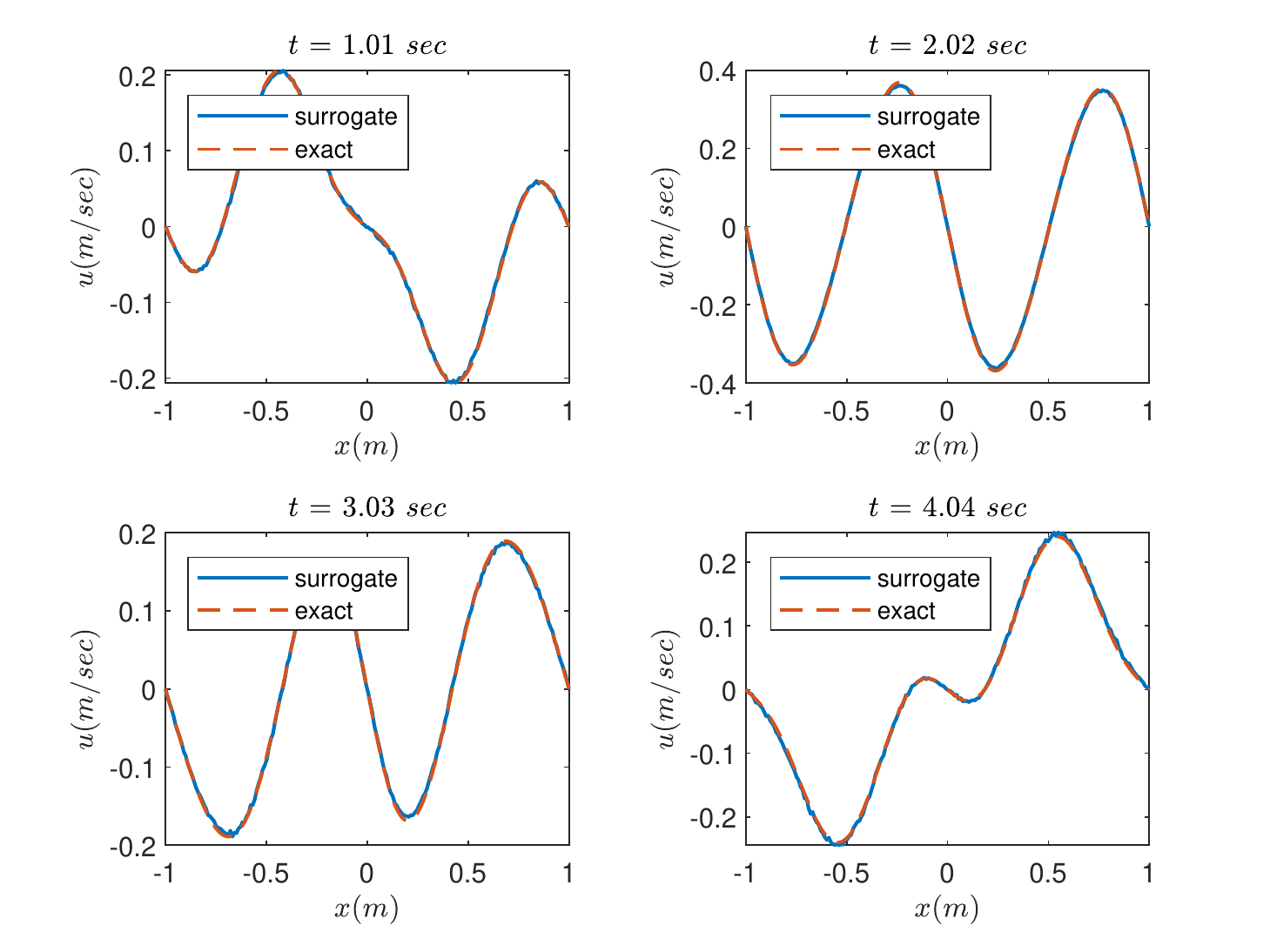}
    \caption{Solution profiles $u(x,t)$ at specific time instants for $\nu = 0.2$}
    \label{fig:solutionProfiles0.2}
\end{figure}

\begin{figure}[H]
\centering
\begin{subfigure}[b]{0.48\textwidth}
         \centering
         \includegraphics[width=\textwidth, height = 0.3\textheight]{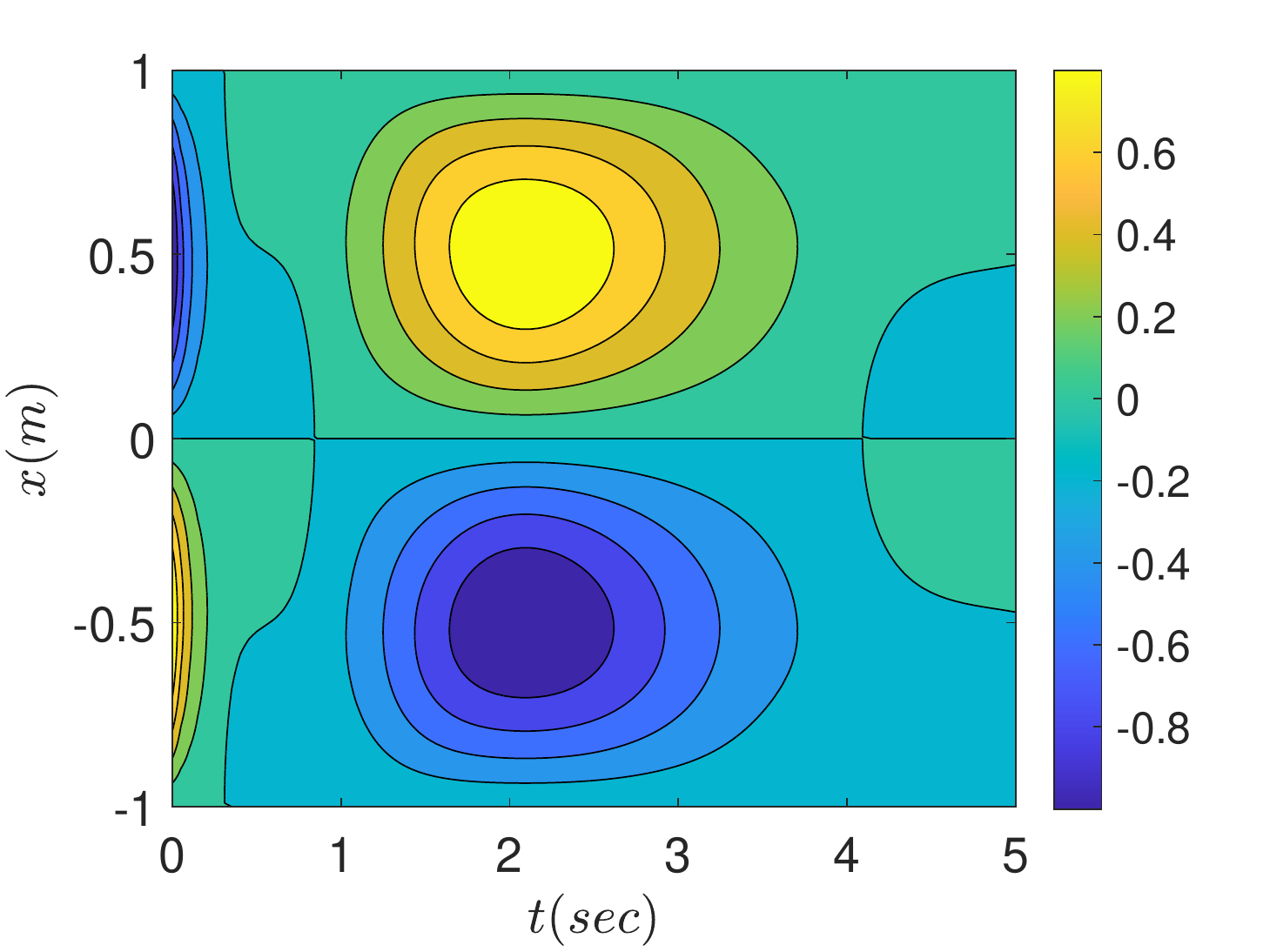}
         \caption{exact model}
     \end{subfigure}
     \hfill
     \begin{subfigure}[b]{0.48\textwidth}
         \centering
         \includegraphics[width=\textwidth, height = 0.3\textheight]{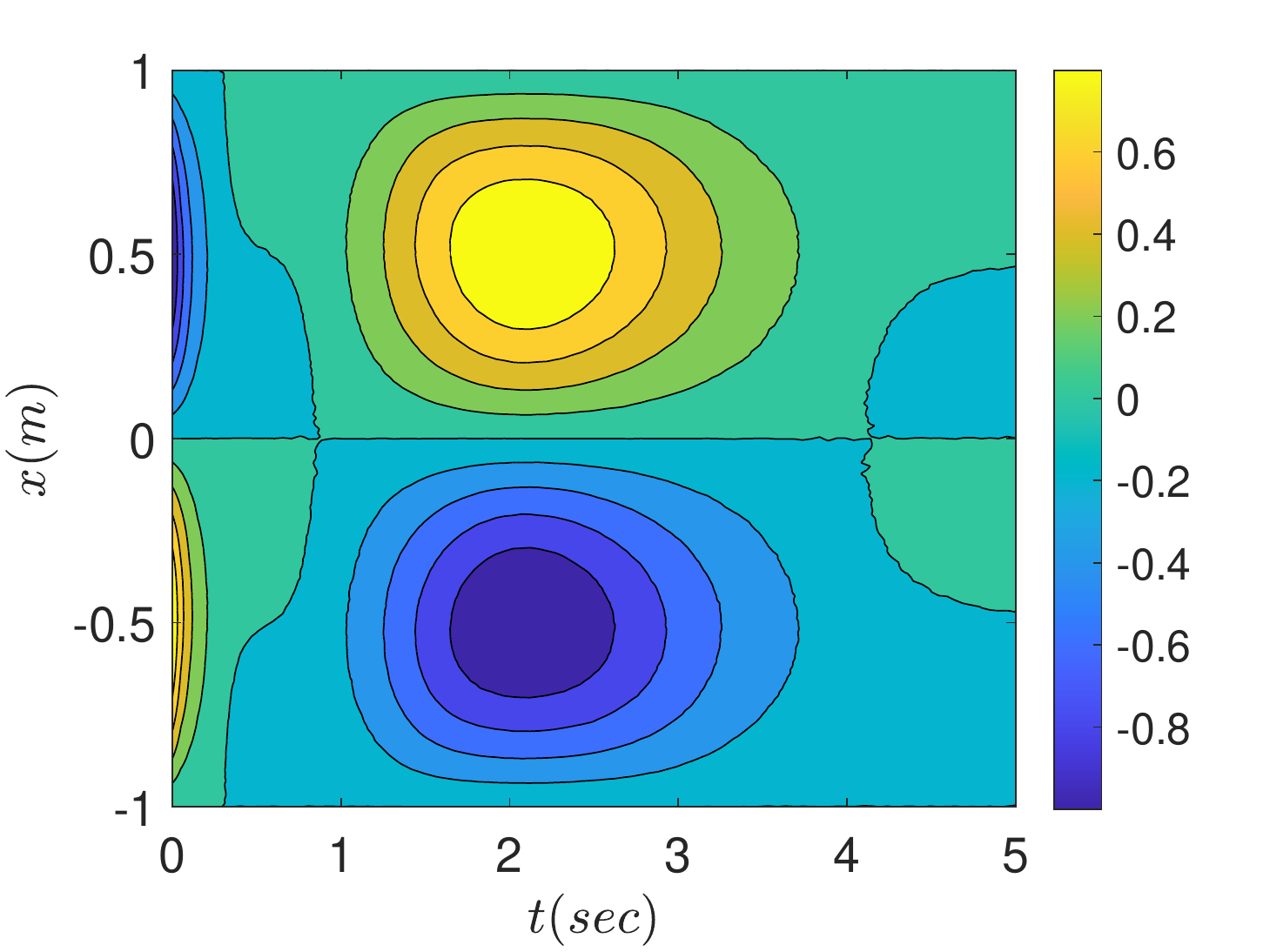}
         \caption{surrogate model}
     \end{subfigure}  
    \caption{Solution profile $u(x,t)$ for $\nu = 0.8$ predicted by (a) the exact model and (b) the surrogate model}
    \label{fig:solution08}
\end{figure}

\begin{figure}[H]
    \centering
    \includegraphics[width=0.82\textwidth]{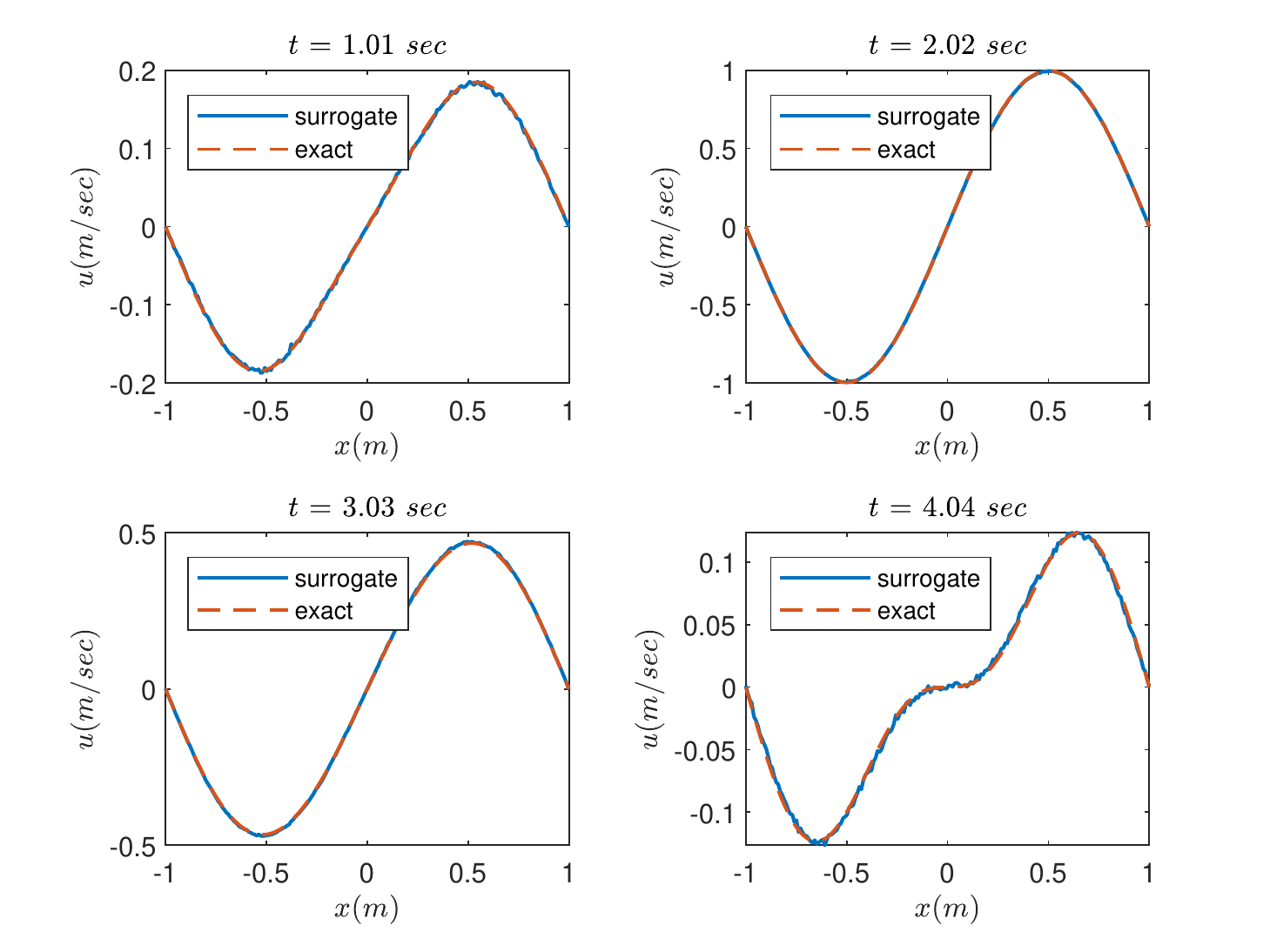}
    \caption{Solution profiles $u(x,t)$ at specific time instants for $\nu = 0.8$}
    \label{fig:solutionProfiles0.8}
\end{figure}

 The CAE-FFNN model accuracy is tested on the solutions for the values of $\nu=0.2, \ 0.8$ that were not included in the initial training data set and compared with those predicted by the finite differences model. Figures \ref{fig:solution02} and \ref{fig:solution08} present the total solution field, while figures \ref{fig:solutionProfiles0.2} and \ref{fig:solutionProfiles0.8} illustrate the solution profiles for specific time steps. From these results it can be seen that the predictions of the proposed surrogate model are almost identical to those of the exact solution. The normalized error between the solution matrices $\boldsymbol{U}_{FD}$ and $\boldsymbol{U}_{SUR}$ of the finite differences model and the surrogate model, respectively, given by $\widehat{err}=\Vert\boldsymbol{U}_{FD}-\boldsymbol{U}_{SUR}\Vert_2/ \Vert\boldsymbol{U}_{FD}\Vert_2$, with $\| \cdot \|_2$ being the $L_{2}$ matrix norm, was found equal to $1.23\%$ for the case of $\nu = 0.2$ and $0.53\%$ for $\nu = 0.8$. 
 
\begin{figure}[H]
\centering
\begin{subfigure}[b]{0.48\textwidth}
         \centering
         \includegraphics[width=\textwidth, height = 0.3\textheight]{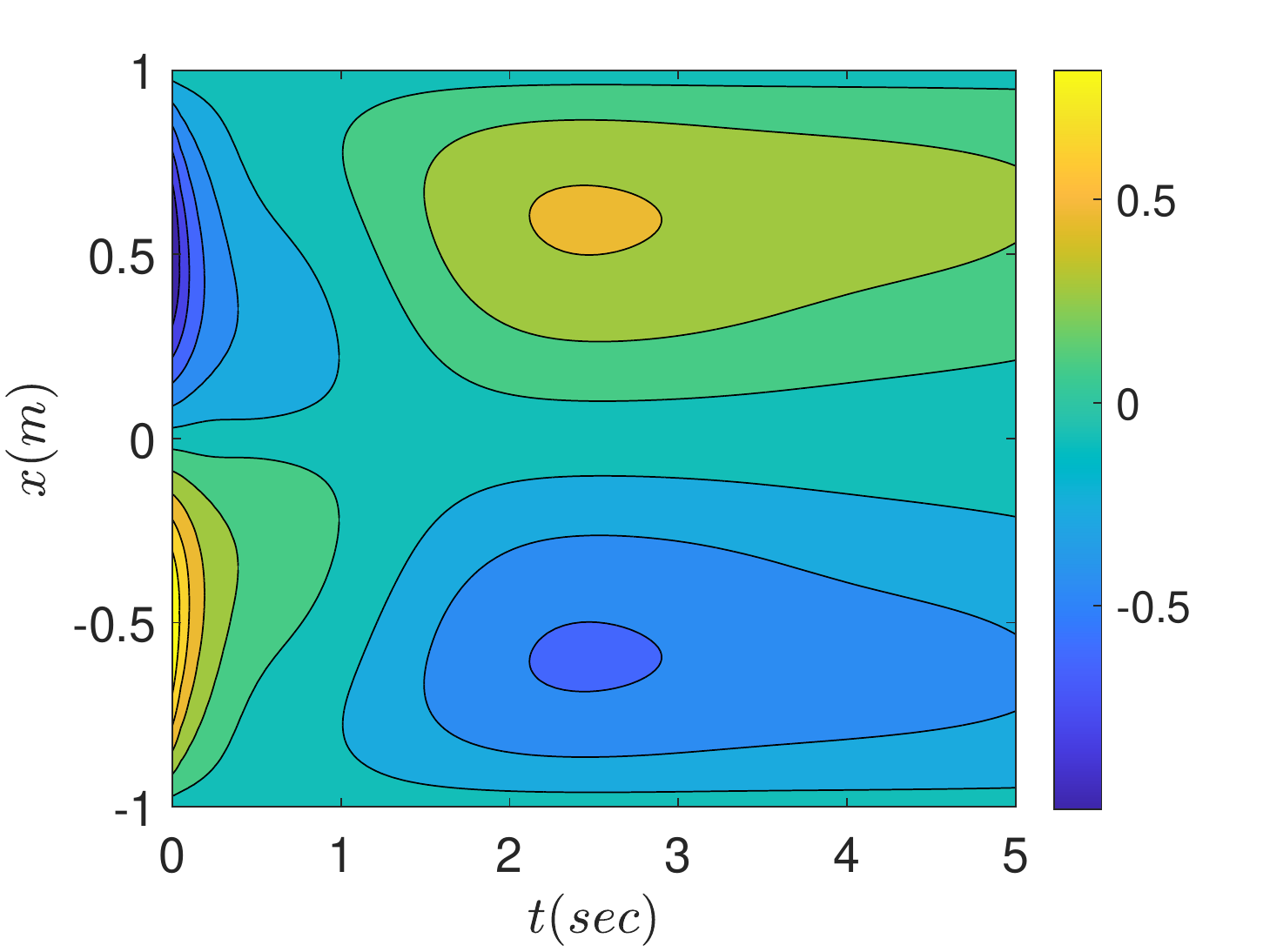}
         \caption{exact model}
     \end{subfigure}
     \hfill
     \begin{subfigure}[b]{0.48\textwidth}
         \centering
         \includegraphics[width=\textwidth, height = 0.3\textheight]{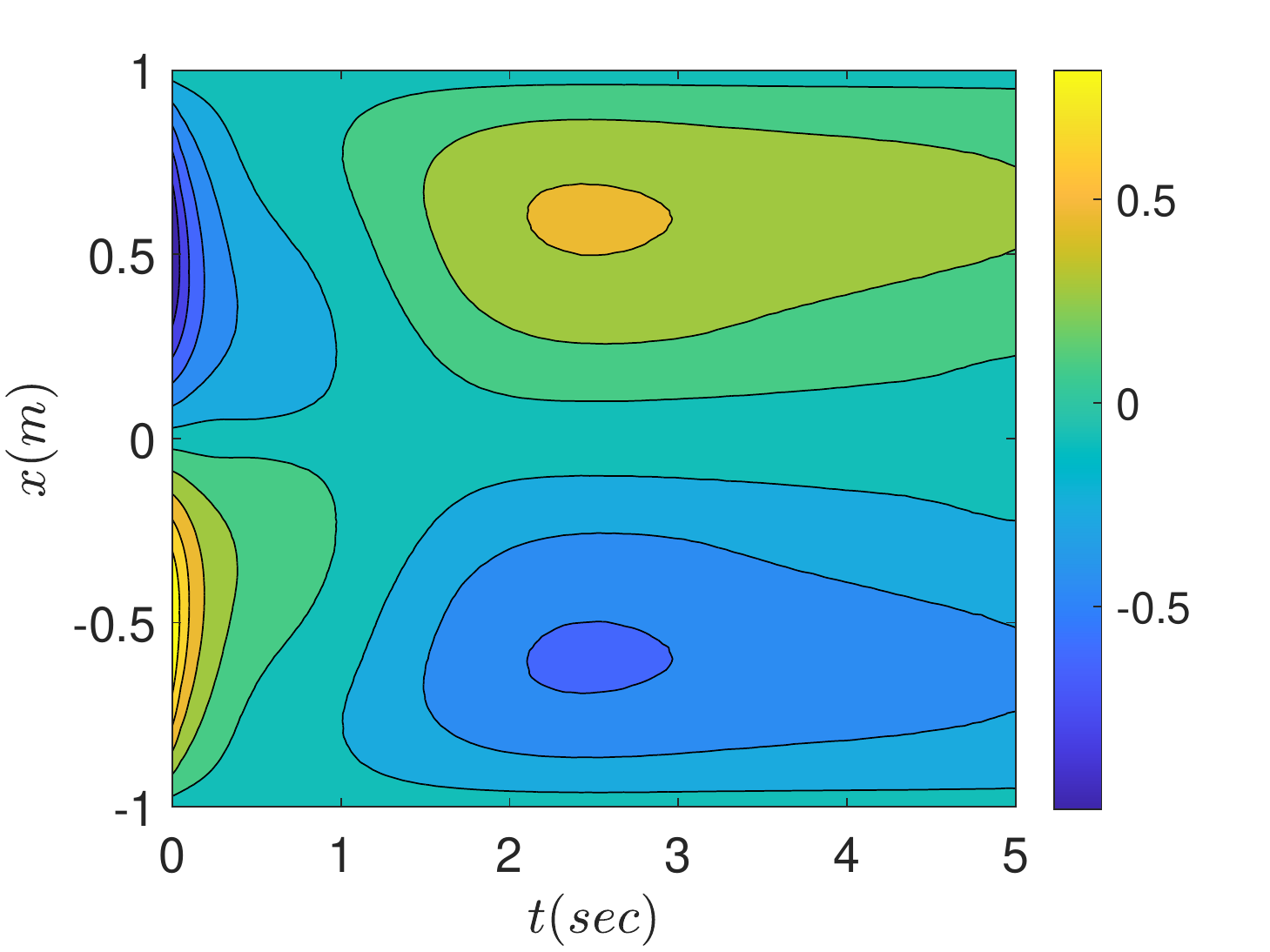}
         \caption{surrogate model}
     \end{subfigure}  
    \caption{Mean value of $u(x,t)$ predicted by (a) the exact model and (b) the surrogate model}
    \label{fig:MeanBurgers}
\end{figure}

\begin{figure}[H]
    \centering
    \includegraphics[width=0.82\textwidth]{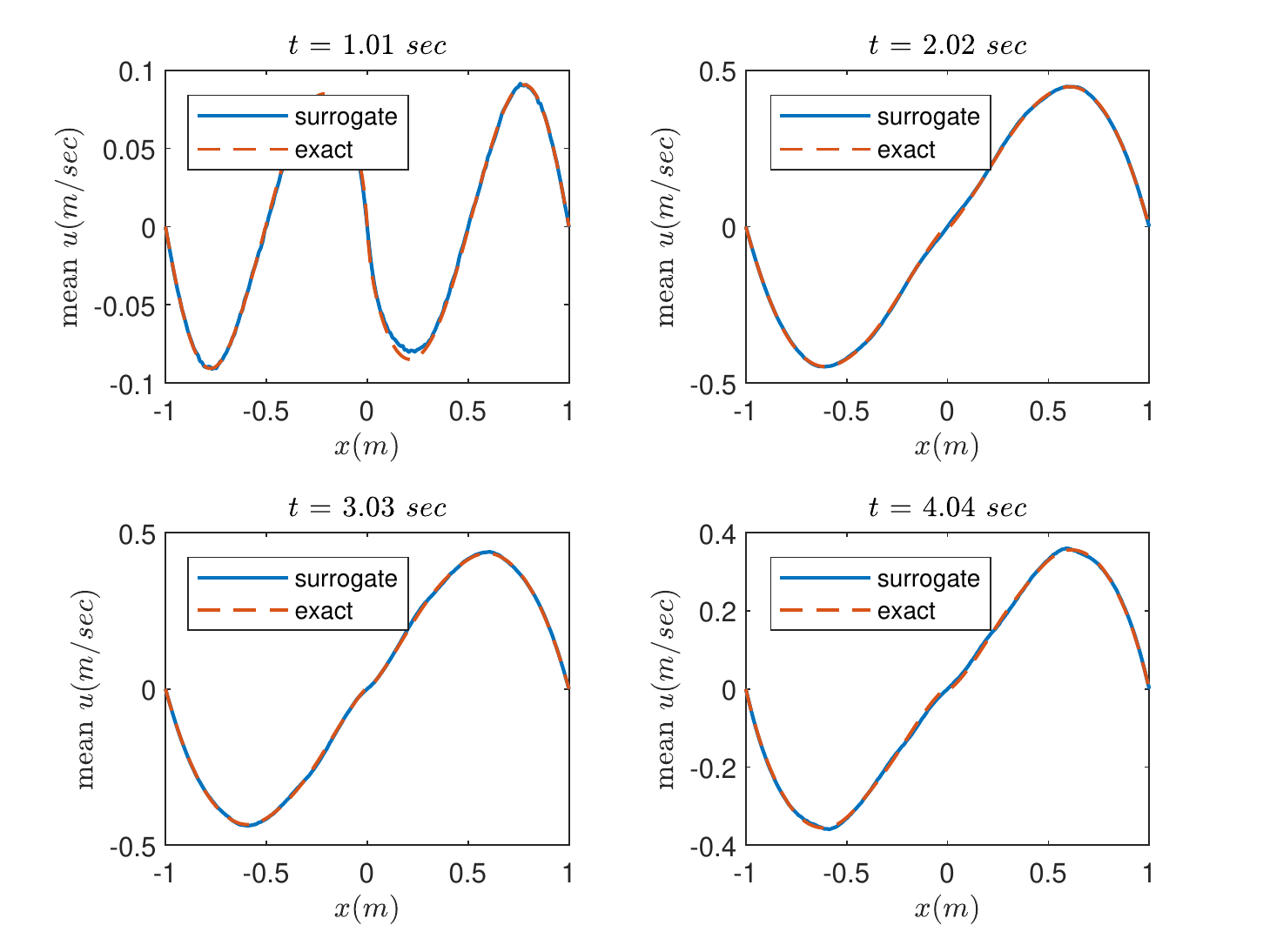}
    \caption{Mean value of $u(x,t)$ at specific time instants}
    \label{fig:ProfilesMeanBurgers}
\end{figure}

\begin{figure}[H]
\centering
\begin{subfigure}[b]{0.48\textwidth}
         \centering
         \includegraphics[width=\textwidth, height = 0.3\textheight]{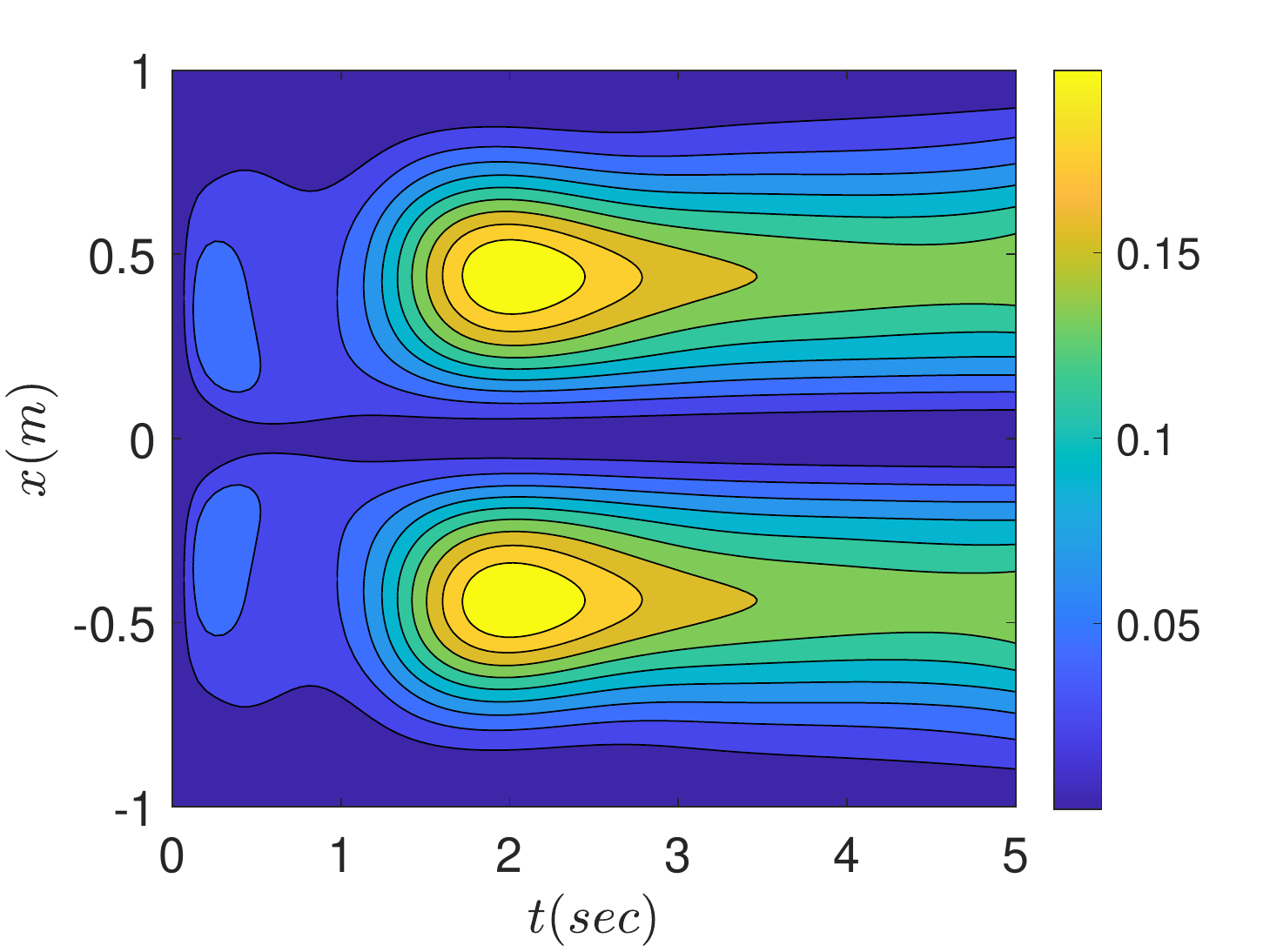}
         \caption{exact model}
     \end{subfigure}
     \hfill
     \begin{subfigure}[b]{0.48\textwidth}
         \centering
         \includegraphics[width=\textwidth, height = 0.3\textheight]{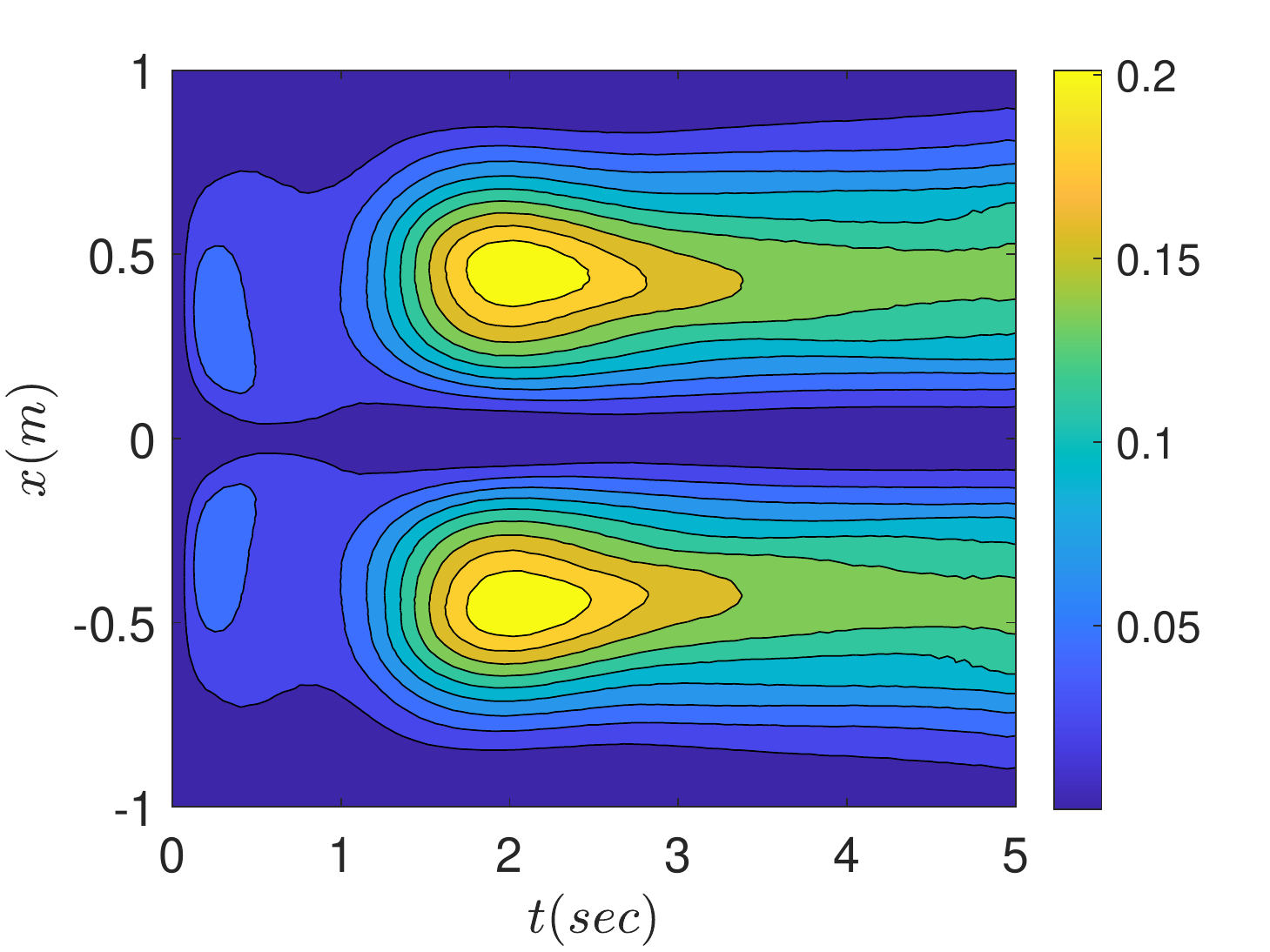}
         \caption{surrogate model}
     \end{subfigure}  
    \caption{Variance of $u(x,t)$ predicted by (a) the exact model and (b) the surrogate model}
    \label{fig:VarBurgers}
\end{figure}

\begin{figure}[H]
    \centering
    \includegraphics[width=0.82\textwidth]{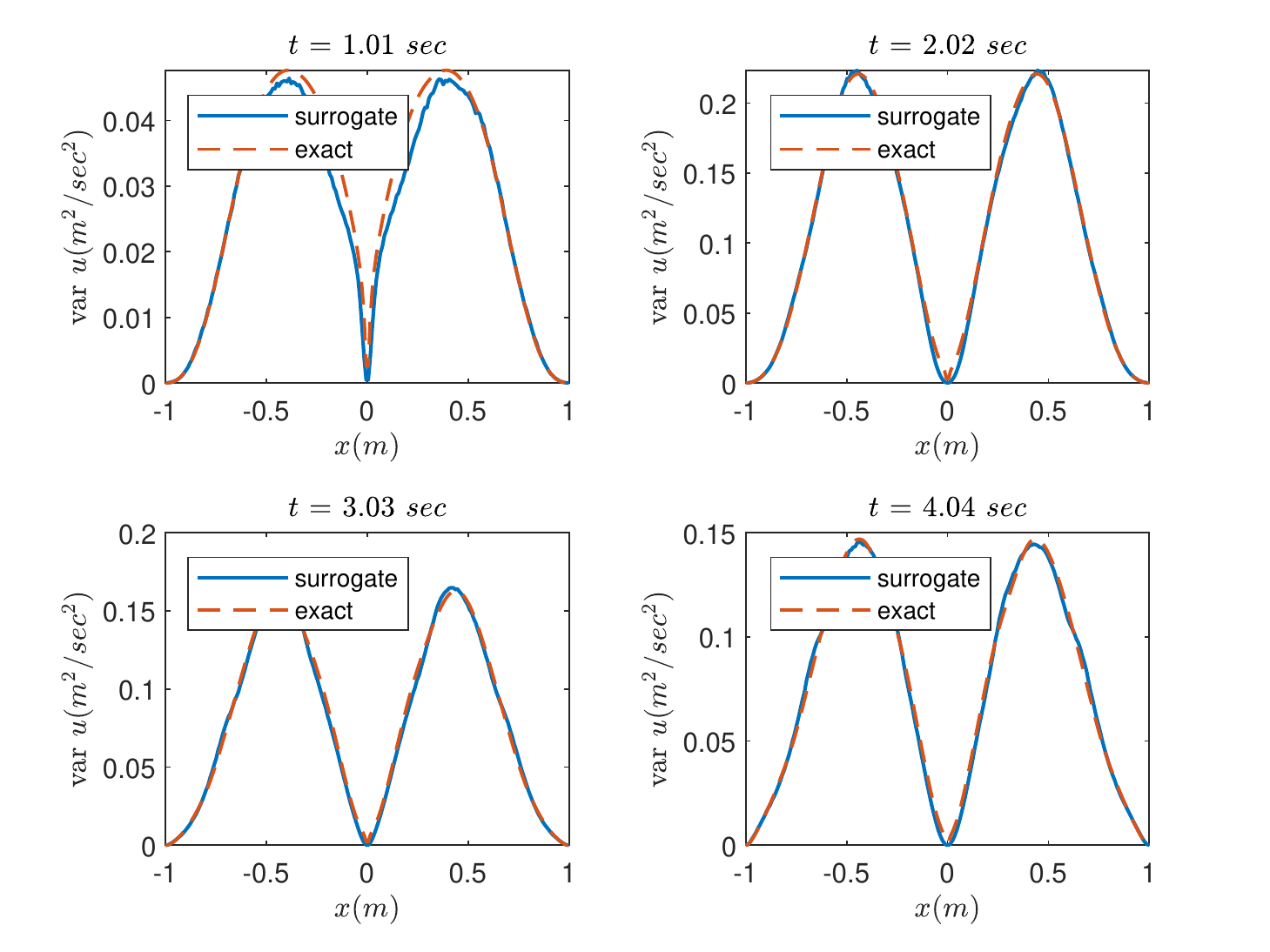}
    \caption{Variance of $u(x,t)$ at specific time instants}
    \label{fig:ProfilesVarBurgers}
\end{figure}

Subsequently, in the context of the MC analysis, $N_{MC} = 3000$ values of $\nu \sim \mathcal{U}[0,1]$ are generated according to their distribution and the corresponding PDEs are solved by the exact and the surrogate model, respectively. The mean value and variance of $u(x,t)$ obtained by the two models are depicted in figures \ref{fig:MeanBurgers} and \ref{fig:VarBurgers}, respectively, while figures \ref{fig:ProfilesMeanBurgers} and \ref{fig:ProfilesVarBurgers} present a comparison between the two models in the mean value and the variance of $u(x,t)$ at specific time instants. As evidenced by these results, the surrogate and the exact model are in very close agreement. The normalized error between the mean solution matrices $\boldsymbol{M}_{FD}$ and $\boldsymbol{M}_{SUR}$ was found equal to 0.96\%, while the same error for the variance matrices $\boldsymbol{V}_{FD}$ and $\boldsymbol{V}_{SUR}$ was 2.54\%.

\begin{figure}[H]
\centering
\begin{subfigure}[b]{0.48\textwidth}
         \centering
         \includegraphics[width=\textwidth]{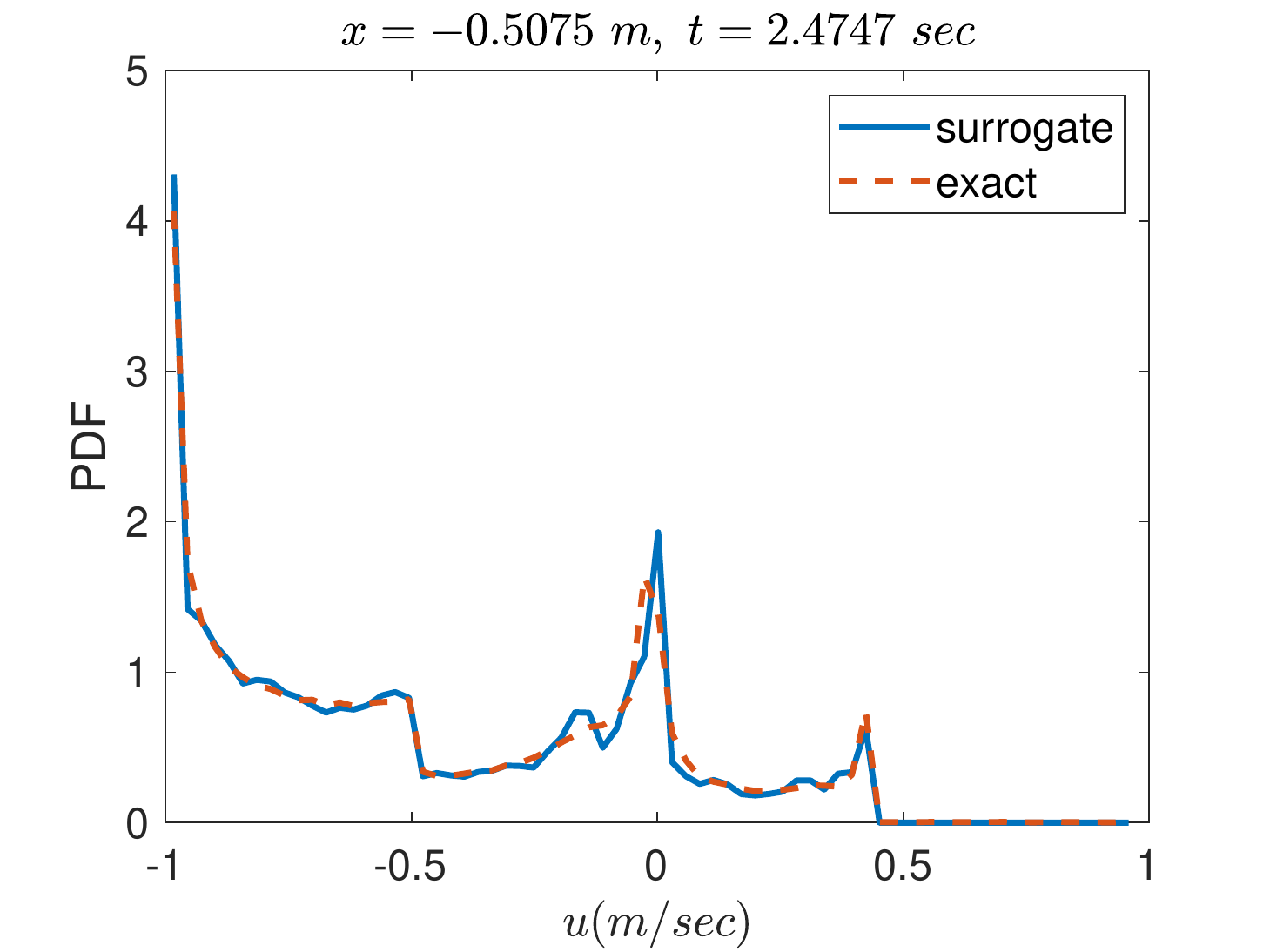}
         \caption{}
     \end{subfigure}
     \hfill
     \begin{subfigure}[b]{0.48\textwidth}
         \centering
         \includegraphics[width=\textwidth]{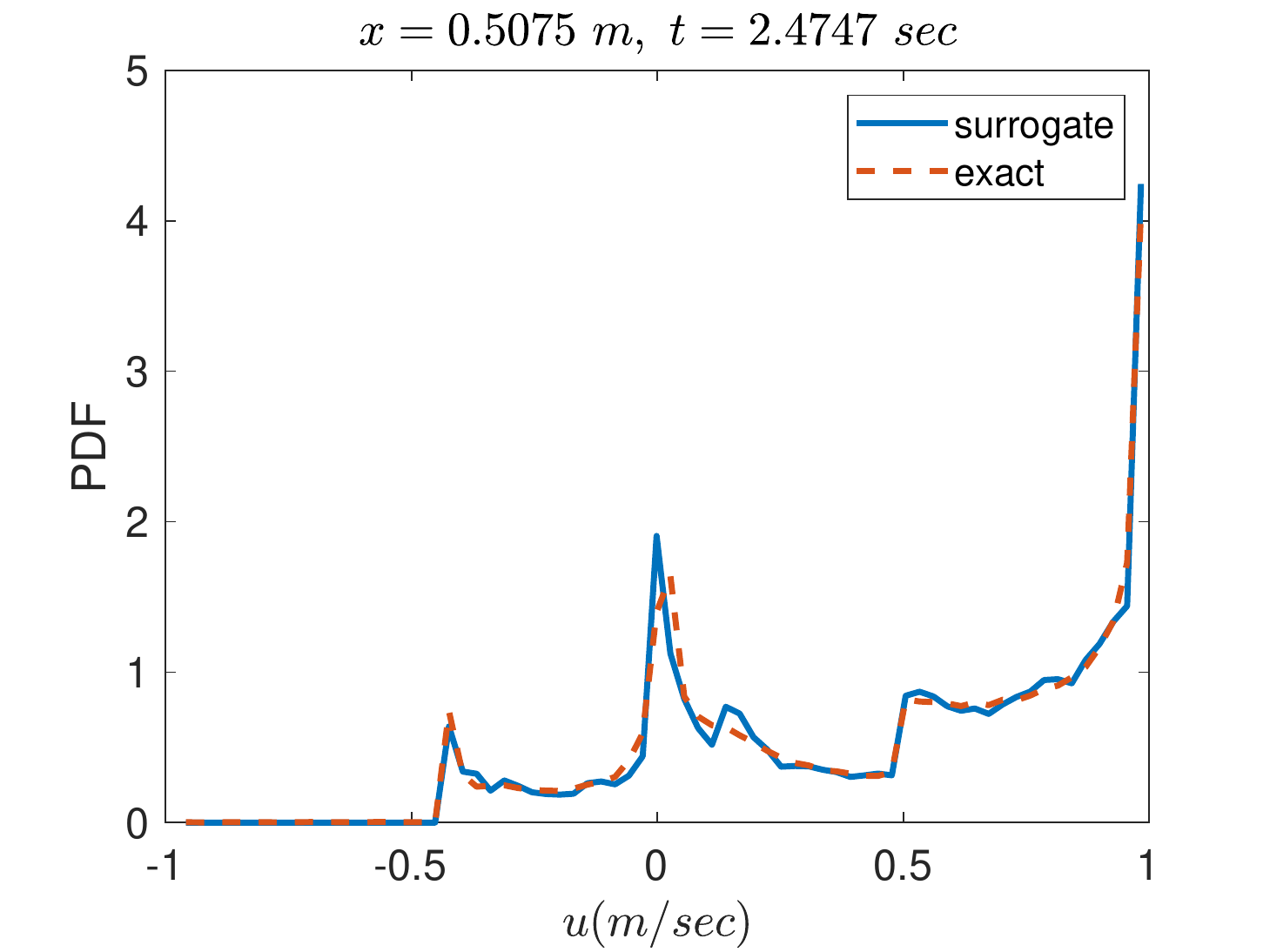}
         \caption{}
     \end{subfigure}  
    \caption{PDF of $u(x,t)$ predicted by the exact model and the surrogate model}
    \label{fig:pdf_burgers}
\end{figure}

To further test the method's interpolation capabilities, a significantly larger number of MC simulations are  performed, in order to acquire the probability density function (PDF) of $u(x, t)$. Specifically, $N_{MC} = 300000$ simulations are carried out by the two models and the results pertaining to positions $x = - 0.5075 \ m$ and $x = 0.5075 \ m$ at $t = 2.4747  \ sec$ are depicted in figure \ref{fig:pdf_burgers}. It becomes apparent from this figure that the surrogate model is able to predict the PDF of $u(x, t)$ with satisfactory accuracy.

\begin{figure}[H]
    \centering
    \includegraphics[width=0.70\textwidth]{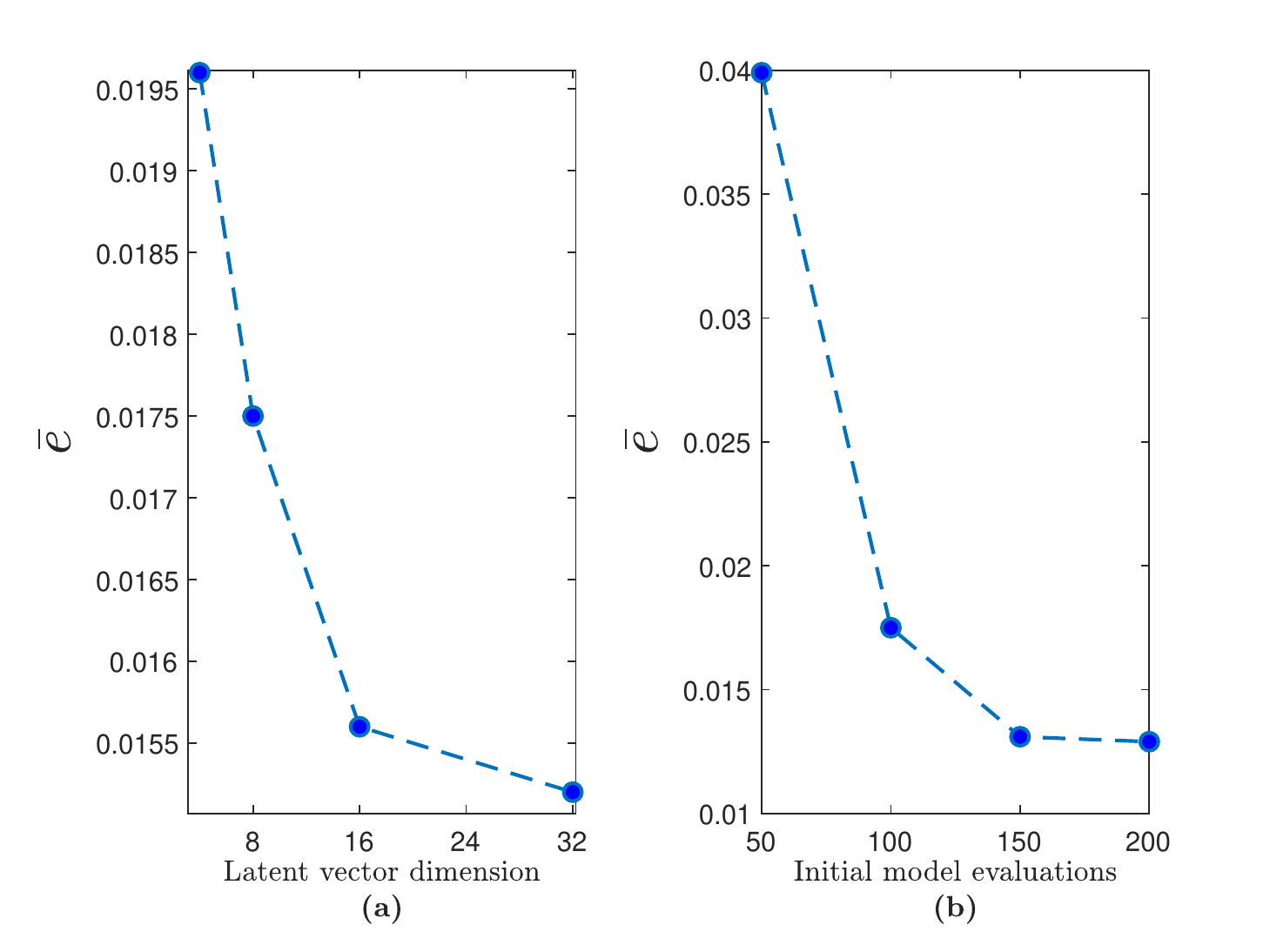}
    \caption{mean error $\Bar{e}$ with respect to (a) the latent vector dimension and (b) the initial data set size}
    \label{fig:convergence}
\end{figure}

Finally, a convergence study with respect to the dimension of the latent vectors and the size of the initial data set is presented in figure \ref{fig:convergence}. The average normalized error is defined as:

\begin{equation} \label{errorNorm}
\Bar{e} = \frac{1}{N_{MC}} \sum_{j=1}^{N_{MC}}\frac{||\boldsymbol{U}_{FD}^{j} - \boldsymbol{U}_{SUR}^{j}||}{||\boldsymbol{U}_{FD}^{j}||}
\end{equation}

\noindent where $\boldsymbol{U}_{FD}^{j}$ and $\boldsymbol{U}_{SUR}^{j}$ are the solution matrices of the $j$-th MC simulation obtained by the 'exact' model and the surrogate model, respectively. From these results, it becomes apparent that a selection of a higher dimensional latent vector representation reduces the amount of information lost in the decoding process, thus is linked to improved accuracy. Subsequently, as the initial data set size increases a reduced mean error $\Bar{e}$ is achieved and converges close to the value $\Bar{e}_{lim} \approx 0.012$.

\subsection{Coupled Shear Walls under Seismic Loading}
\begin{figure}[H]
    \centering
    \includegraphics[width=0.64\textwidth]{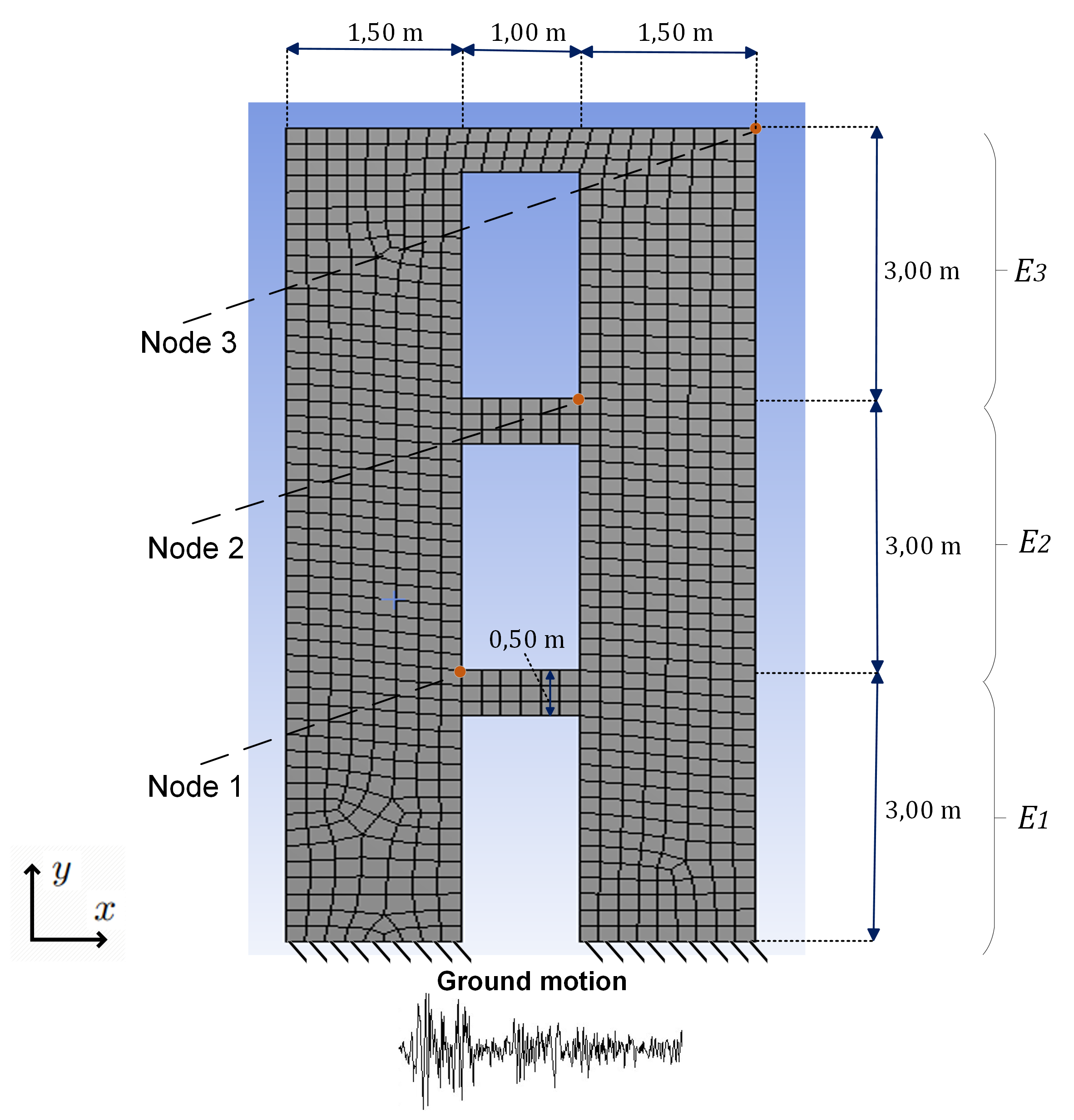}
    \caption{Geometry and finite element meshing of the coupled shear walls}
    \label{fig:shear_walls}
\end{figure}

\begin{figure}[H]
    \centering
    \includegraphics[width=0.50\textwidth]{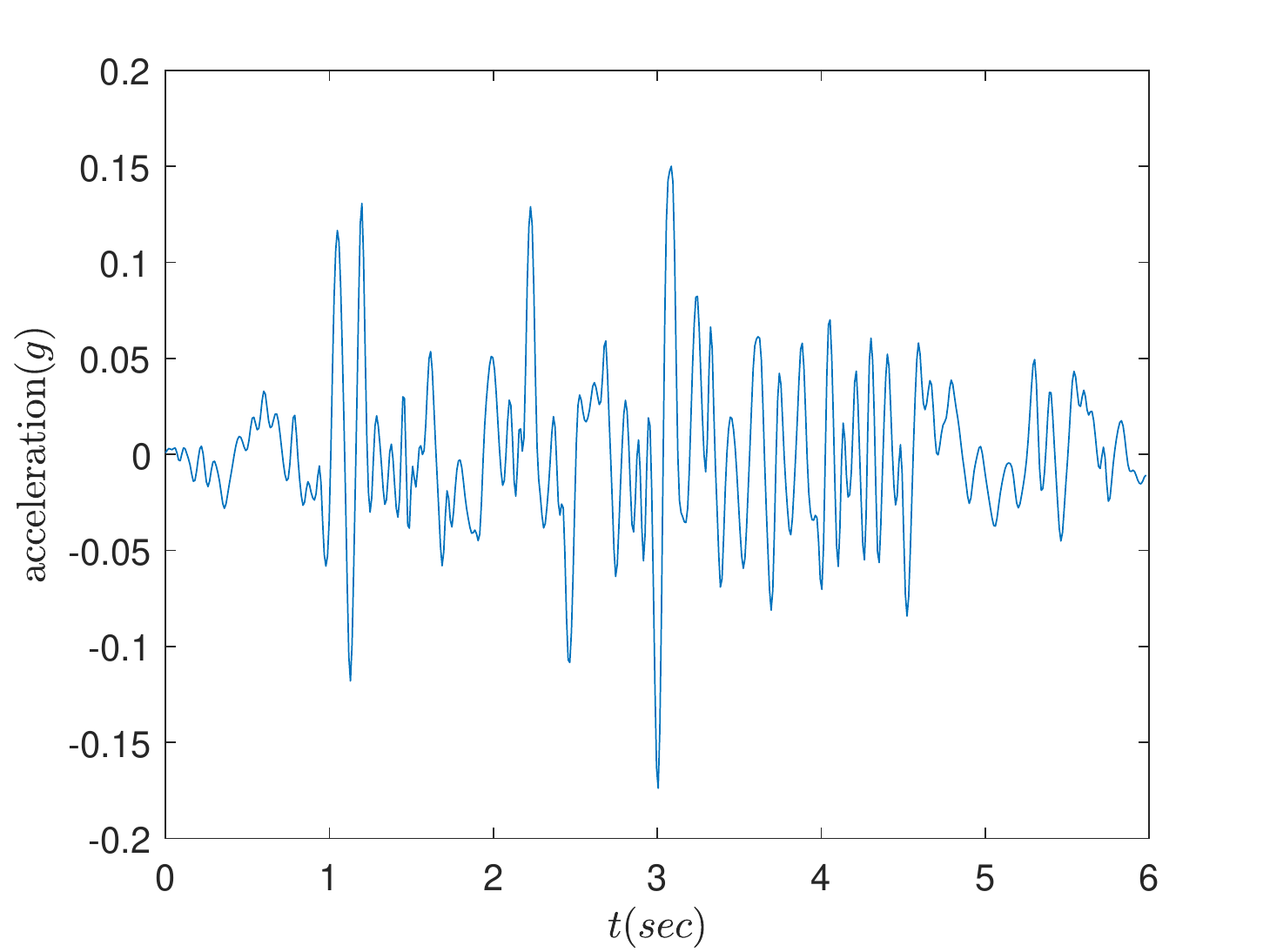}
    \caption{Acceleration data of the selected ground motion}
    \label{fig:accel}
\end{figure}

A transient plane stress structural problem is considered as the second test example. The problem is governed by the equations of motion of 2-D linear elasticity:

\begin{align} 
\frac{E}{2(1+\nu)}\nabla^2 u_{x} + \frac{E}{2(1-\nu)} \frac{\partial}{\partial x}\bigg(\frac{\partial u_{x}}{\partial x} + \frac{\partial u_{y}}{\partial y}\bigg) + p_{x} = \rho  \frac{\partial{u_{x}}^2}{\partial{t}^2}  \label{planeStressProblem1}\\ 
\frac{E}{2(1+\nu)}\nabla^2 u_{y} + \frac{E}{2(1-\nu)} \frac{\partial}{\partial y}\bigg(\frac{\partial u_{x}}{\partial x} + \frac{\partial u_{y}}{\partial y}\bigg) + p_{y} = \rho  \frac{\partial{u_{y}}^2}{\partial{t}^2} \label{planeStressProblem2}
\end{align}

\noindent where $u_{x} \equiv u_{x}(x,y,t)$ and $u_{y}\equiv u_{y}(x,y,t)$ are the displacement fields, $E$ is the modulus of elasticity, $\nu$ is the Poisson ratio, $\rho$ is the material's mass density and $p_{x}$ and $p_{y}$ are the body forces. Specifically, the three-story reinforced concrete coupled shear walls of figure \ref{fig:shear_walls} were subjected to seismic loading, that of the  accelerogram of 1972 Kefalonia earthquake \cite{Kefalonia} (figure \ref{fig:accel}) with a total duration of 6.00 $sec$. The Poisson ratio is assumed $\nu = 0.2$, the mass density of the wall is taken as $\rho = 2500$ $kg/m^{3}$, the thickness of the wall is considered $\tau = 1$ $m$, while body forces $p_{x}$ and $p_{y}$ were assumed zero. The Young moduli $E_1$, $E_2$ and $E_3$ of each story are considered as uncorrelated random variables following the log-normal distribution with mean value $\mu = 30$ $GPa$ and standard deviation $\sigma = 0.25\mu = 7.5$ $GPa$. This phenomenon occurs in reinforced concrete structures, where the construction of each story is initiated several days after the completion of the previous story so that the concrete achieves at least 95\% of its design strength capacity. As a consequence, the concrete mixture used in each construction phase is different, which justifies the lack of correlation in the random variables describing the mechanical properties of each storey. The selection of the log-normal distribution with such a high value for the standard deviation $\sigma$ is purely for academic purposes in order to illustrate the capabilities of the proposed method.

The 'exact' solutions of the problem are obtained by solving eq. \eqref{planeStressProblem1} and \eqref{planeStressProblem2} with the finite element method using plane stress elements. More specifically, the walls are spatially discretized with 876 quadrilateral elements, while for the time discretization the Newmark integration scheme \cite{NOH2019106079} is applied with time step size $\Delta t = 0.01$ $sec$, leading to a total of 1966 degrees of freedom and 600 time steps for the spatial and time domain, respectively.

In order to efficiently span the parametric space, $N = 500$ triplets of parameters $\{[E_{1}^i, E_{2}^i, E_{3}^i]\}_{i=1}^{N}$ are generated with the aid of Latin Hybercube Sampling (LHS) \cite{OLSSON200347}. For each triplet of parameters, the corresponding dynamic problem is solved with the above mentioned numerical procedure and the solution matrices are stored in a 3D matrix $\mathbf{S} = [\boldsymbol{U}_{1}, \boldsymbol{U}_{2},....,\boldsymbol{U}_{N}] \in\mathbb{R}^{500 \times 1966 \times 600}$. 

\begin{figure}[H]
    \centering
    \includegraphics[width=1.0\textwidth]{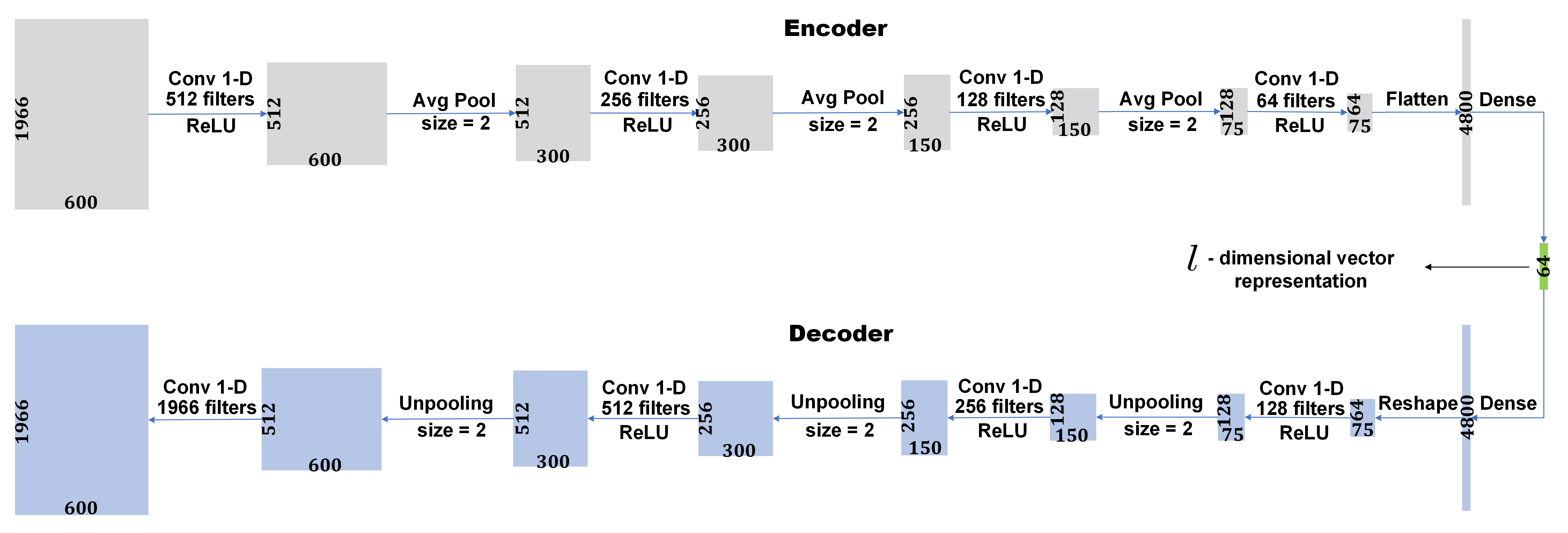}
    \caption{CAE architecture for the solution of the structural dynamic problem}
    \label{fig:CAE_structural}
\end{figure}

A CAE is subsequently trained over this data set for 500 epochs with learning rate 1e-4 and a batch size of 8. The mean square reconstruction error of $\boldsymbol{U}_{i}$ is minimized again by the Adam optimizer. The selected CAE's architecture is presented in figure \ref{fig:CAE_structural}.  An encoded $500 \times 64$ training matrix $\mathbf{S}^{e} = [\boldsymbol{z}_{1}, \boldsymbol{z}_{2},...,\boldsymbol{z}_{N}]$ is then obtained via the encoder, where each column $\boldsymbol{z}_{i}$ is the $64 \times 1$ latent vector representation of the solution time history matrix $\boldsymbol{U}_{i}$.  The above encoded training matrix $\mathbf{S}^{e}$ along with the stored parameter triplets $\{[E_{1}^i, E_{2}^i, E_{3}^i]\}_{i=1}^{N}$ from the previous step are used as outputs and inputs, respectively, in the training process of the FFNN in order to construct a mapping from the parametric to the encoded solution space. The FFNN is trained  for 10000 epochs with learning rate 1e-4 and a batch size of 100. The selected  architecture is shown in table \ref{table:2}.

 \begin{table}[H]
 \centering
\begin{tabular}{ |c|c|c|c| } 
\hline
\textbf{Layer} & \textbf{Nodes} & \textbf{Activation} \\
\hline
Input & 3 & - \\ 
\hline
Hidden 1 & 256 & ReLU \\ 
\hline
Hidden 2 & 256 & ReLU \\ 
\hline
Hidden 3 & 256 & ReLU \\ 
\hline
Hidden 4 & 256 & ReLU \\ 
\hline
Hidden 5 & 256 & ReLU \\ 
\hline
Hidden 6 & 256 & ReLU\\ 
\hline
Output & 64 & -  \\ 
\hline
\end{tabular}
\caption{FFNN architecture for the solution of the structural dynamic problem}
\label{table:2}
\end{table}

\begin{figure}[H]
\centering
\begin{subfigure}[b]{0.48\textwidth}
         \centering
         \includegraphics[width=0.80\textwidth]{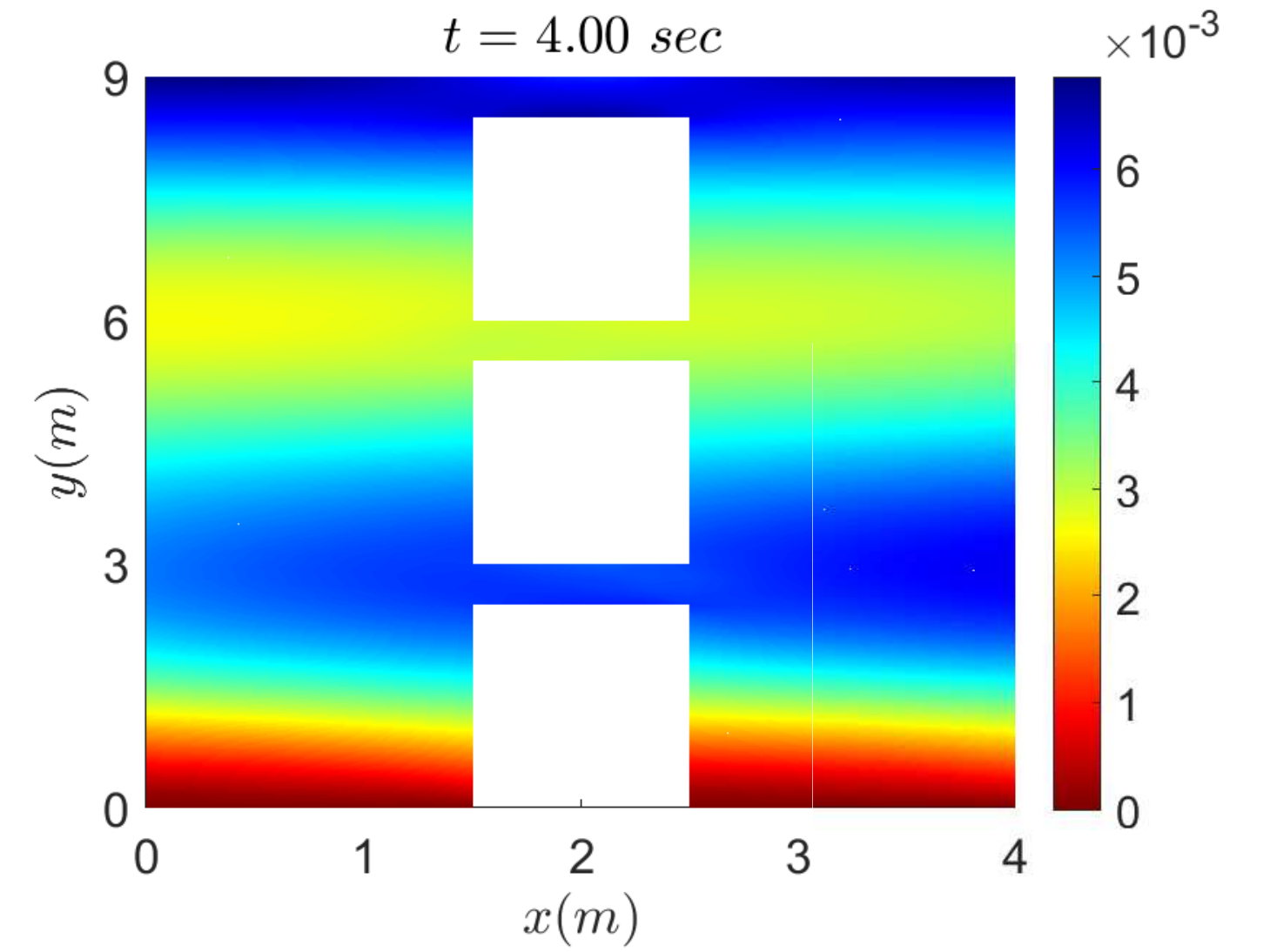}
         \caption{exact model}
     \end{subfigure}
     \hfill
     \begin{subfigure}[b]{0.48\textwidth}
         \centering
         \includegraphics[width=0.80\textwidth]{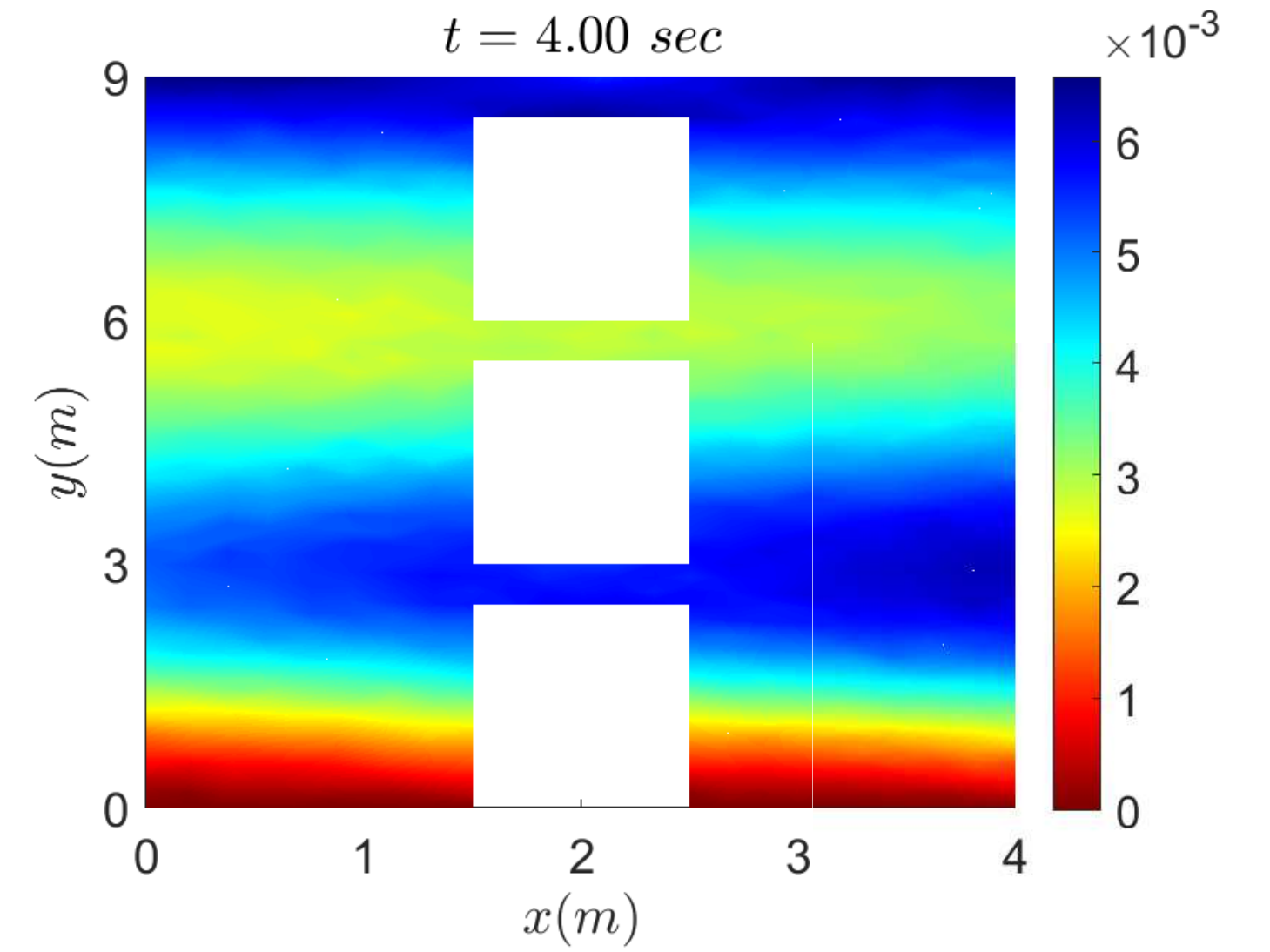}
         \caption{surrogate model}
     \end{subfigure}  
    \caption{$u_{x}$ at $t = 4.00$ $sec$ predicted by (a) the exact model and (b) the surrogate model}
    \label{fig:ux_t4}
\end{figure}

\begin{figure}[H]
\centering
\begin{subfigure}[b]{0.48\textwidth}
         \centering
         \includegraphics[width=0.80\textwidth]{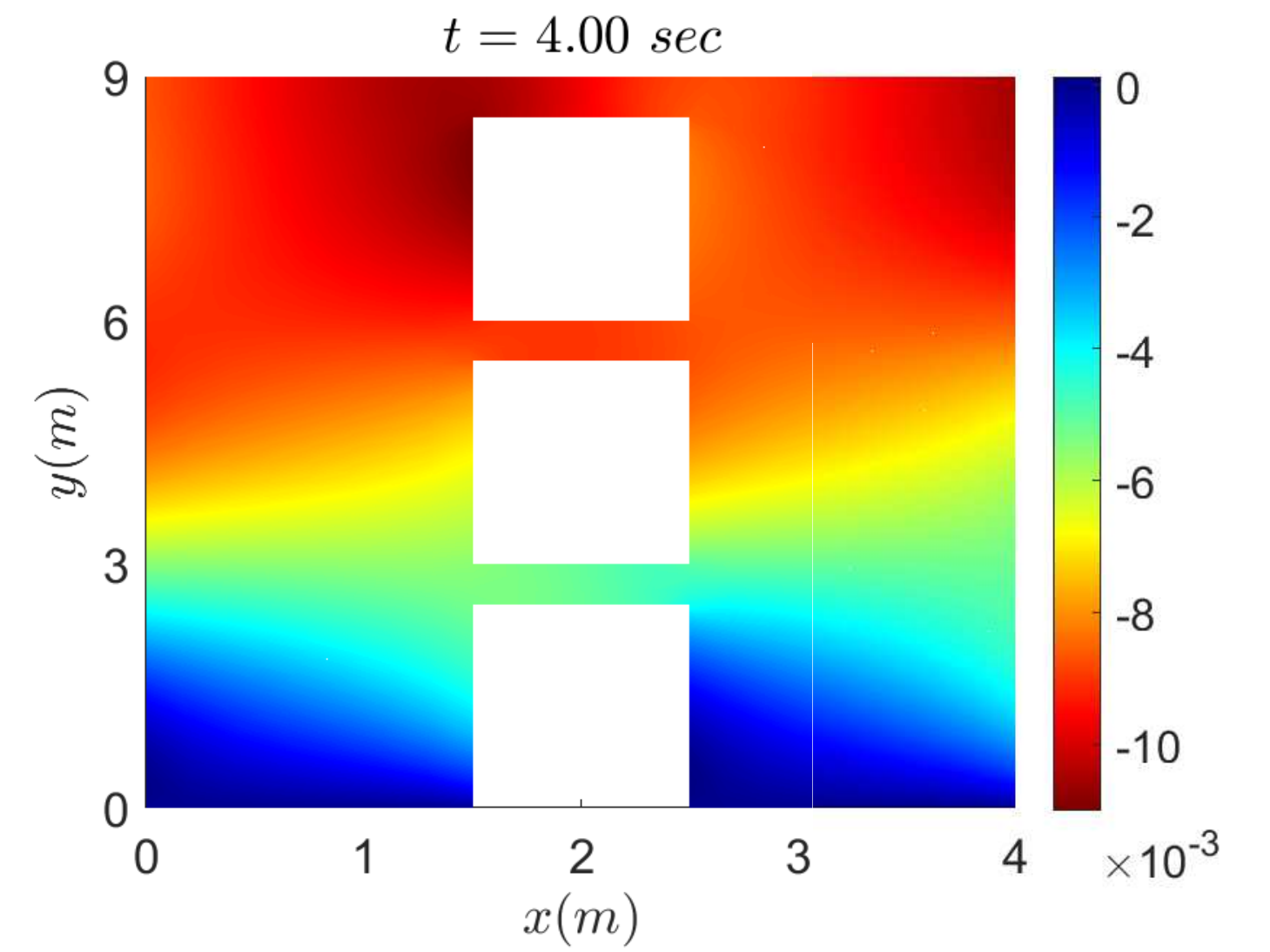}
         \caption{exact model}
     \end{subfigure}
     \hfill
     \begin{subfigure}[b]{0.48\textwidth}
         \centering
         \includegraphics[width=0.80\textwidth]{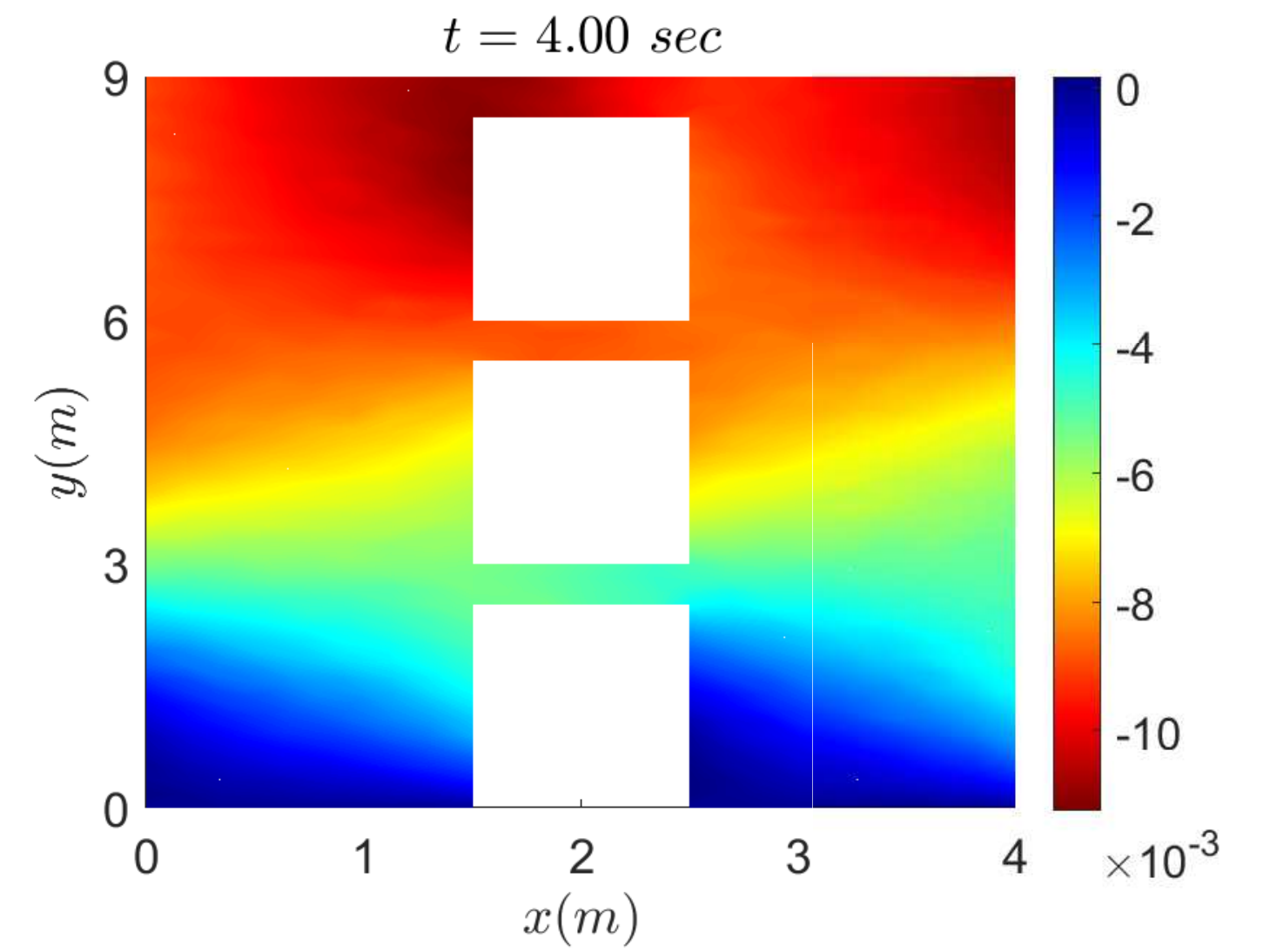}
         \caption{surrogate model}
     \end{subfigure}  
    \caption{$u_{y}$ at $t = 4.00$ $sec$ predicted by (a) the exact model and (b) the surrogate model}
    \label{fig:uy_t4}
\end{figure}

\begin{figure}[H]
     \centering
     \begin{subfigure}[b]{0.325\textwidth}
         \centering
         \includegraphics[width=\textwidth]{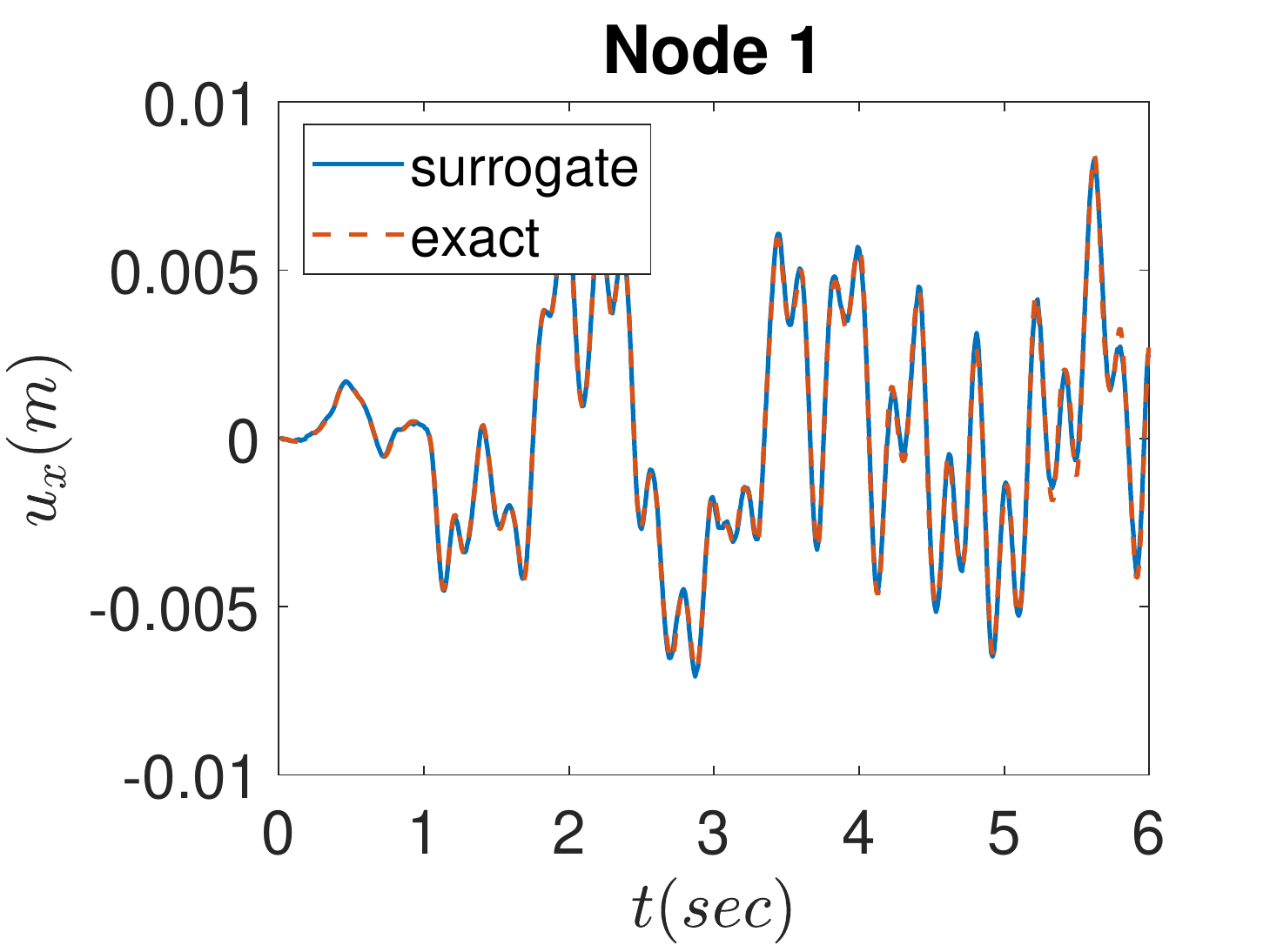}
     \end{subfigure}
     \hfill
     \begin{subfigure}[b]{0.325\textwidth}
         \centering
         \includegraphics[width=\textwidth]{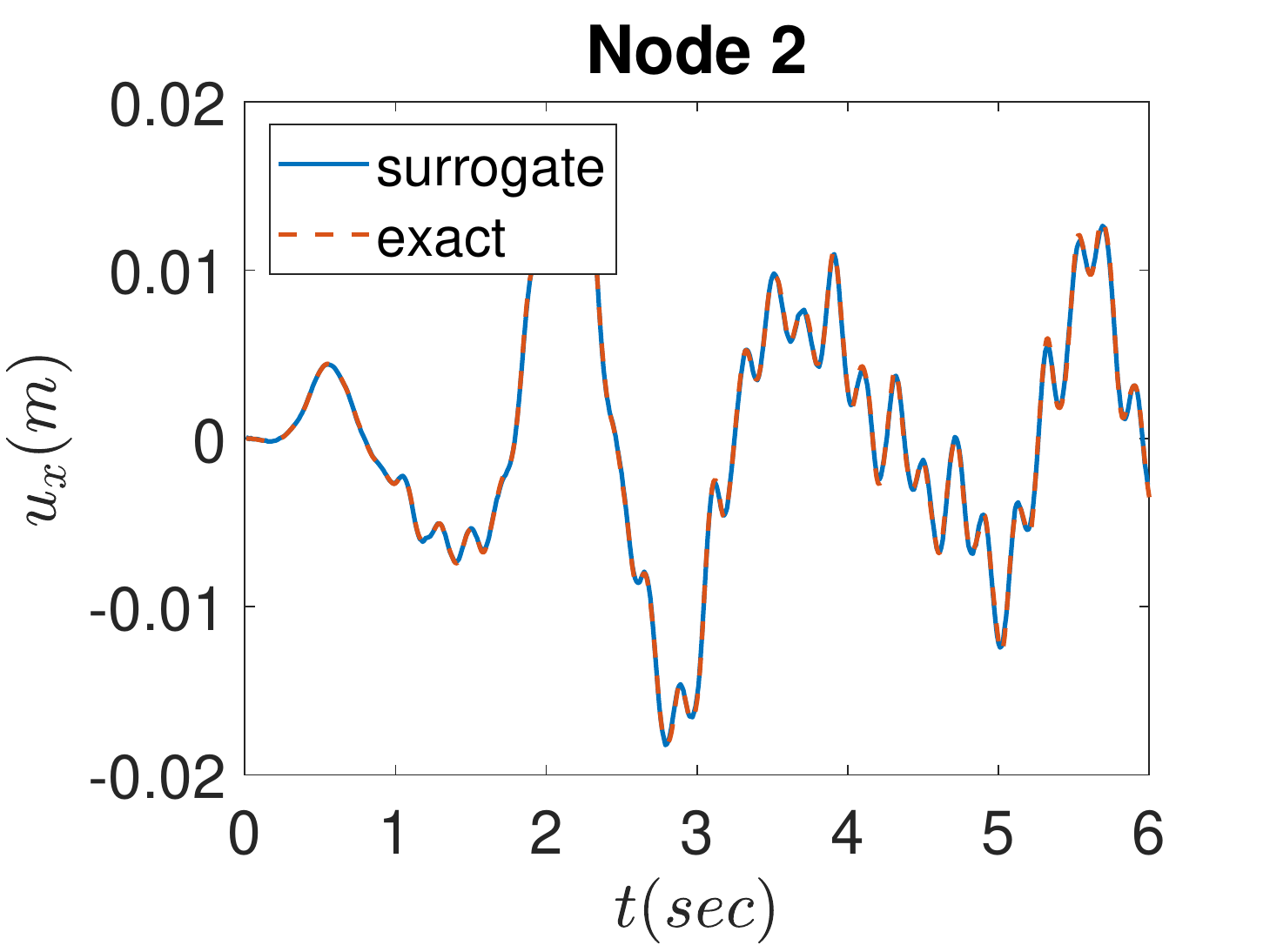}
     \end{subfigure}
     \hfill
     \begin{subfigure}[b]{0.325\textwidth}
         \centering
         \includegraphics[width=\textwidth]{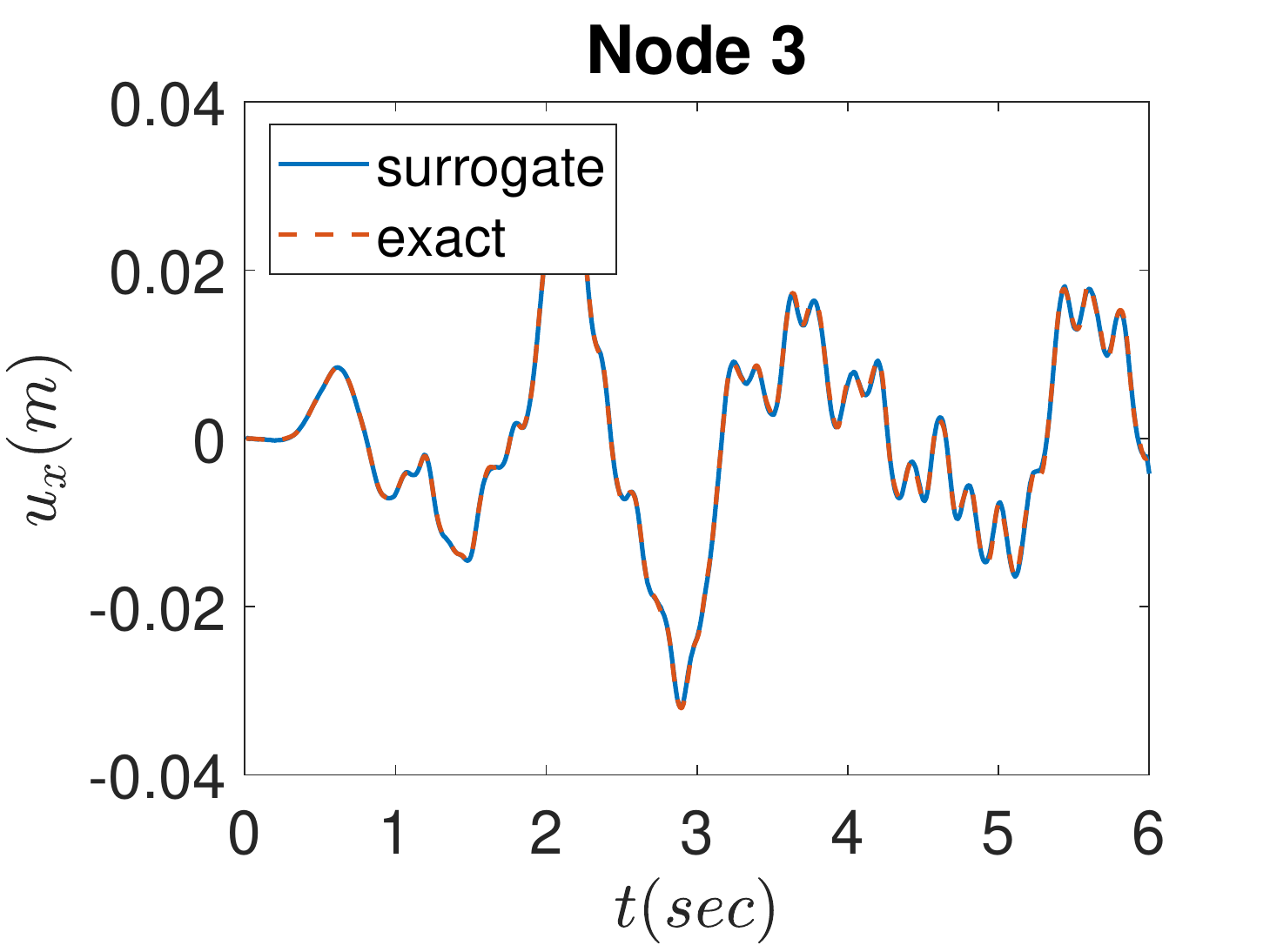}
     \end{subfigure}
     \vfill   
     \begin{subfigure}[b]{0.325\textwidth}
         \centering
         \includegraphics[width=\textwidth]{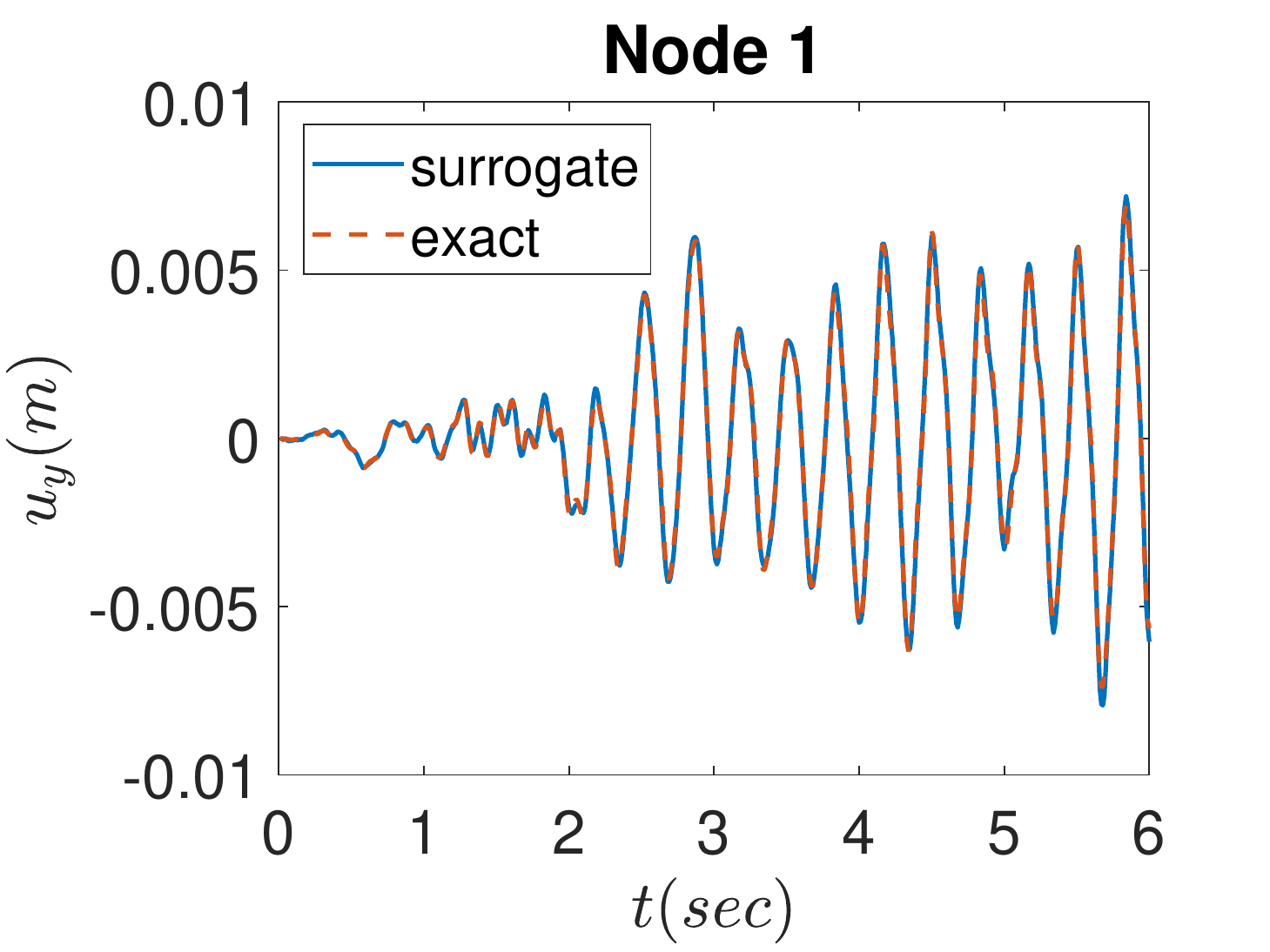}
     \end{subfigure}
      \hfill
     \begin{subfigure}[b]{0.325\textwidth}
         \centering
         \includegraphics[width=\textwidth]{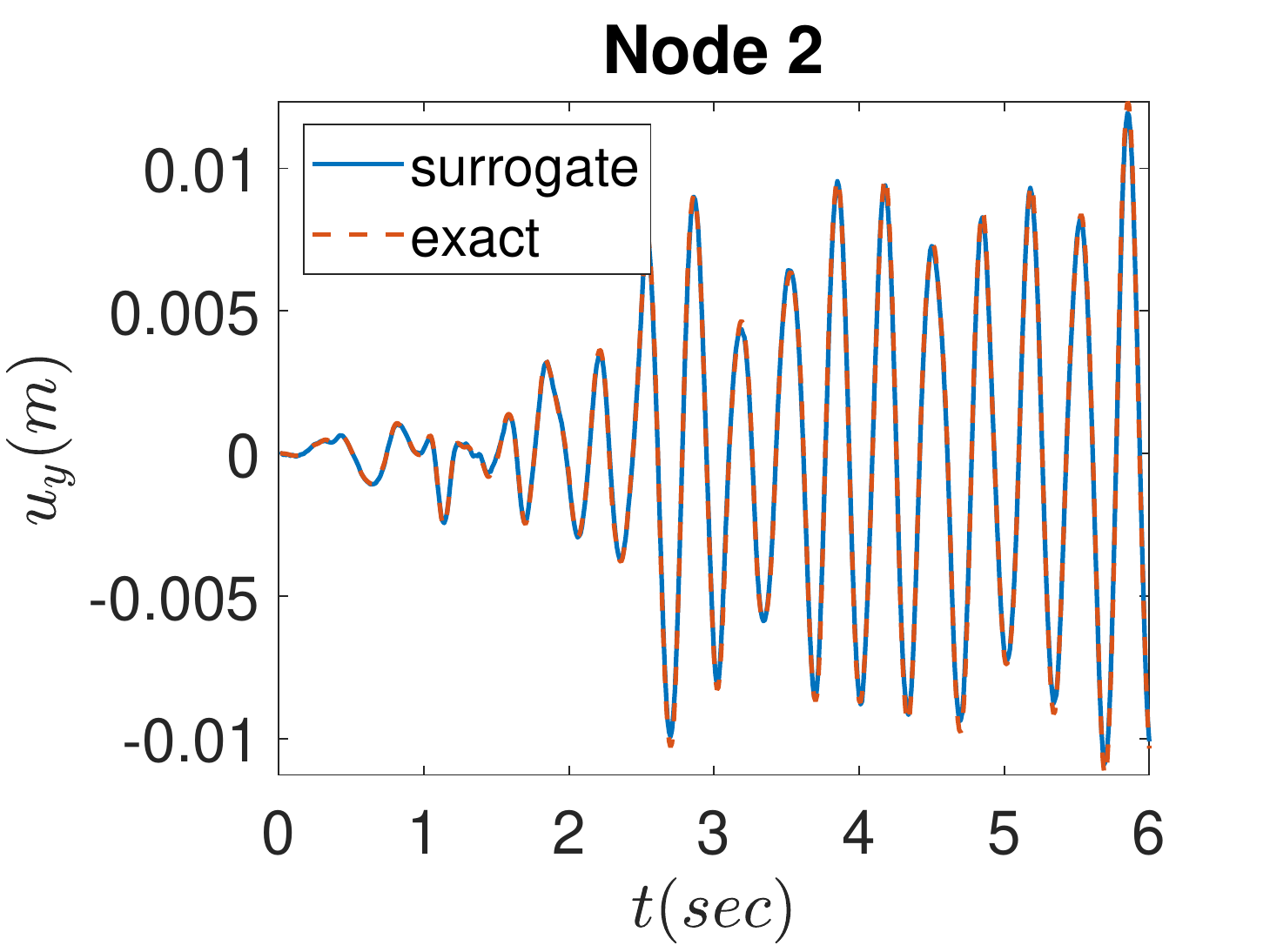}
     \end{subfigure}
     \hfill     
     \begin{subfigure}[b]{0.325\textwidth}
         \centering
         \includegraphics[width=\textwidth]{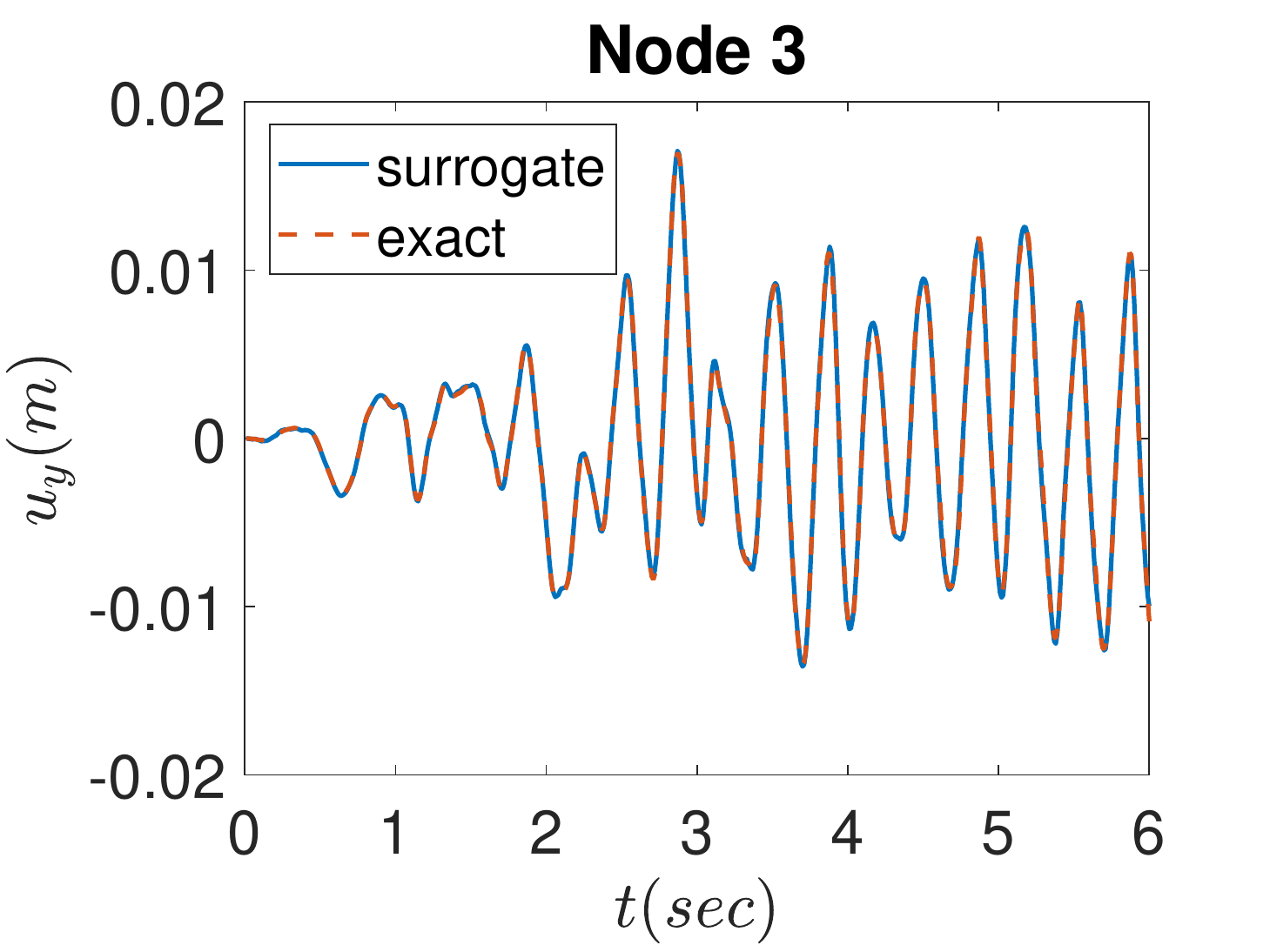}
     \end{subfigure}
        \caption{Displacements $u_{x}$ and $u_{y}$ of monitored nodes predicted by the exact and the surrogate model }
        \label{fig:monitored_nodes}
\end{figure}

In order to test the surrogate's generalization capabilities to 'unseen' parameter values, a random triplet of parameters that was not included in the training data set, $[E_{1}, E_{2}, E_{3}]$ = $[18.66, 27.02, 21.65]$ $GPa$ is selected. Figures   \ref{fig:ux_t4} and \ref{fig:uy_t4} present contour plots for the displacement fields $u_{x}$ and $u_{y}$ of the whole structure at $t = 4.00 \ sec$, predicted by the exact and the surrogate model, respectively, while figure \ref{fig:monitored_nodes} depicts a comparison between the exact and the surrogate model in the displacements $u_x$ and $u_y$ of the monitored nodes 1 through 3 (see figure \ref{fig:shear_walls}). From these results it can be observed that the predictions obtained by the surrogate model are in a near perfect match with those of the exact model. The normalized error between the solution matrices $\boldsymbol{U}_{FEM}$ and $\boldsymbol{U}_{SUR}$ of the FEM and the surrogate model, respectively, given by $\widehat{err}=\Vert\boldsymbol{U}_{FEM}-\boldsymbol{U}_{SUR}\Vert_2/ \Vert\boldsymbol{U}_{FEM}\Vert_2$, was found equal to $1.53\%$.

\begin{figure}[H]
\centering
\begin{subfigure}[b]{0.48\textwidth}
         \centering
         \includegraphics[width=0.85\textwidth]{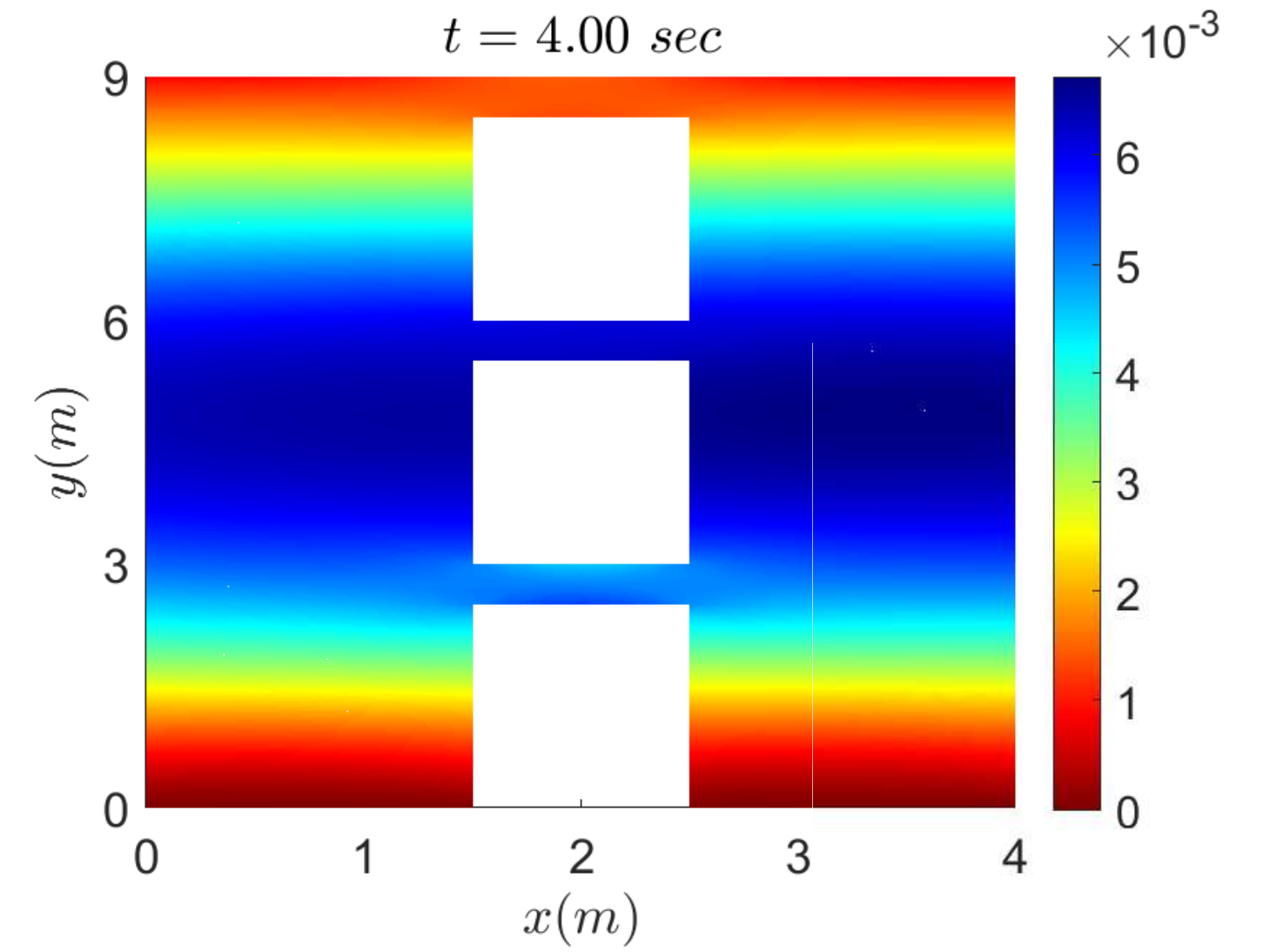}
         \caption{exact model}
     \end{subfigure}
     \hfill
     \begin{subfigure}[b]{0.48\textwidth}
         \centering
         \includegraphics[width=0.85\textwidth]{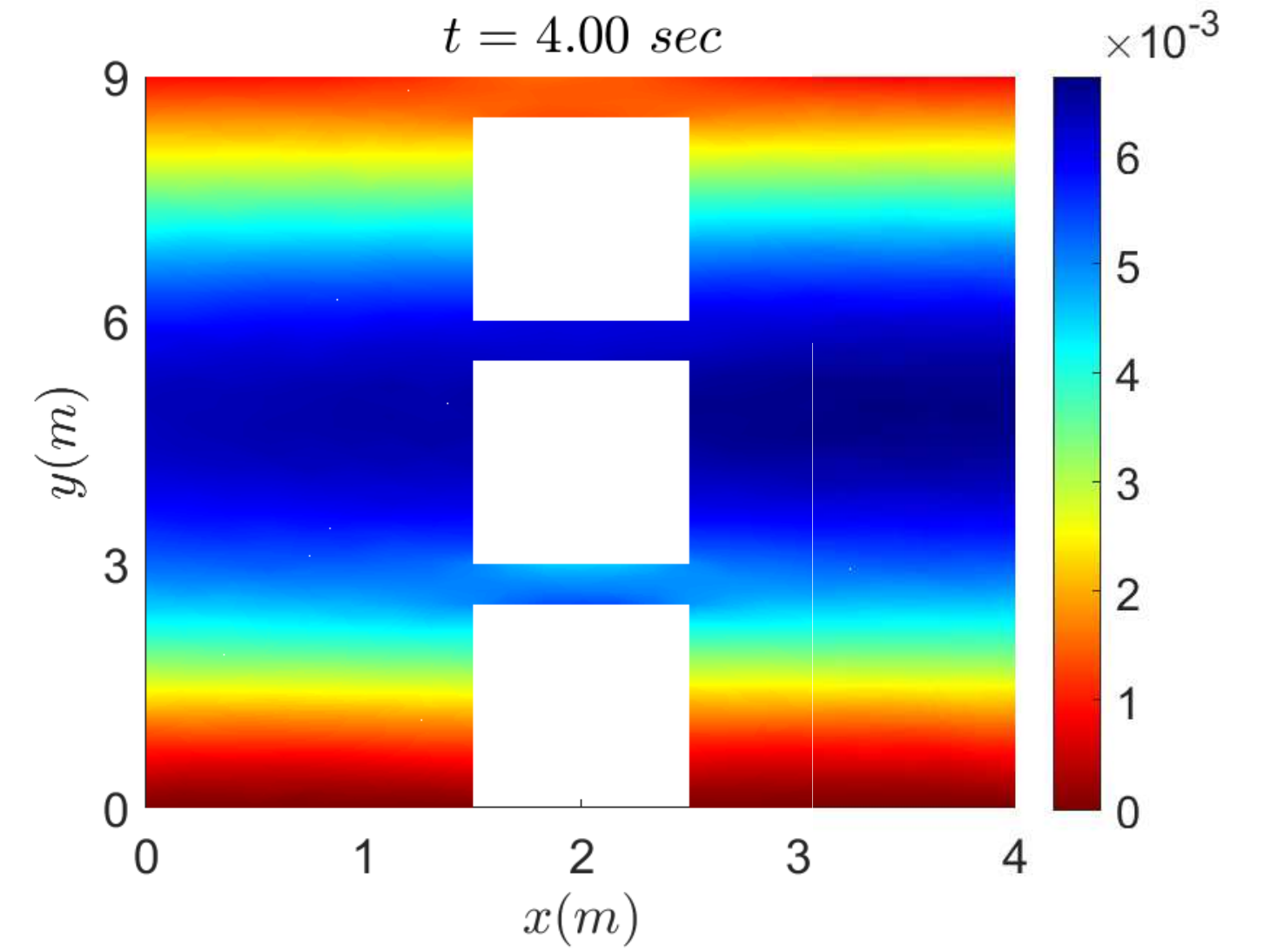}
         \caption{surrogate model}
     \end{subfigure}  
    \caption{Mean value of $u_{x}$ at $t = 4.00$ $sec$ predicted by (a) the exact model and (b) the surrogate model}
    \label{fig:mean_ux_t4}
\end{figure}

\begin{figure}[H]
\centering
\begin{subfigure}[b]{0.48\textwidth}
         \centering
         \includegraphics[width=0.85\textwidth]{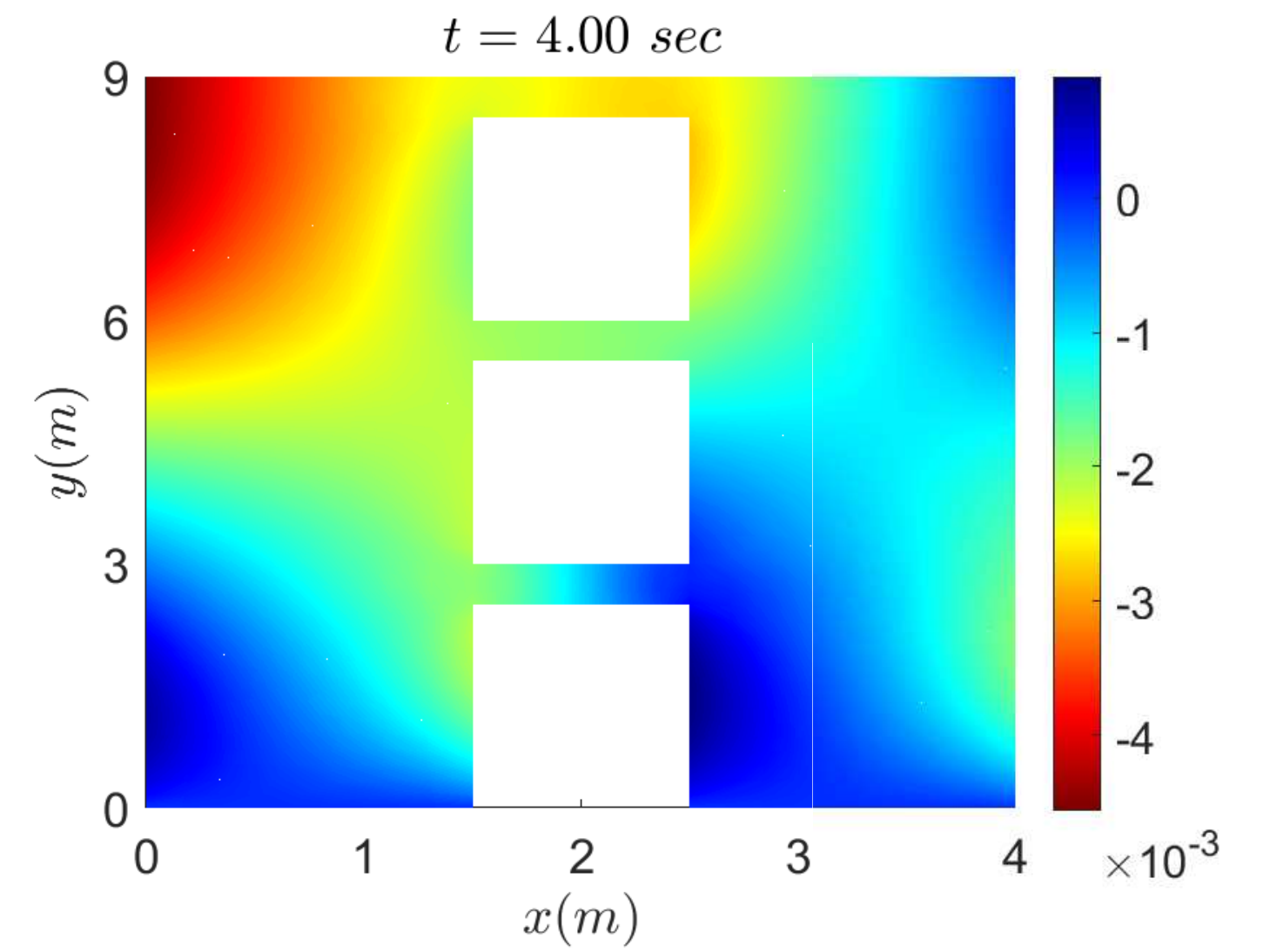}
         \caption{exact model}
     \end{subfigure}
     \hfill
     \begin{subfigure}[b]{0.48\textwidth}
         \centering
         \includegraphics[width=0.85\textwidth]{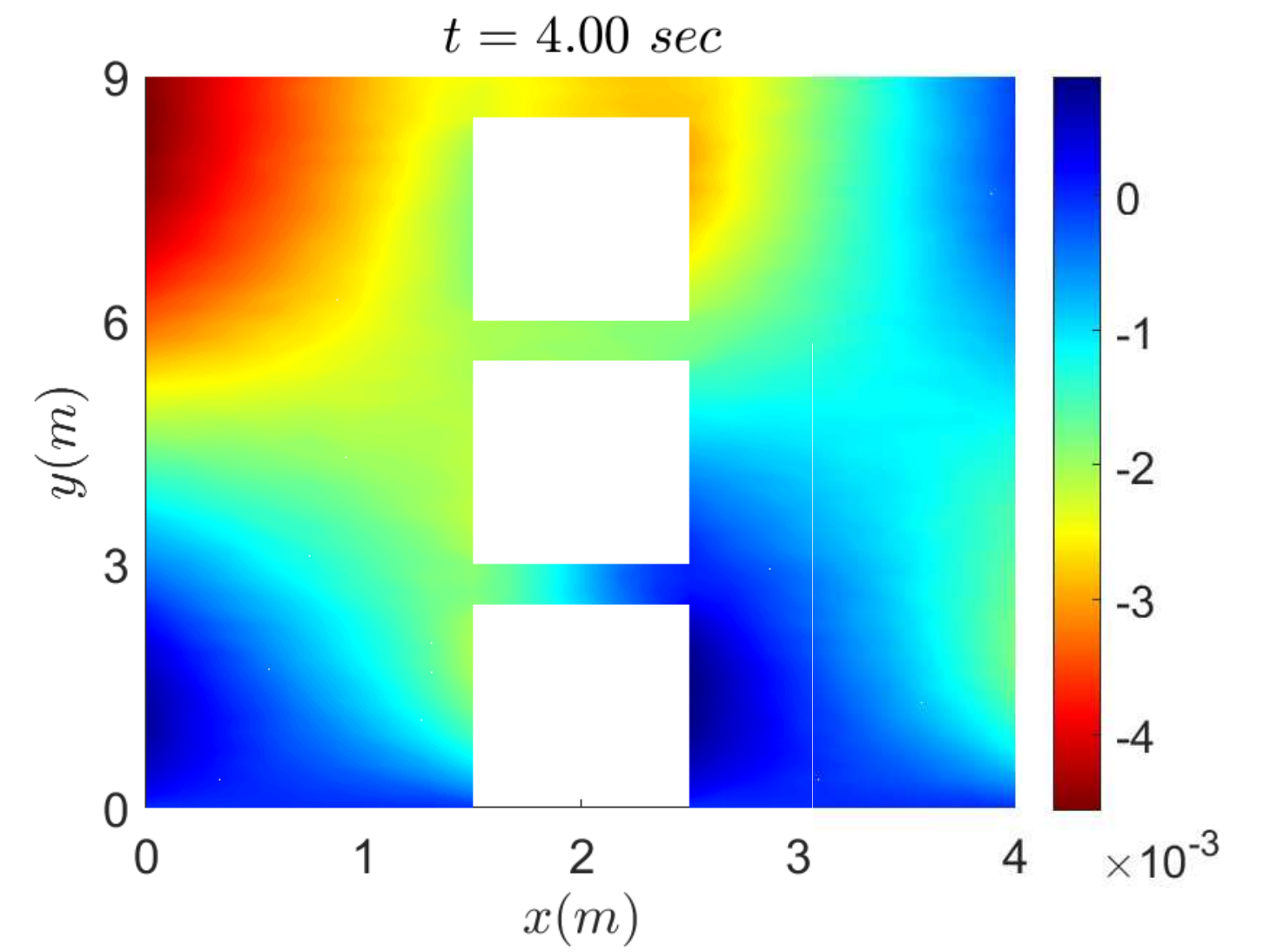}
         \caption{surrogate model}
     \end{subfigure}  
    \caption{Mean value of $u_{y}$ at $t = 4.00$ $sec$ predicted by (a) the exact model and (b) the surrogate model}
    \label{fig:mean_uy_t4}
\end{figure}

\begin{figure}[H]
\centering
\begin{subfigure}[b]{0.48\textwidth}
         \centering
         \includegraphics[width=0.85\textwidth]{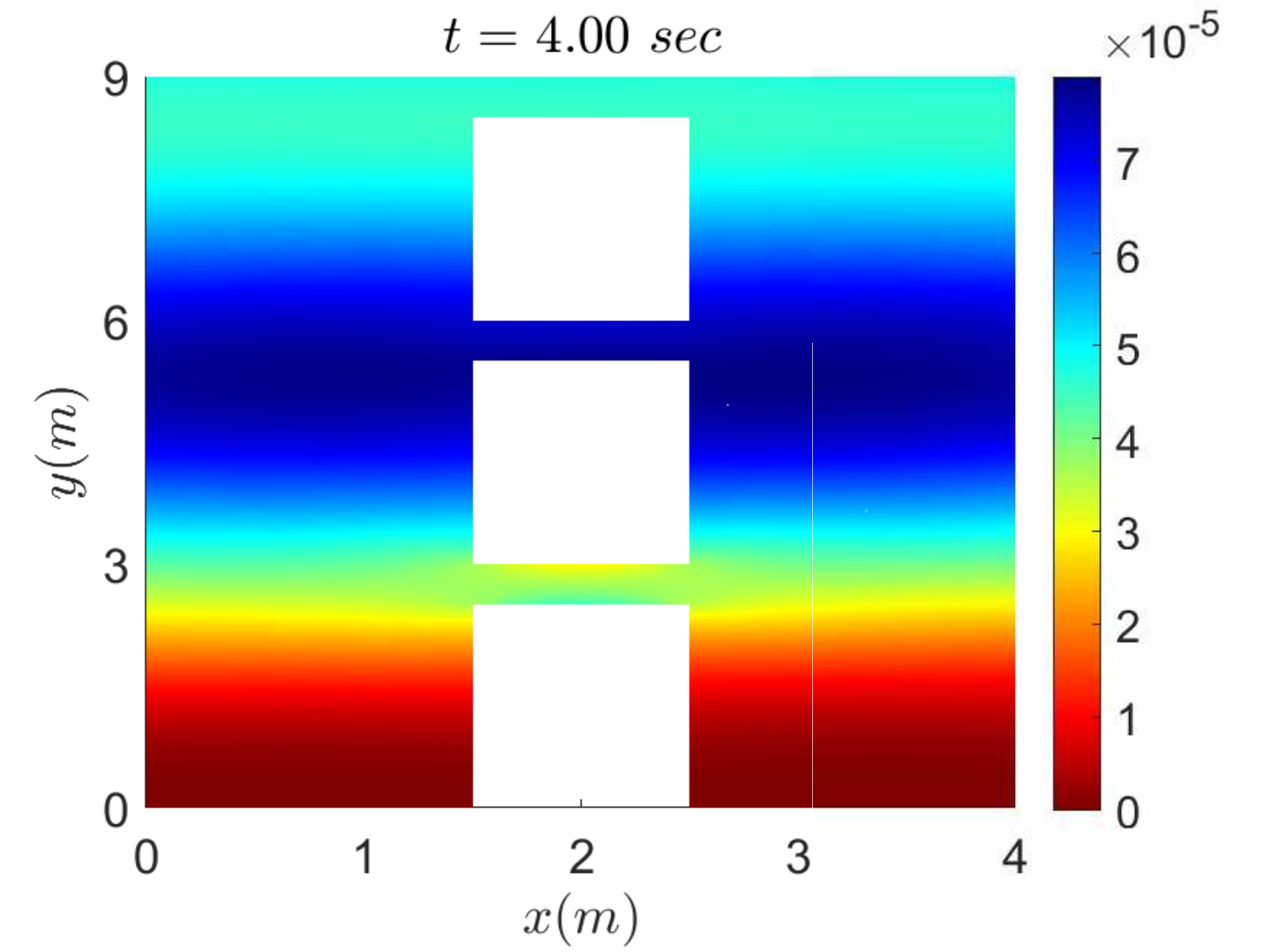}
         \caption{exact model}
     \end{subfigure}
     \hfill
     \begin{subfigure}[b]{0.48\textwidth}
         \centering
         \includegraphics[width=0.85\textwidth]{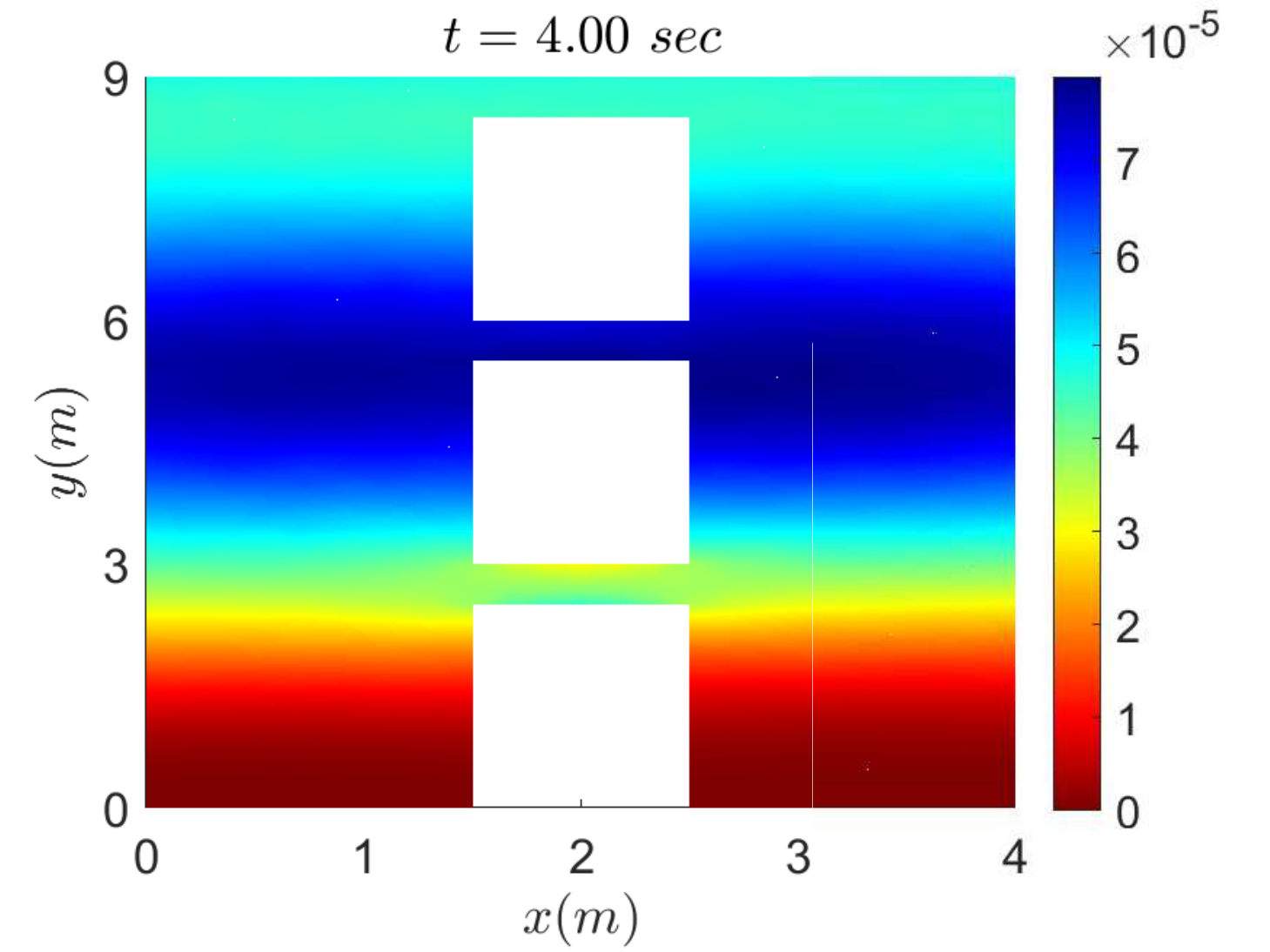}
         \caption{surrogate model}
     \end{subfigure}  
    \caption{Variance of $u_{x}$ at $t = 4.00$ $sec$ predicted by (a) the exact model and (b) the surrogate model}
    \label{fig:var_ux_t4}
\end{figure}

\begin{figure}[H]
\centering
\begin{subfigure}[b]{0.48\textwidth}
         \centering
         \includegraphics[width=0.85\textwidth]{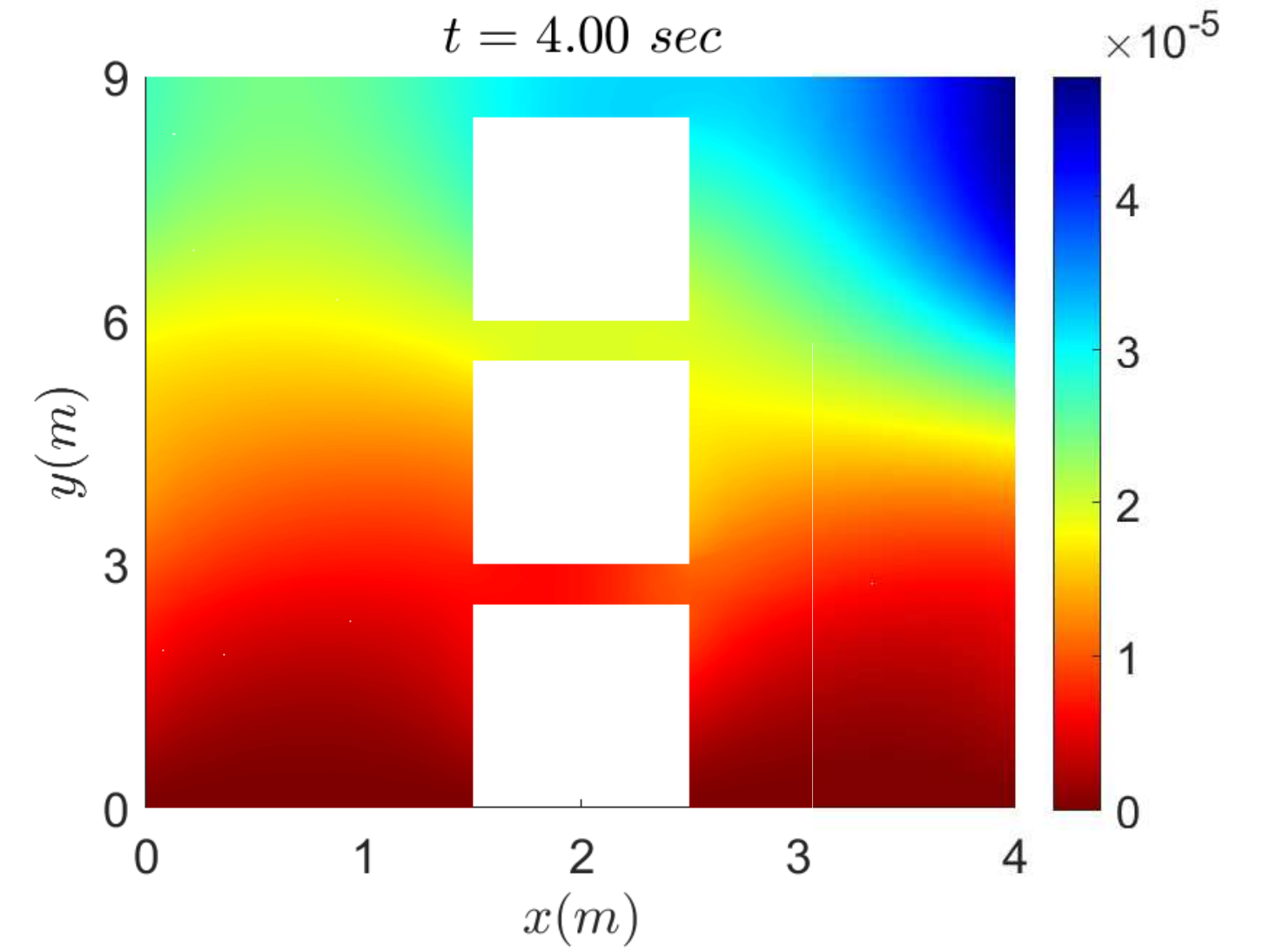}
         \caption{exact model}
     \end{subfigure}
     \hfill
     \begin{subfigure}[b]{0.48\textwidth}
         \centering
         \includegraphics[width=0.85\textwidth]{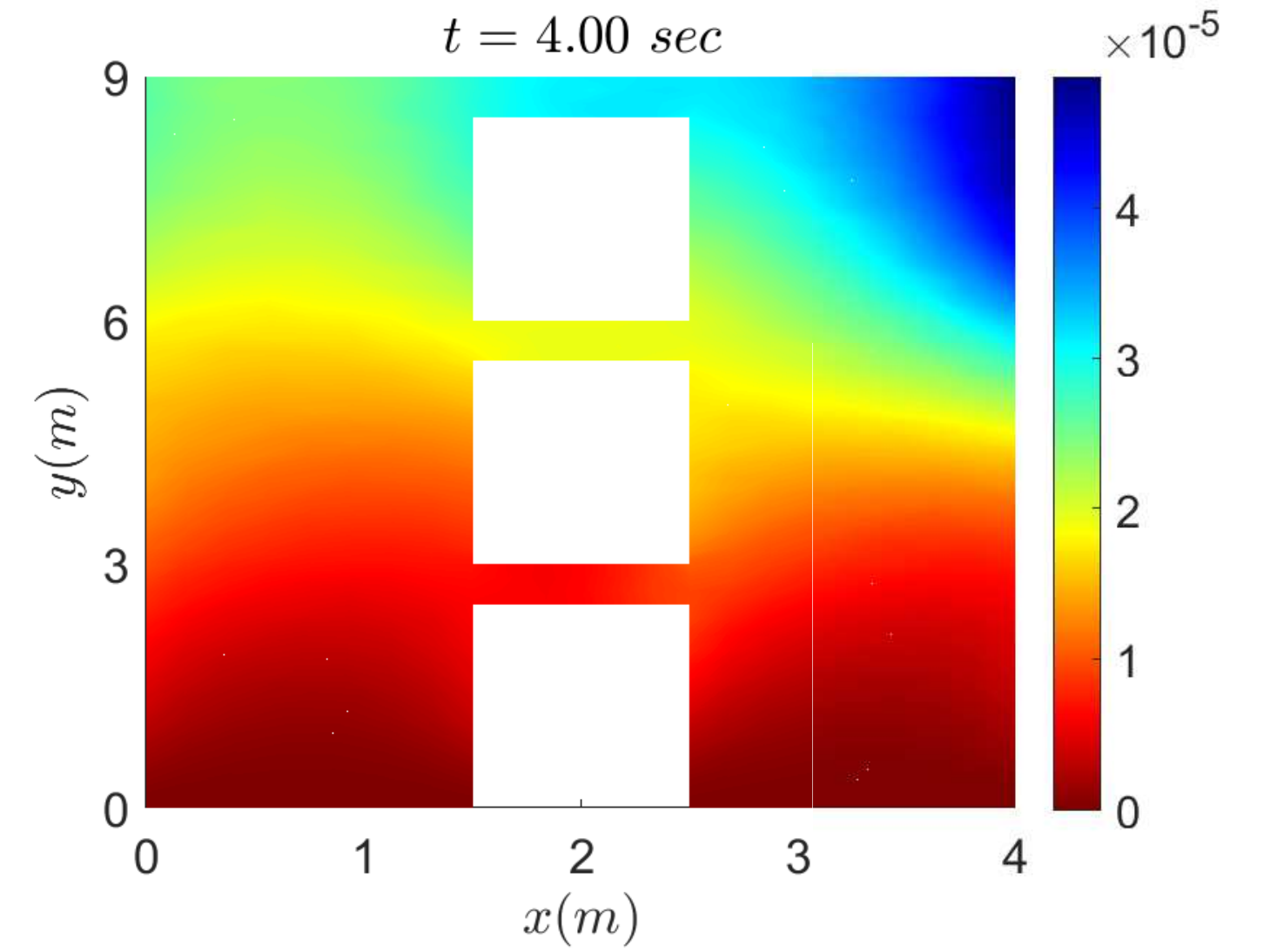}
         \caption{surrogate model}
     \end{subfigure}  
    \caption{Variance of $u_{y}$ at $t = 4.00$ $sec$ predicted by (a) the exact model and (b) the surrogate model}
    \label{fig:var_uy_t4}
\end{figure}

\begin{figure}[H]
     \centering
     \begin{subfigure}[b]{0.325\textwidth}
         \centering
         \includegraphics[width=\textwidth]{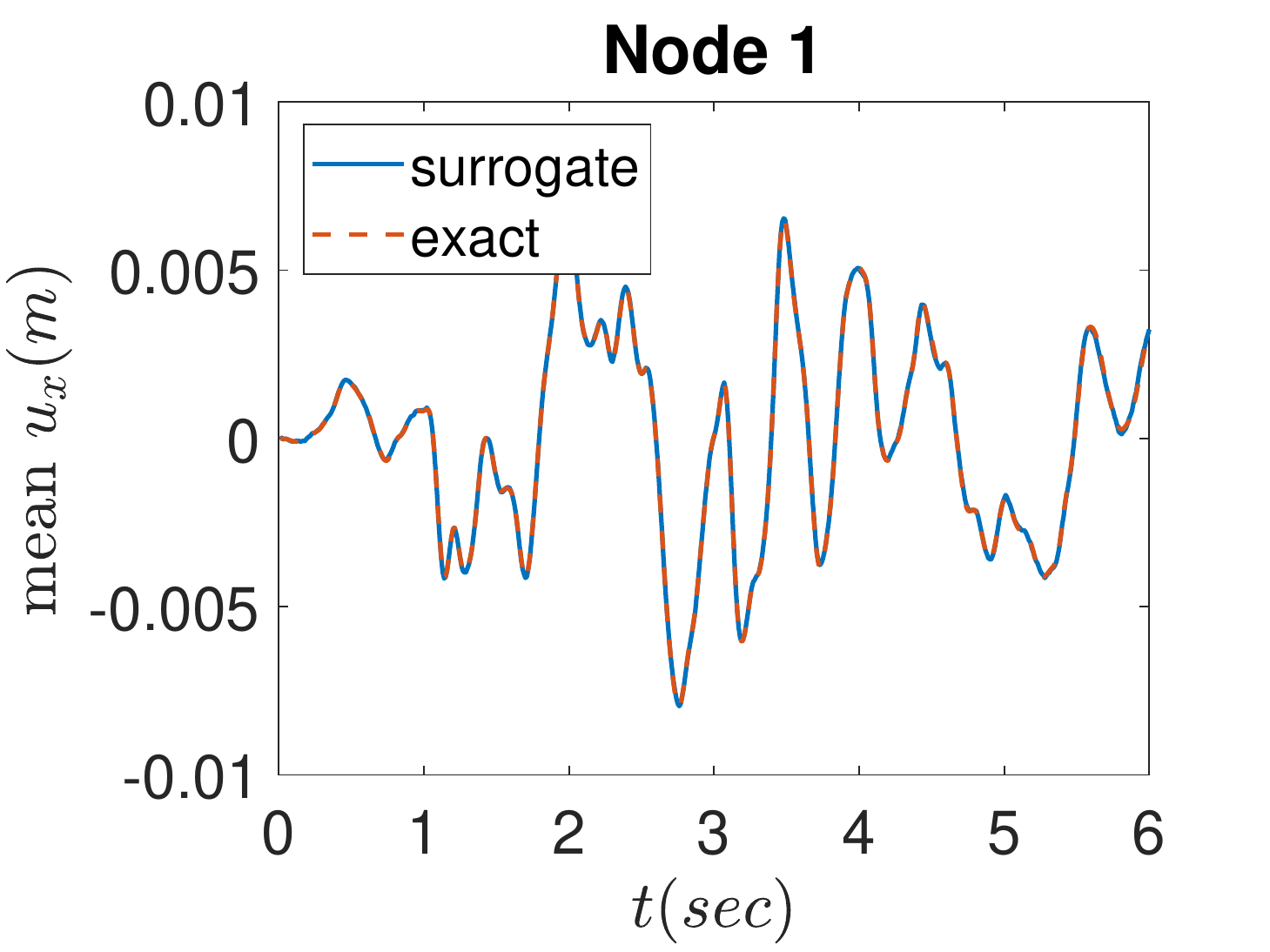}
     \end{subfigure}
     \hfill
     \begin{subfigure}[b]{0.325\textwidth}
         \centering
         \includegraphics[width=\textwidth]{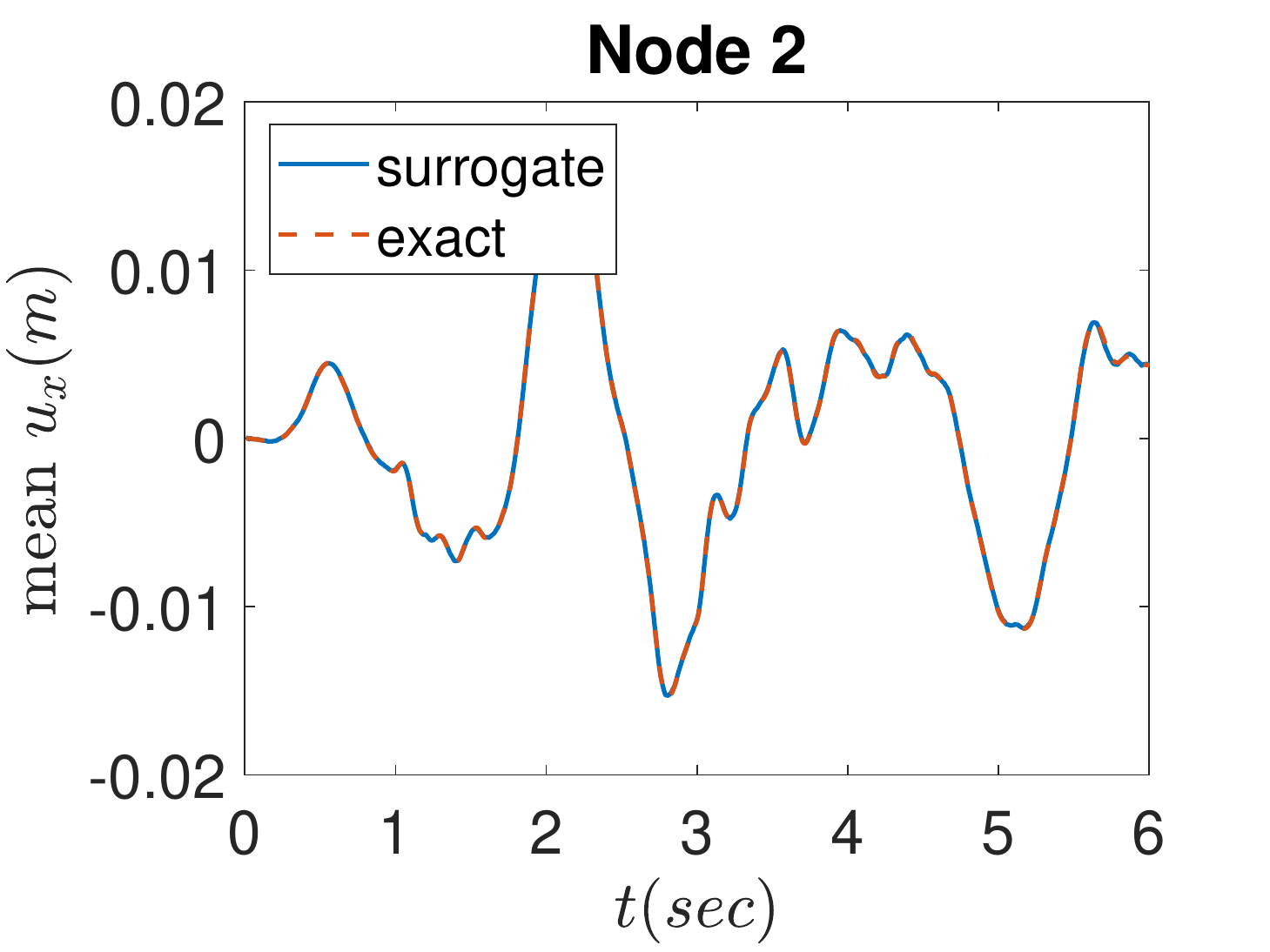}
     \end{subfigure}
     \hfill
     \begin{subfigure}[b]{0.325\textwidth}
         \centering
         \includegraphics[width=\textwidth]{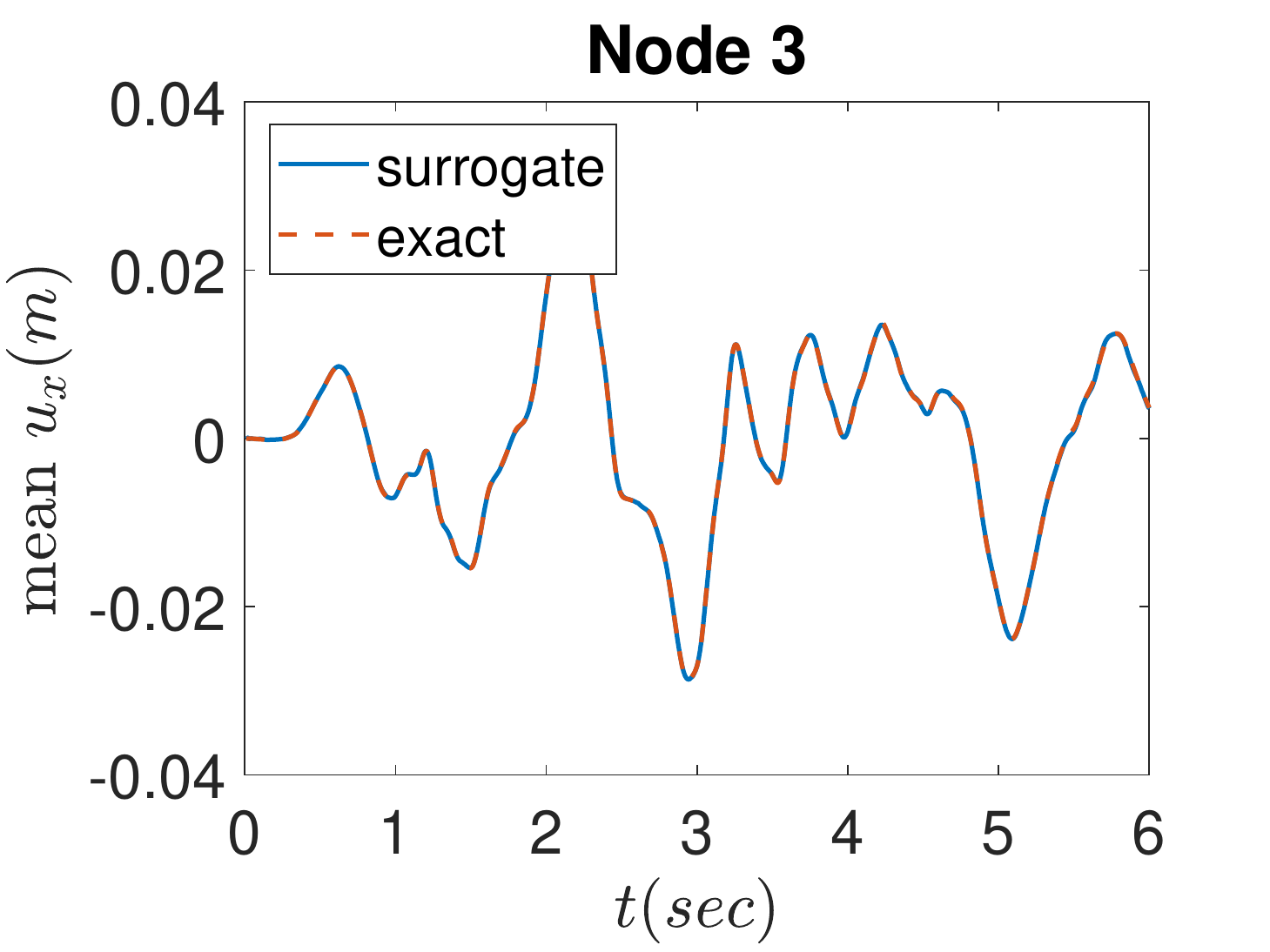}
     \end{subfigure}
     \vfill    
     \begin{subfigure}[b]{0.325\textwidth}
         \centering
         \includegraphics[width=\textwidth]{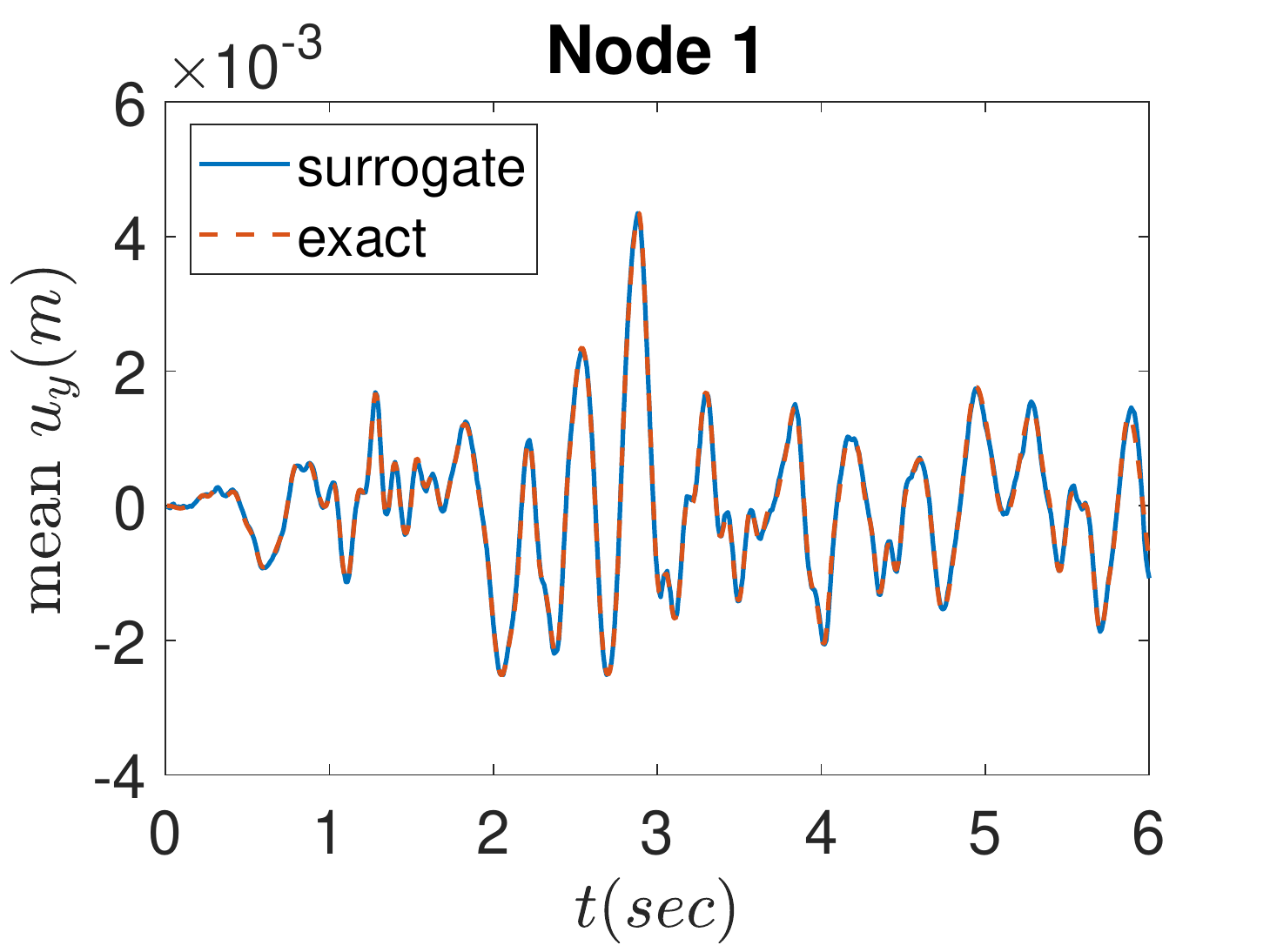}
     \end{subfigure}
      \hfill
     \begin{subfigure}[b]{0.325\textwidth}
         \centering
         \includegraphics[width=\textwidth]{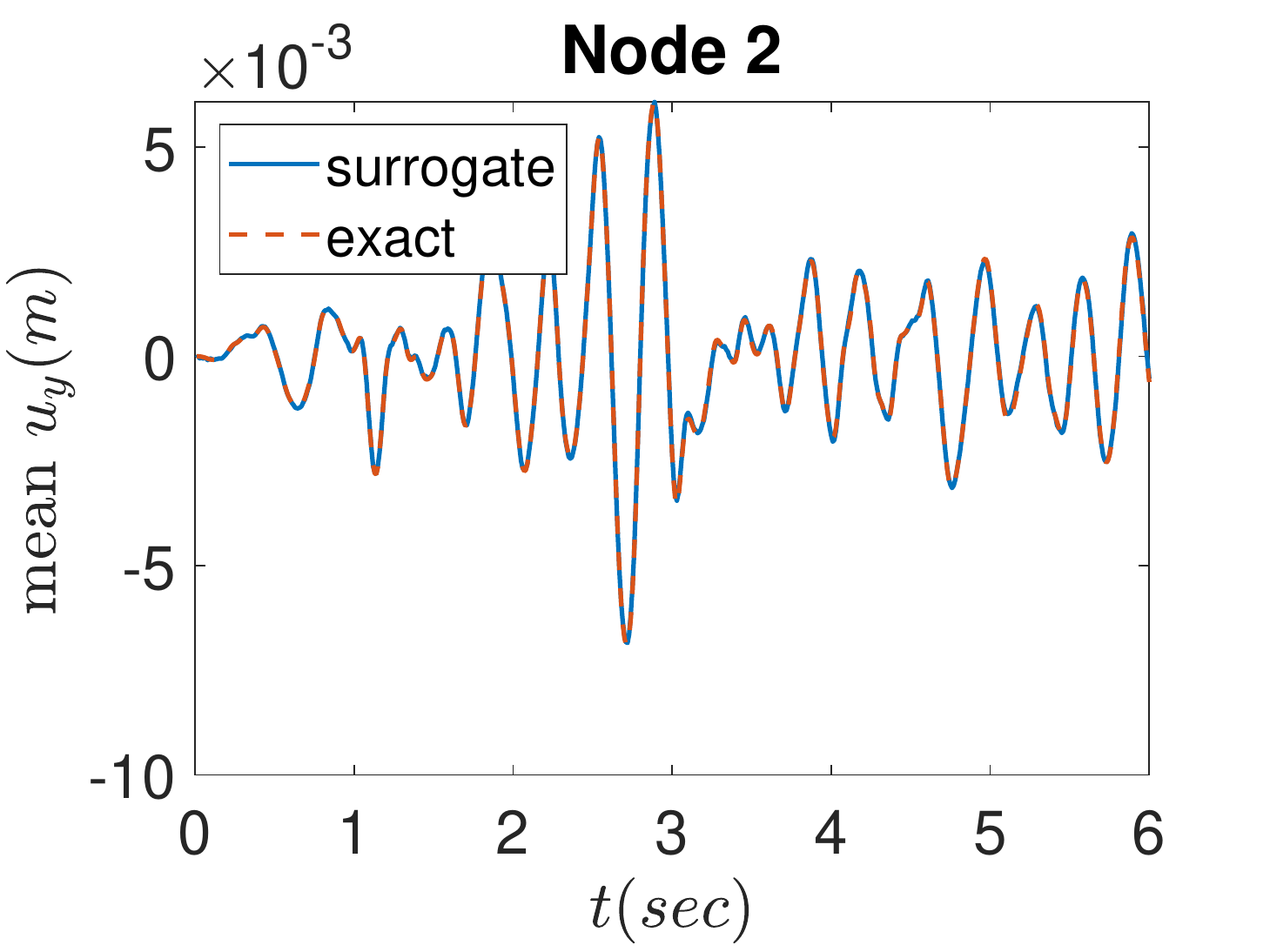}
     \end{subfigure}
     \hfill     
     \begin{subfigure}[b]{0.325\textwidth}
         \centering
         \includegraphics[width=\textwidth]{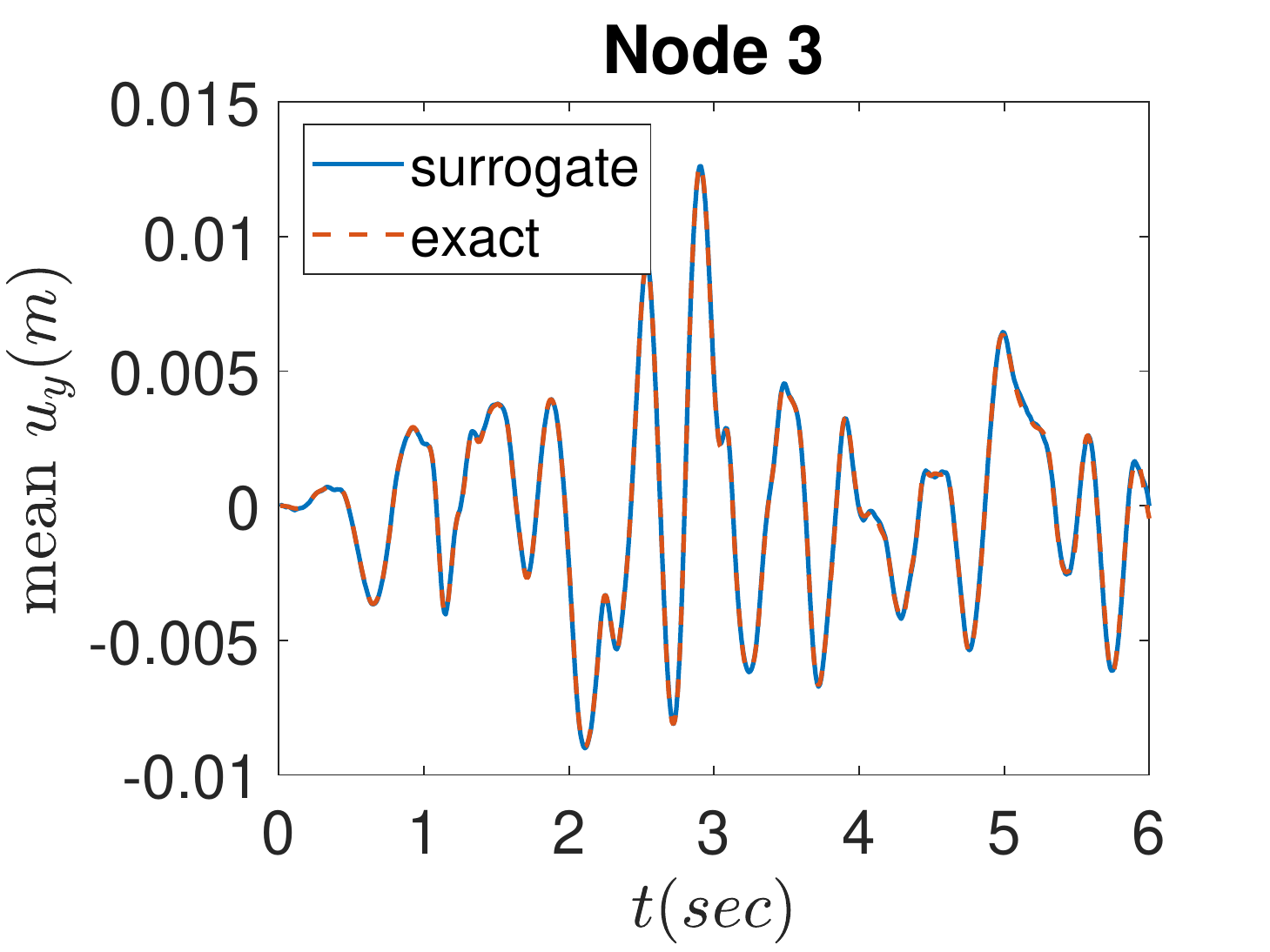}
     \end{subfigure}
        \caption{Mean $u_{x}$ and $u_{y}$ of monitored nodes predicted by the exact and the surrogate model }
        \label{fig:monitored_nodes_mean}
\end{figure}

\begin{figure}[H]
     \centering
     \begin{subfigure}[b]{0.325\textwidth}
         \centering
         \includegraphics[width=\textwidth]{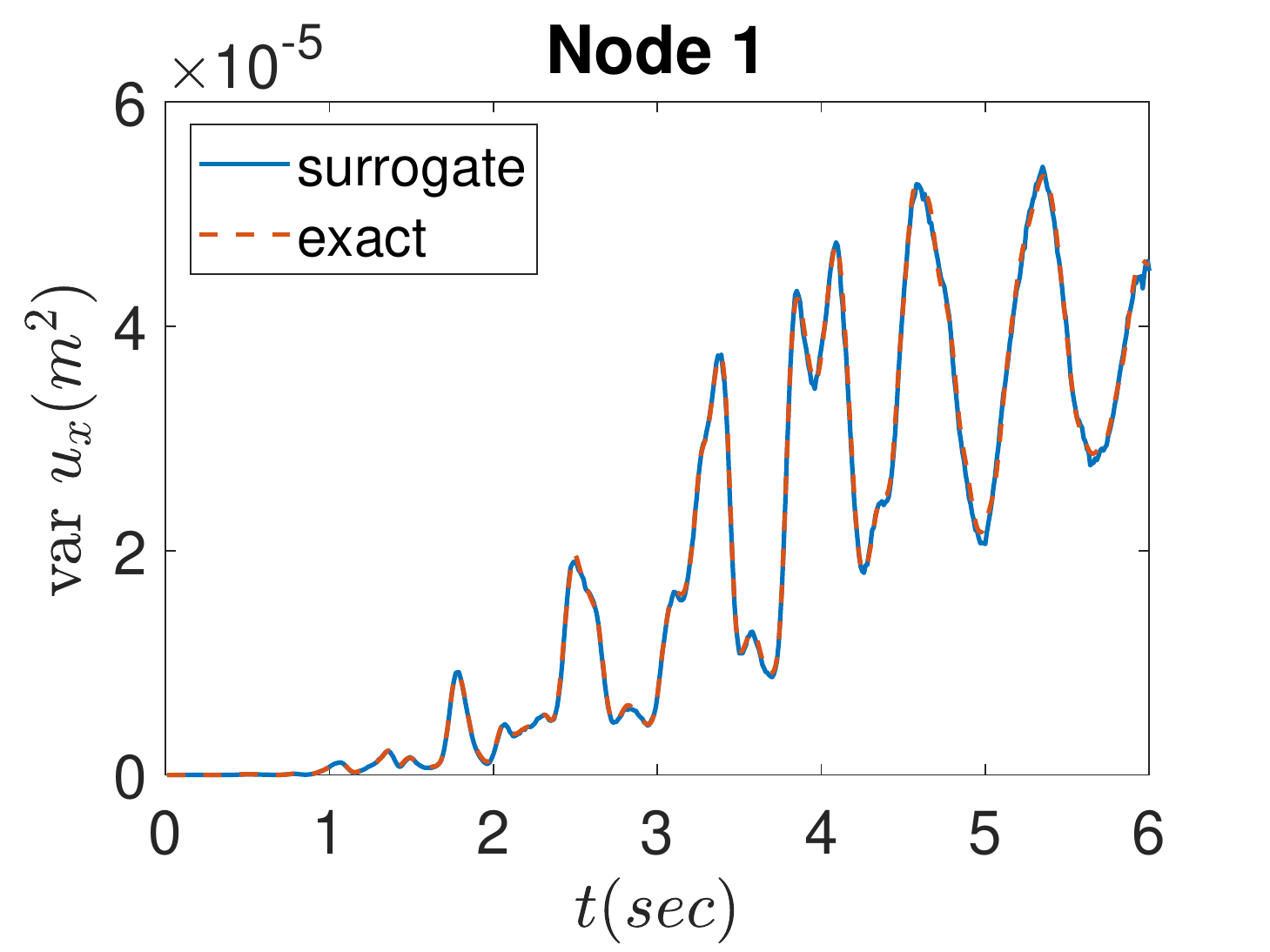}
     \end{subfigure}
     \hfill
     \begin{subfigure}[b]{0.325\textwidth}
         \centering
         \includegraphics[width=\textwidth]{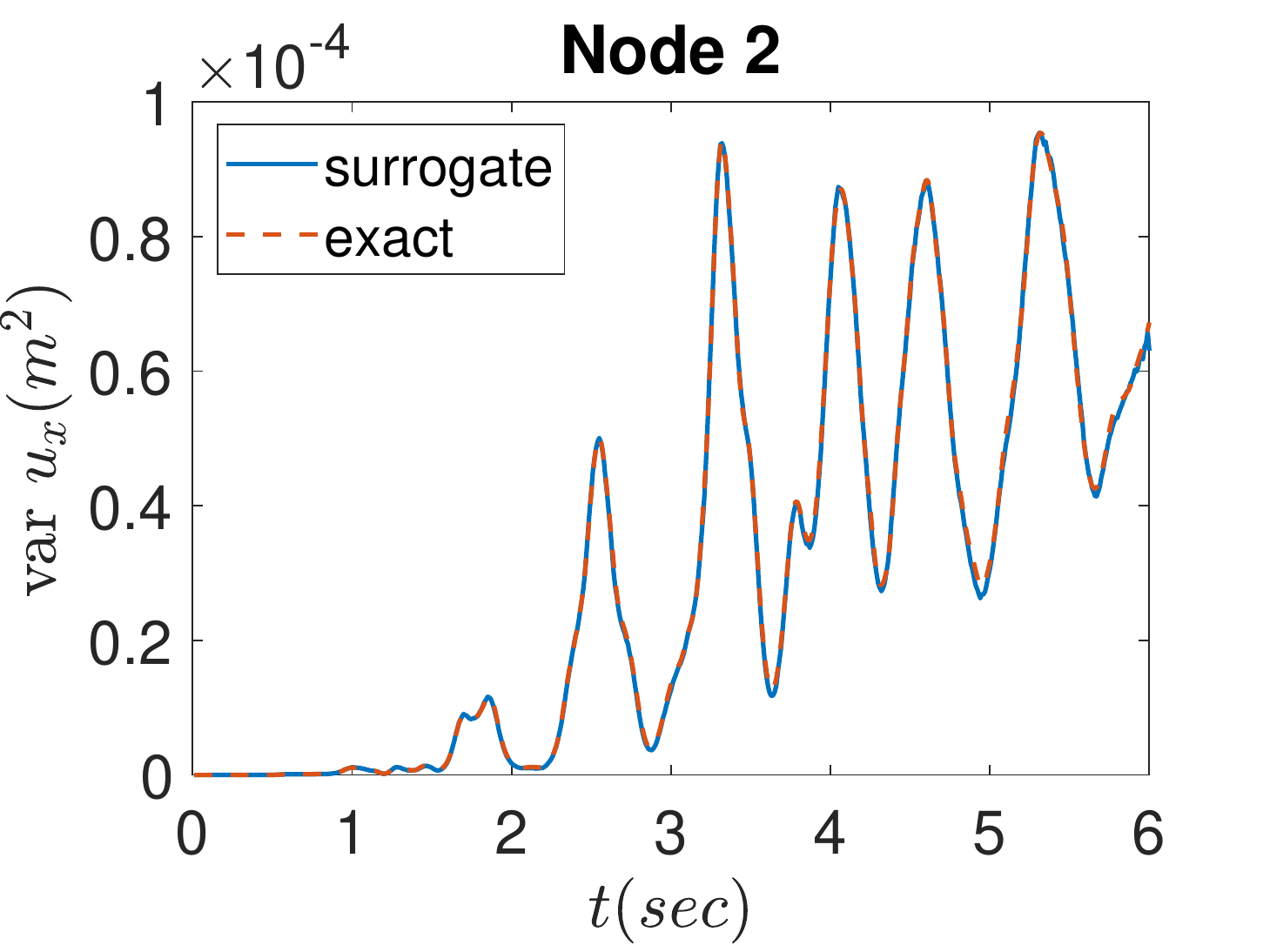}
     \end{subfigure}
     \hfill
     \begin{subfigure}[b]{0.325\textwidth}
         \centering
         \includegraphics[width=\textwidth]{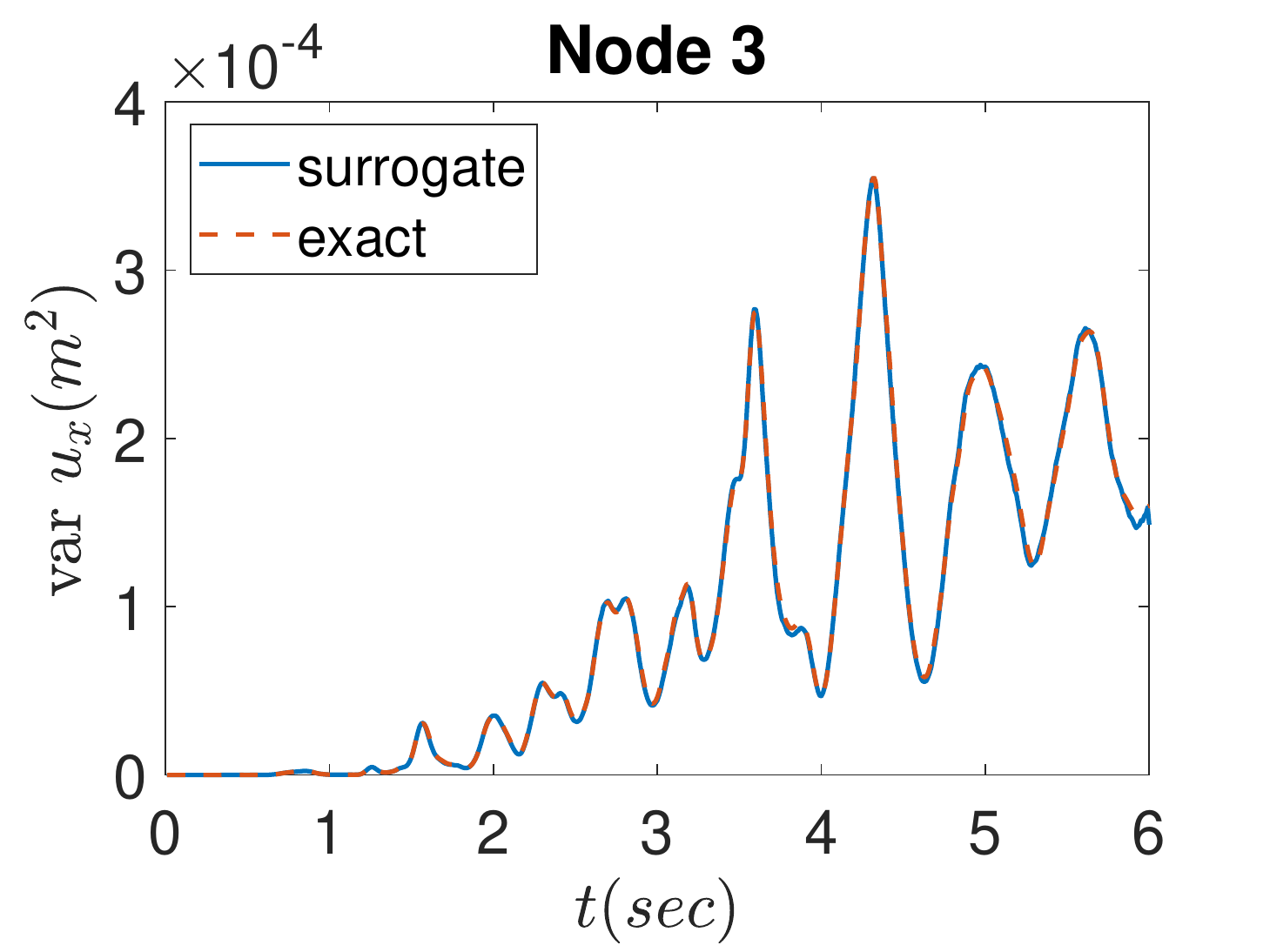}
     \end{subfigure}
     \vfill     
     \begin{subfigure}[b]{0.325\textwidth}
         \centering
         \includegraphics[width=\textwidth]{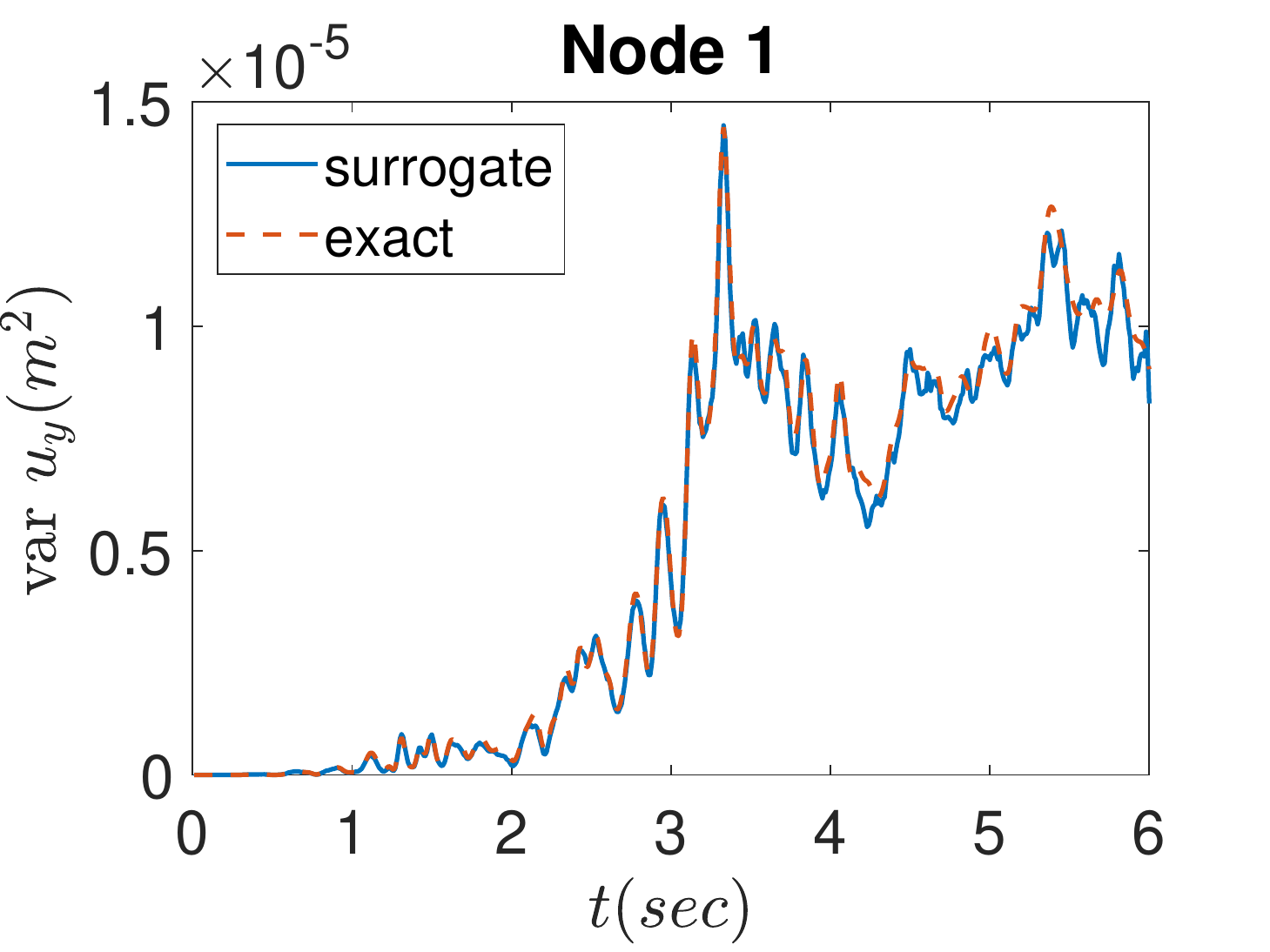}
     \end{subfigure}
      \hfill
     \begin{subfigure}[b]{0.325\textwidth}
         \centering
         \includegraphics[width=\textwidth]{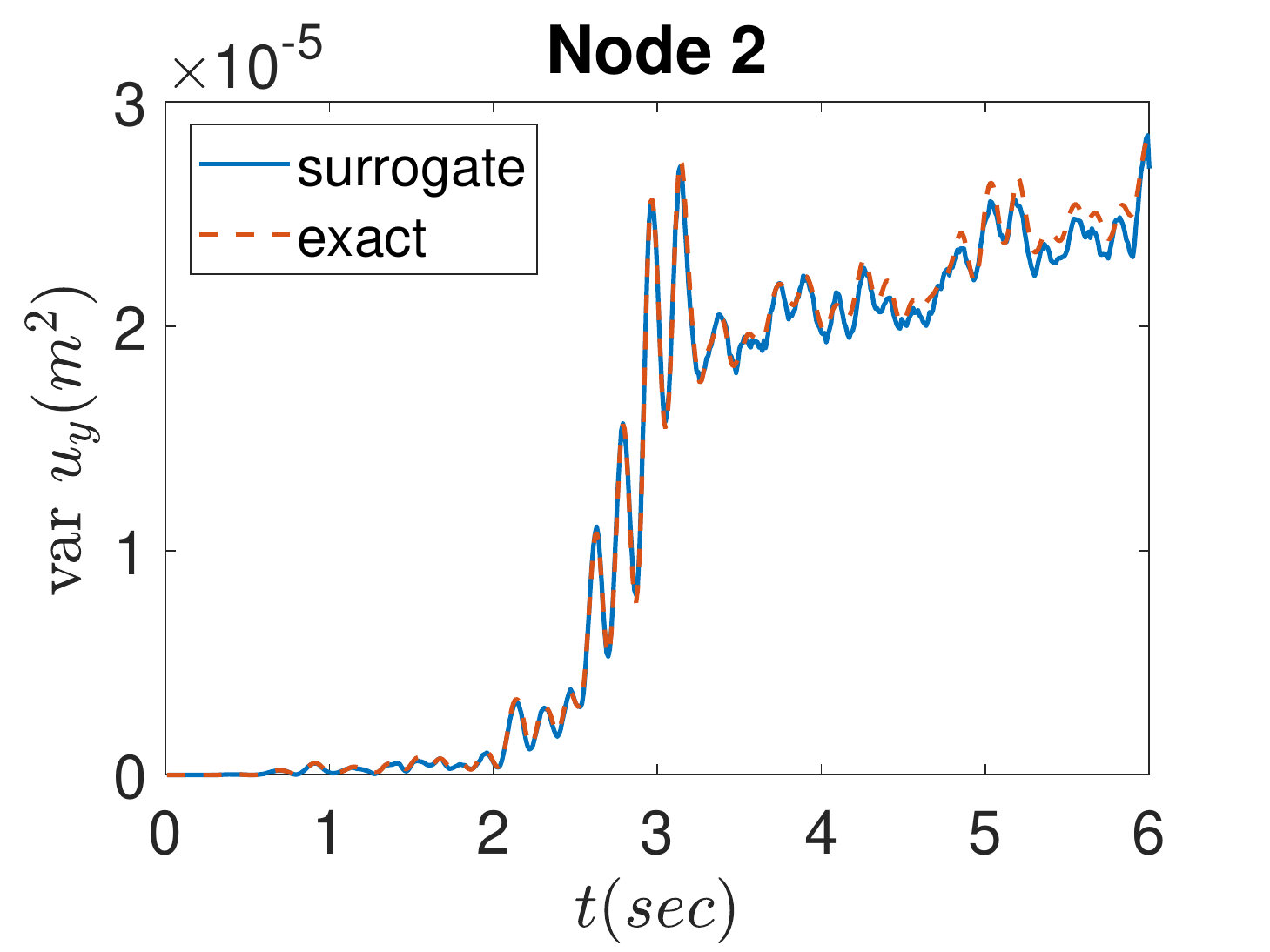}
     \end{subfigure}
     \hfill     
     \begin{subfigure}[b]{0.325\textwidth}
         \centering
         \includegraphics[width=\textwidth]{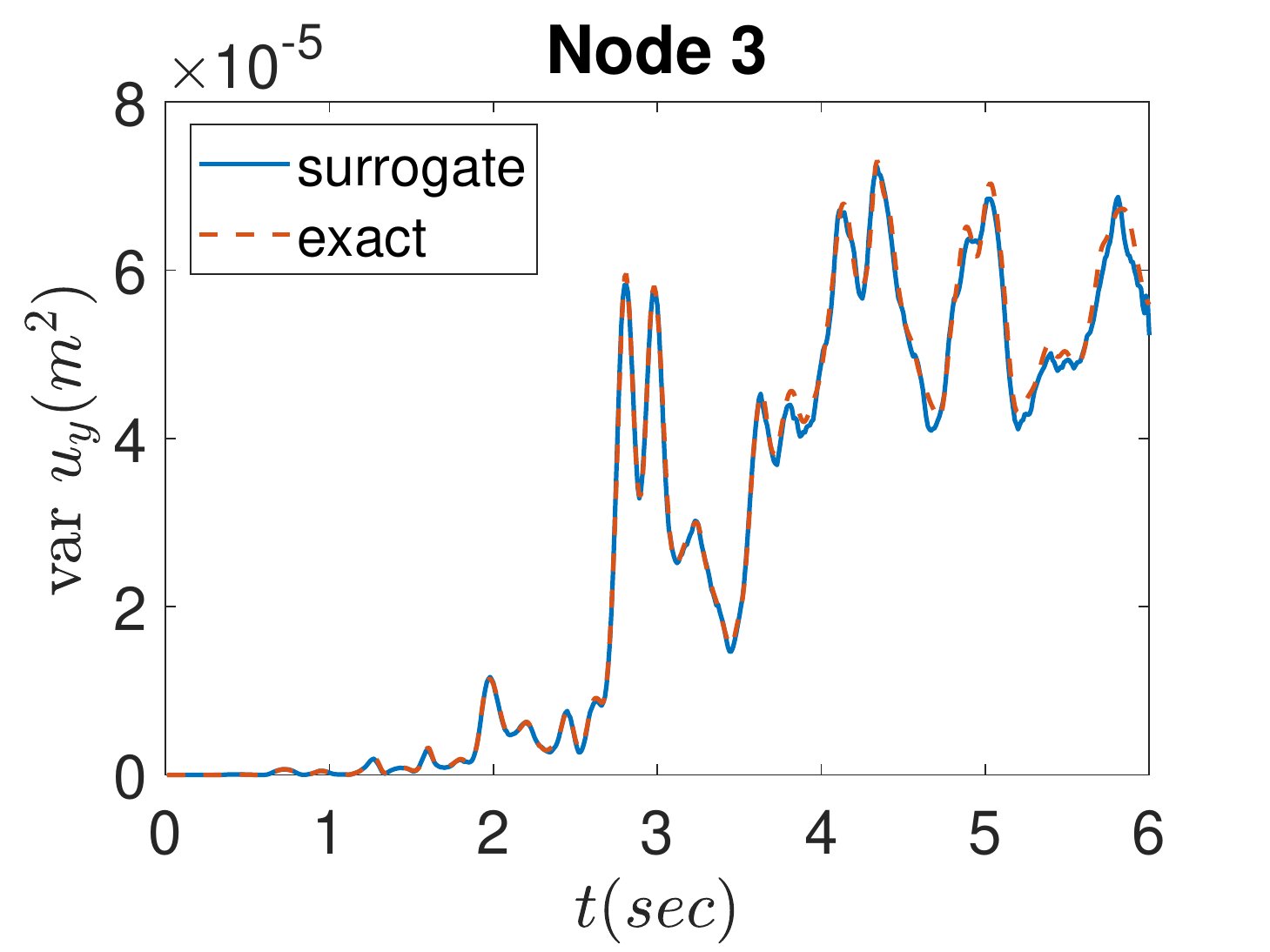}
     \end{subfigure}
        \caption{Variance of $u_{x}$ and $u_{y}$ of monitored nodes predicted by the exact and the surrogate model }
        \label{fig:monitored_nodes_var}
\end{figure}

 Subsequently, $N_{MC} = 3000$ triplets $\{[E_1^j, E_2^j, E_3^j]\}_{j=1}^{N_{MC}}$ are generated according to the above described log-normal distribution and an MC analysis is performed for both the exact and the surrogate model. Figure \ref{fig:mean_ux_t4} depicts contour plots for the mean value of $u_{x}$ at $t = 4.00 \ sec$ predicted by the two models, while figure \ref{fig:mean_uy_t4} shows the same contour plots for the mean value of $u_{y}$. In addition, figures  \ref{fig:var_ux_t4} and \ref{fig:var_uy_t4} display the variance contours of these displacement fields. Furthermore, figures \ref{fig:monitored_nodes_mean} and \ref{fig:monitored_nodes_var} display a comparison between the two models in the mean value and the variance of the displacements $u_{x}$ and $u_{y}$ of the monitored nodes 1 through 3. Again, the predictions obtained by the proposed CAE-FFNN model are in very close agreement with those computed by the FEM model. The normalized error between the mean solution matrices $\boldsymbol{M}_{FEM}$ and $\boldsymbol{M}_{SUR}$ is equal to 0.62\%, while the same error for the variance matrices $\boldsymbol{V}_{FEM}$ and $\boldsymbol{V}_{SUR}$ is 1.37\%.

\begin{figure}[H]
    \centering
    \includegraphics[width=0.70\textwidth]{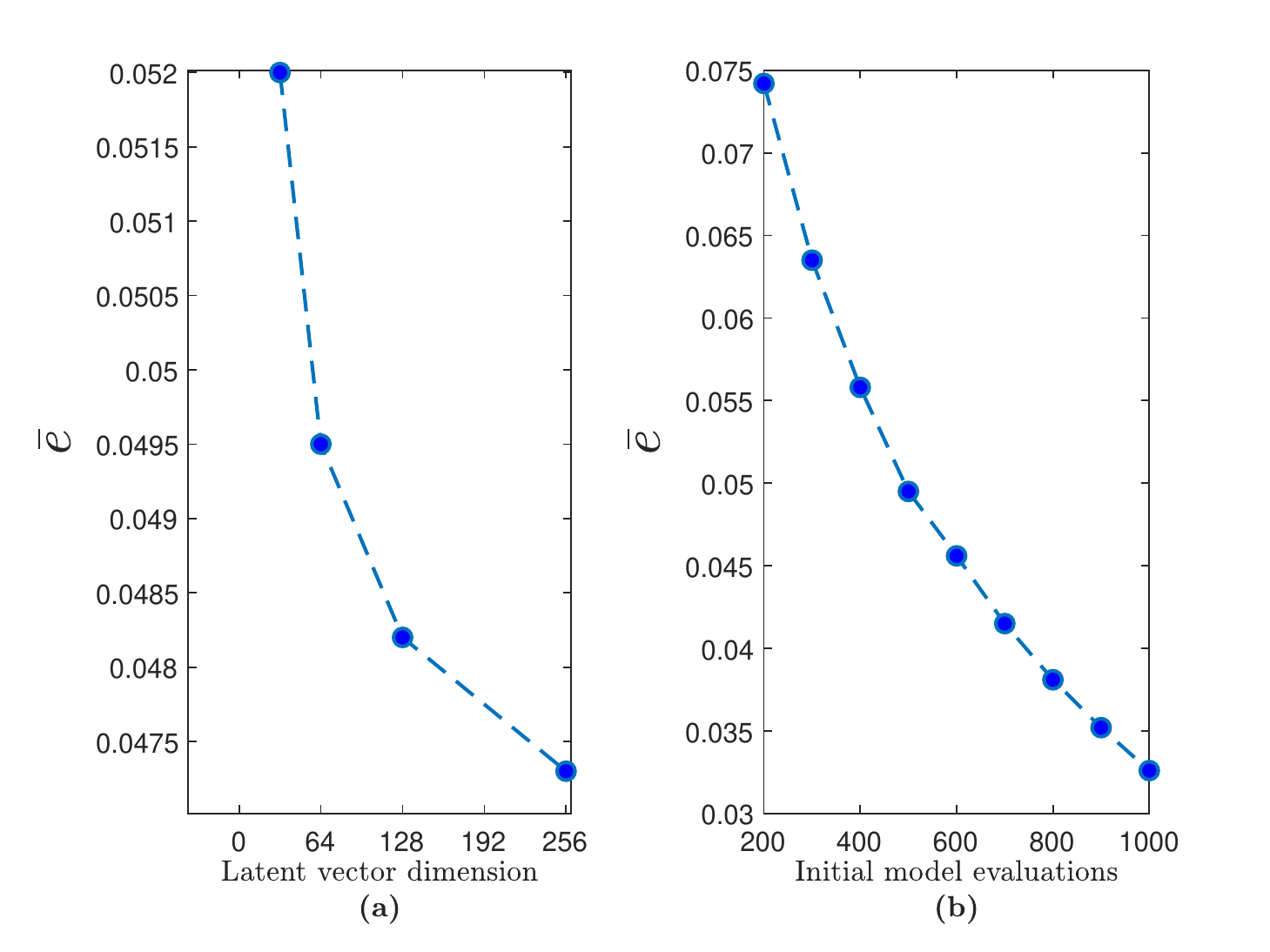}
    \caption{mean error $\Bar{e}$ with respect to (a) the latent vector dimension and (b) the initial data set size}
    \label{fig:error_samples_structural}
\end{figure}

Furthermore, in figure \ref{fig:error_samples_structural} a convergence study with respect to the dimension of the latent vectors and the initial data set size is provided. The average normalized error $\Bar{e}$ of the 3000 MC simulations is given by:
 
\begin{equation}
\Bar{e} = \frac{1}{N_{MC}} \sum_{j=1}^{N_{MC}}\frac{||\boldsymbol{U}_{FEM}^{j} - \boldsymbol{U}_{SUR}^{j}||}{||\boldsymbol{U}_{FEM}^{j}||}
\end{equation}

\noindent with $U_{FEM}^{j}$, $U_{SUR}^{j}$ being the solution matrices of the $j$-th MC simulation obtained by the FEM and the surrogate model, respectively.

\noindent These results indicate that a choice of a higher dimensional latent vector representation leads to improved accuracy, as in the previous example. Furthermore, the average error $\Bar{e}$ decreases as the initial data set size increases and converges close to the value $\Bar{e}_{lim } \approx 0.03$. It is worth mentioning that an optimized set of hyperparameters (latent vector dimension, number of hidden layers, learning rate, etc.) or a different architecture of the CAE and the FFNN  could potentially further reduce the value of $\Bar{e}_{lim}$, but the accuracy achieved for $N = 500$ samples is already deemed adequate for the purposes of this analysis.

\begin{figure}[H]
\centering
\begin{subfigure}[b]{0.48\textwidth}
         \centering
    \includegraphics[width=\textwidth, height = 0.30\textheight]{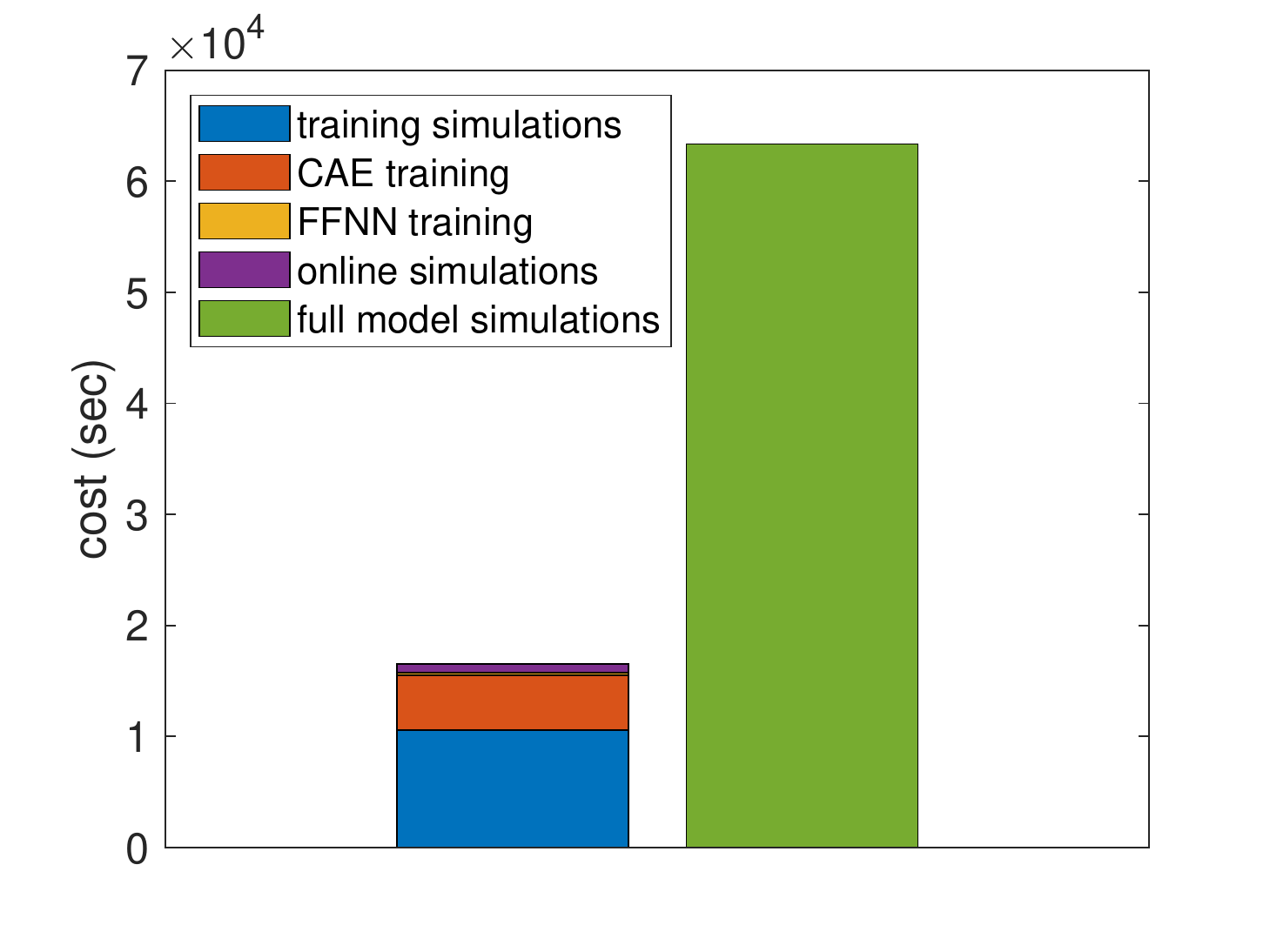}
         \caption{$N_{MC} = 3000$}
             \label{fig:cost_3000}
     \end{subfigure}
     \hfill
     \begin{subfigure}[b]{0.48\textwidth}
         \centering
    \includegraphics[width=\textwidth, height = 0.3\textheight]{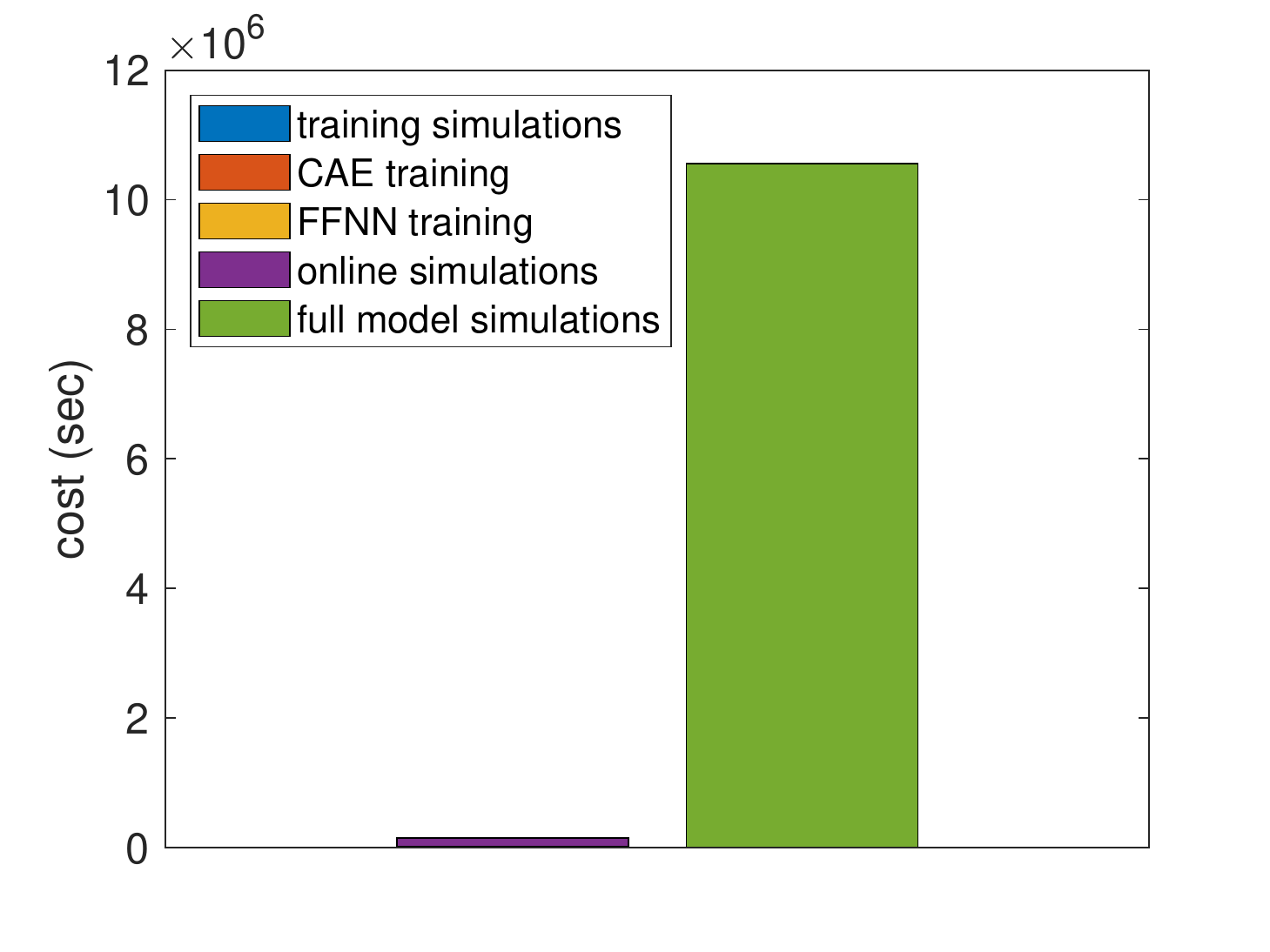}
         \caption{$N_{MC} = 500000$}
     \label{fig:cost_500000}
     \end{subfigure}  
    \caption{Comparison of computational cost between the surrogate and the exact model}
    \label{fig:cost}
\end{figure}

Regarding the computational cost, the results are very promising. Specifically, an MC simulation required an average of 21.12 $sec$ to complete with the exact model, while it only needed 0.26 $sec$ with the surrogate model, which translates to a speed up of $\times81.23$. This drastic cost reduction is the outcome of the 'simulation free' approach of the proposed novel method that eliminates the need of formulating and solving multiple linear systems of equations during the solution procedure of each simulation and is expected to be even greater as the problem's dimensionality increases. The training of the CAE and the FFNN was performed using the GPU version of the Tensorflow framework \cite{tensorflow2015-whitepaper} on an NVIDIA GeForce GTX TITAN X GPU, while the online computations were performed on an Intel\textsuperscript{\textregistered}
 CORE\textsuperscript{TM} -i7 X 980 CPU. Figure \ref{fig:cost_3000} illustrates the computational costs required by the FEM model and the CAE-FFNN model to complete the 3000 MC simulations. This figure also displays the offline cost for training the surrogate and how it was allocated. In particular, the cost for obtaining the 500 initial solutions was 10560 $sec$, the training of the CAE required 4970 $sec$ and the training of the FFNN 211 $sec$. The cost of the 3000 online simulations was only 780 $sec$, which led to a total cost for the surrogate of 16521 $sec$.  On the other hand, the full model MC simulations required 63360 $sec$, almost 4 times that of the surrogate.

\begin{figure}[H]
\centering
\begin{subfigure}[b]{0.48\textwidth}
         \centering
    \includegraphics[width=\textwidth, height = 0.30\textheight]{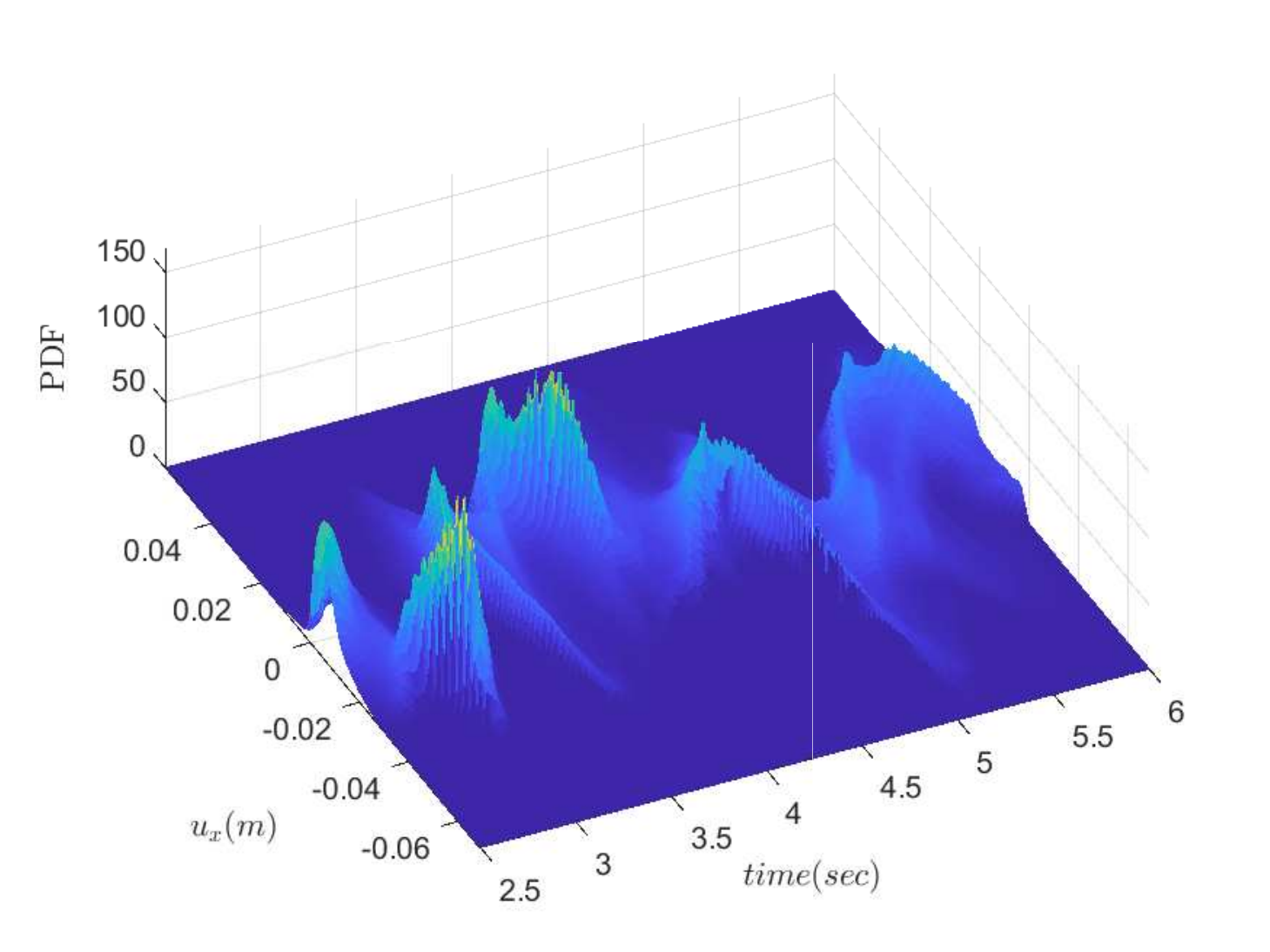}
         \caption{}
     \end{subfigure}
     \hfill
     \begin{subfigure}[b]{0.48\textwidth}
         \centering
    \includegraphics[width=\textwidth, height = 0.30\textheight]{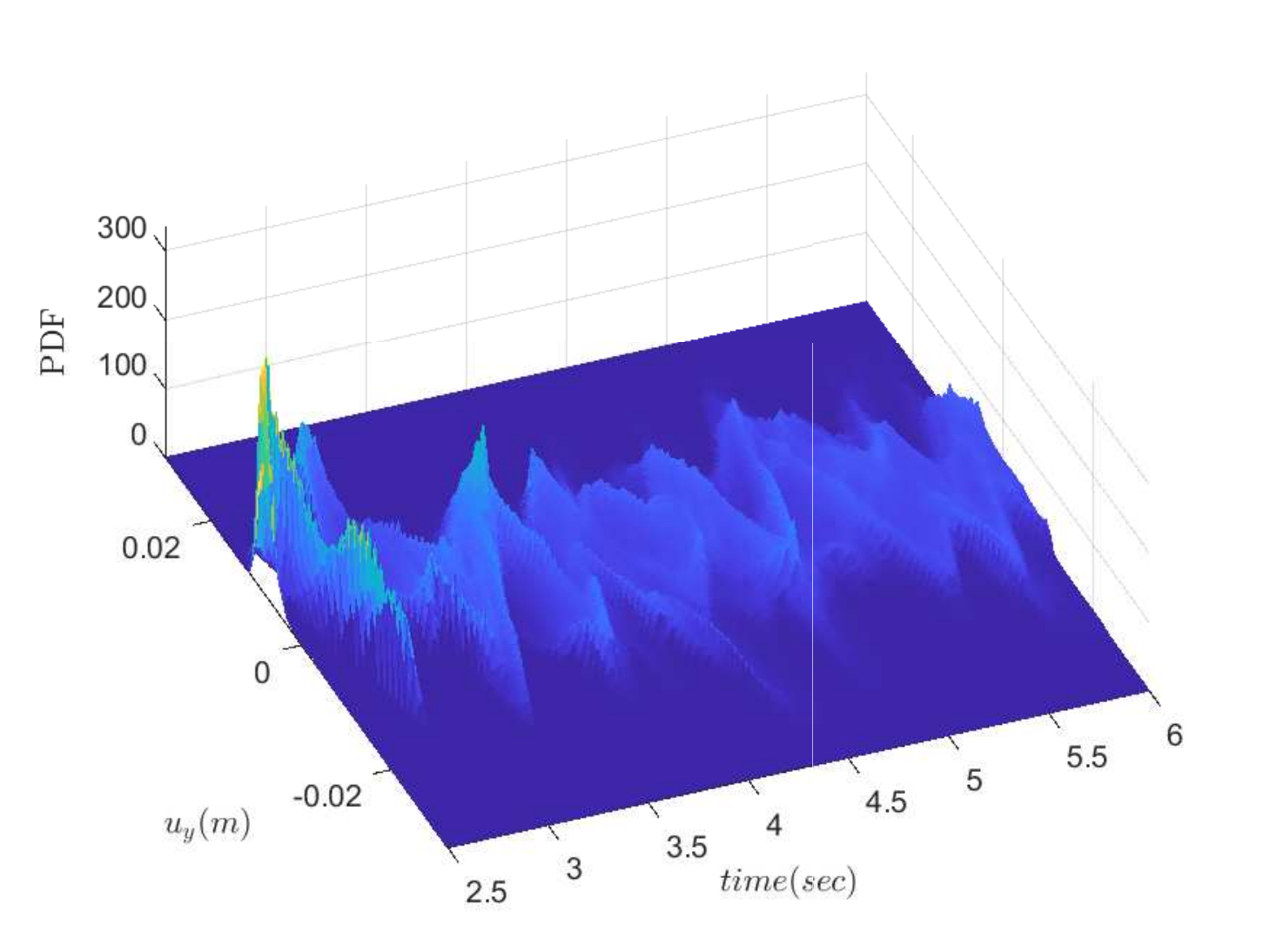}
         \caption{}
     \end{subfigure}  
    \caption{Time evolution of PDF for the displacements (a) $u_{x}$ and (b) $u_{y}$ of the monitored node 3}
    \label{fig:PDF}
\end{figure}

 Finally, the tested surrogate model is utilized to perform  $N_{MC} = 500000$ simulations in order to calculate the time evolution of the probability density function (PDF) of the displacements $u_{x}$ and  $u_{y}$ of the monitored node 3. These results are presented in figure \ref{fig:PDF}. Needless to say, that this analysis would be infeasible without using the proposed surrogate method. In particular, the FEM model would have required approximately 122 $days$ to complete the MC simulation, while the surrogate model required only $40.7$ $hours$, including the offline cost. This remarkable decrease in computational cost is equivalent to a speed up of $\times81.10$. A comparison between the two models is schematically represented in figure \ref{fig:cost_500000}.
 
\section{Conclusions} \label{section5}
This work presented a novel surrogate modeling method based on Convolutional Autoencoders in conjunction with feed-forward Neural Networks. Using a reduced set of system solutions as data set, the CAE provided a low-dimensional representation of this high-dimensional data set through its encoder, as well as the inverse map through its decoder. Then, a FFNN was trained to map points from the problem's parametric space to the encoded solution space and the decoder map was used to reconstruct the system solutions to their original dimension. By composing the FFNN with the decoder, a ’simulation-free’ approach was established to obtain the complete system solutions at a very low cost, rendering this approach ideal for problems requiring multiple model evaluations or ’on-the-fly’ calculations. The method was demonstrated on the solution of time-dependent stochastic PDEs, parametrized by random variables, in the context of the Monte Carlo simulation. The results obtained exhibited high accuracy and remarkable computational gains. Future investigations are focused towards the application of the method to more complex non-linear PDEs.

\section*{Acknowledgments}
\noindent This research has been co-financed by the European Regional Development Fund of the European Union and Greek national funds through the Operational Program Competitiveness, Entrepreneurship and Innovation, under the call RESEARCH-CREATE-INNOVATE (project code: T1EDK-04320).

\newpage
\bibliographystyle{elsarticle-num} 
\bibliography{main}
\end{document}